\documentclass[12pt]{amsart}
\usepackage{subfiles}
\def\insertbibliography{\bibliographystyle{amsplain}\bibliography{thesis_collection}}
\def\biblio{\insertbibliography}

\usepackage{pdfpages}
\usepackage{caption}
\usepackage{subcaption}
\usepackage{amsmath,amsthm,amssymb,amsfonts}
\usepackage{tikz}
\usepackage{float}
\usepackage{todo}
\usepackage{nicefrac}
\usepackage{hyperref}
\usepackage{cleveref}

\theoremstyle{plain}
\newtheorem{thm}{Theorem}[section]
\newtheorem{lemma}[thm]{Lemma}
\newtheorem{prop}[thm]{Proposition}
\newtheorem{cor}[thm]{Corollary}

\newtheorem{Qst}{Question}[section]

\theoremstyle{definition}

\newtheorem{Def}[thm]{Definition}
\newtheorem{Construction}{Construction}
\newtheorem{Example}{Example}
\newtheorem{remark}{Remark}

\usepackage{babel}

\makeatother

\usepackage{babel}

\usepackage[margin=1in]{geometry}

\DeclareRobustCommand{\SkipTocEntry}[5]{}

\newcommand{\mbb}[1]{\mathbb{#1}}
\newcommand{\der}[1]{\frac{\partial}{\partial {#1}}}

\newcommand{\bs}[1]{\boldsymbol{#1}}
\newcommand{\xpoints}{\overline{\bs{\chi}}}
\newcommand{\bphi}{\bs{\varphi}}
\newcommand{\bchi}[1]{\bs{\chi}^{(#1)}}
\newcommand{\conv}{\text{Conv}}
\newcommand{\pair}[2]{\langle {#1}, {#2} \rangle}
\newcommand{\trop}{\text{trop}}
\newcommand{\val}{\text{val}}

\newcommand{\pr}{\operatorname{pr}}
\newcommand{\Sp}{\operatorname{Span}}

\newcommand{\sing}{\text{Sing}}
\begin{document}
	\title{Enumeration of algebraic and tropical singular hypersurfaces}
	
	\author{Uriel Sinichkin}
	\address{School of Mathematical Sciences, Tel Aviv
		University,
		Tel Aviv 69978,
		Israel}
	\email{sinichkin@mail.tau.ac.il}
	
	\thanks
	{\emph{2020 Mathematics Subject Classification}:
		Primary
		14N10. 
		Secondary
		14P05, 
		14T90, 
		14J17, 
		05E14, 
		52B20. 
	}
	
	\date{\today}

\begin{abstract}
	\begin{sloppypar}
We develop a version of Mikhalkin's lattice path algorithm for projective hypersurfaces of arbitrary degree and dimension, which enumerates singular tropical hypersurfaces passing through appropriate configuration of points.
By proving a correspondence theorem combined with the lattice path algorithm, we construct a $ \delta $ dimensional linear space of degree $ d $ real hypersurfaces containing $ \frac{1}{\delta!}(\gamma_nd^n)^{\delta}+O(d^{n\delta-1}) $ hypersurfaces with $ \delta $ real nodes, where $ \gamma_n $ are positive and given by a recursive formula.
This is asymptotically comparable to the number $ \frac{1}{\delta!} \left( (n+1)(d-1)^n \right)^{\delta}+O\left(d^{n(\delta-1)} \right) $ of complex hypersurfaces having $ \delta $ nodes in a $ \delta $ dimensional linear space.
In the case $ \delta=1 $ we give a slightly better leading term.
\end{sloppypar}
\end{abstract}

\def\biblio{}

\keywords{Real singular hypersurface, tropical singular hypersurface, high dimensional tropical geometry.}

\ \vspace{-0.2in}

\maketitle

\vspace{-.25in}

\tableofcontents

\vspace{-.35in}

\newpage

\section{Introduction}

It is classically known that the set of degree $ d $ singular hypersurfaces in $ \mbb{P}^n $ form a hypersurface in the total space of degree $ d $ hypersurfaces $ \mbb{P}^N $ ($ N = \binom{n+d}{n}-1 $).
It is therefore natural to ask what is its degree, which happens to be $ (n+1)(d-1)^n $ (see, for example, \cite[Chapter 9, Example 2.10(a)]{Gelfand_kapranov}).
We can restate this result as an enumerative problem, by picking a generic line in $ \mbb{P}^N $ (also called a \emph{pencil}) and asking how many of the hypersurfaces in this line are singular.
A natural way to pick a pencil of hypersurfaces is to fix a generic collection of $ N-1 $ points $ \overline{\bs{p}}\subseteq \mbb{P}^n $ and consider the hypersurfaces passing through all of them.
Since the condition to pass through a point is linear, this indeed will give a line in $ \mbb{P}^N $.
We thus get, if we work over $ \mbb{C} $, that the number of singular degree $ d $ hypersurfaces through $ N-1 $ points in $ \mbb{CP}^n $ in general position is  $ (n+1)(d-1)^n $, and in particular does not depend on the choice of points (as long as it is generic).
More generally, a generic $ \delta $ dimensional linear space in $ \mbb{CP}^N $ contains
\[ \frac{1}{\delta!}((n+1)(d-1)^n)^\delta + O(d^{n(\delta -1)}) \]
hypersurfaces with $ \delta $ nodes \cite{Kerner_complex_multinodal}.

Over the real numbers the answer to that question is no longer invariant.
That is, there exist different generic sets of points in $ \mbb{RP}^n $ with different number of degree $ d $ real hypersurfaces passing through them and having a real singularity.
One might ask the following natural questions:
\begin{Qst}\label{qst:uninodal}
	What is the maximal possible number of degree $ d $ real hypersurfaces passing through a configuration of $ N-1 $ real points in general position, and having a real singularity.
\end{Qst}
\begin{Qst}\label{qst:multinodal}
	For a positive integer $ \delta $, what is the maximal possible number of degree $ d $ real hypersurfaces passing through a configuration of $ N-\delta $ real points in general position, and having $ \delta $ real singularities.
\end{Qst}

We present our solution for those questions in Sections \ref{sec:real_count} and \ref{sec:multinodal}.
In both cases we have found a lower bound for the maximal number of hypersurfaces by providing an example of a point collection with the prescribed number of singular (respectively, having $ \delta $ singularities) hypersurfaces passing through them. 
We are mainly concerned with the asymptotic of this number as $ d\to \infty $.
To that end we have proved the following two results:
\begin{thm}\label{thm:real_count_formula}
	For all $ n\ge 2 $ and $ d\gg 0 $ there exist a generic real pencil of hypersurfaces of degree $ d $ in $ \mbb{P}^n $ that contains at least $ \alpha_nd^n +O(d^{n-1}) $ singular hypersurfaces, where $ \alpha_n $ satisfy the recurrence relation
	\begin{equation}\label{eq:lemma_real_count_coef_rec}
	\alpha_n = \frac{1}{n}\left( \alpha_{n-1} + \sum_{r=0}^{n-2}(1+\beta_{n-1-r})\alpha_r \right), \ n \ge 2,\quad \alpha_0=1, \ \alpha_1=2,
	\end{equation}
	and $\beta_r $ satisfy the recurrence 
	\begin{equation}\label{eq:real_count_beta_recurrence}
	\beta_r = \frac{1}{r}\left( \beta_{r-1} + \sum_{r_0=1}^{r-2}\beta_{r_0}\beta_{r-r_0-1} \right), \ r\ge 1, \quad \beta_1=1. 
	\end{equation}
\end{thm}
\begin{thm}\label{thm:multinode_count}
	For all $ n\ge 2 $ and $ d\gg 0 $ there exist a configuration $ \overline{\boldsymbol{p}} $ of $ N-\delta = \binom{n+d}{n}-1-\delta $ real points in $ \mbb{P}^n $ s.t. there are at least $ \frac{1}{\delta!}\left( \gamma_n d^n \right)^\delta + O(d^{n\delta-1}) $ real hypersurfaces passing through $ \overline{\boldsymbol{p}} $ and having $ \delta $ real nodes, where $ \gamma_n $ satisfies the recurrence relation \eqref{eq:lemma_real_count_coef_rec}
	\begin{equation*}
	\gamma_n = \frac{1}{n}\left( \gamma_{n-1} + \sum_{r=0}^{n-2}(1+\beta_{n-1-r}\gamma_r) \right), \ n\ge 2, \quad \gamma_0=1, \ \gamma_1=1,
	\end{equation*}
	and $ \beta $ is the sequence defined in Theorem \ref{thm:real_count_formula} and satisfying the relation \eqref{eq:real_count_beta_recurrence}.
	
\end{thm}

One can show that $ \alpha_n \ge 4 $ for all $ n\ge 14 $ (Corollary \ref{cor:at_least_4dn_sols}) and $ \gamma_n\ge \frac{5}{2} $ for all $ n\ge 9 $ (Corollary \ref{cor:gamma_lower_bound}).

Note that in particular, this is asymptotically comparable to the number of complex singular (resp. $ \delta $-nodal) hypersurfaces in a pencil (resp. $ \delta $ dimensional linear space).
When applying the tropical approach presented in the current work to the enumeration of  $ \delta $-nodal complex hypersurfaces, one get the correct leading term $ \frac{1}{\delta !}\left ((n+1)d^n \right )^{\delta} $.
This can be viewed as a confirmation that we have found most of the tropical singular hypersurfaces satisfying the specified conditions.

The introduction of tropical geometry to enumerative geometry is due to Mikhalkin \cite{Mikhalkin2003}, where he provided a way to compute Gromov-Witten invariants of toric surfaces using tropical geometry.
This technique was later generalized by many authors for various enumerative problems.
Most of those results are concerned with the enumeration of curves.
In higher dimensions, tropical geometry was used for the enumeration of real singular surfaces in $ \mbb{P}^3 $ by Markwig, Markwig and Shustin in \cite{Shustin2017_enum} and later for the enumeration of multi-nodal surfaces in $ \mbb{P}^3,\mbb{P}^2\times \mbb{P}^1 $ and $ \mbb{P}^1\times\mbb{P}^1\times\mbb{P}^1 $ by Markwig, Markwig, Shaw and Shustin in \cite{Markwig2019_floor_plan}.
They found arrangements of points in $ \mbb{P}^3 $ with $ \frac{1}{\delta !}(\frac{3}{2}d^3)^{\delta} + O(d^{3\delta-1}) $ surfaces of degree $ d $ with $ \delta $ nodes.
To the best of our knowledge, the higher dimensional problem remained open until now.

The tropical approach amounts to applying a non-Archimedean valuation to the coordinates of the points in $ \mbb{P}^n $.
This transforms the hypersurfaces in question to a piecewise linear objects called \emph{tropical hypersurfaces}.
It is then possible to characterize the singular tropical hypersurfaces which then can be enumerated.
As opposed to the enumeration of singular tropical curves and surfaces, where the classification of empty lattice circuits were used, in higher dimensions no such classification is available (see Remark \ref{rem:no_circ_classification}).
To make the enumeration feasible, one needs to put the points in a special position called \emph{Mikhalkin position} (see Construction \ref{con:Mikhalkin_position}).
Using an approach similar to the floor diagram algorithm presented in \cite{Markwig2019_floor_plan}, we have constructed almost all singular hypersurfaces that pass through collection of points in Mikhalkin position.
Next, we need, for every tropical singular hypersurface, to find all the algebraic singular hypersurfaces tropicalizing to it and satisfying the point conditions.
This is done using a variant of Viro's patchworking technique \cite{Shustin2005, Viro1876}.
We put a special emphasis on writing the explicit equations appearing in the patchworking procedure, which allows us to understand which of their solutions are real.

One should note that there exist additional tools to attack the considered enumerative problem, other than the tropical approach.
Those include integration with respect to Euler characteristic and small deformation of hyperplane arrangement, both described in \cite[Section 2]{Markwig2019_floor_plan}.
Their shortcoming are that they work only for the uninodal problem and produce smaller number of hypersurfaces than our result.
Integration with respect to Euler characteristic works only for even value of $ n $.

The content of the paper is organized as follows. In Section \ref{sec:preliminaries} we define the notions of tropical geometry and recall basic facts about them.
In the same Section we define singular tropical hypersurfaces and characterize them combinatorially.
In Section \ref{sec:initial_data} we explain the choice for the point conditions we consider. 
We then reformulate this choice in terms of lattice paths and prove that a disconnected lattice path can not occur.
Sections \ref{sec:construction_graded_circuits} and \ref{sec:patchworking} are the main sections of this work.
In Section \ref{sec:construction_graded_circuits} we perform most of the combinatorial work, and construct the tropical singularities that can appear on tropical hypersurfaces that pass through points in Mikhalkin position.
In Section \ref{sec:patchworking} we construct the singular algebraic hypersurfaces tropicalizing to each tropical hypersurface constructed in Section \ref{sec:construction_graded_circuits}.
Finally, in Section \ref{sec:results} we combine the results of previous sections to prove Theorems \ref{thm:real_count_formula} and \ref{thm:multinode_count}.
We also show in Section \ref{sec:discriminant_degree} that the leading term of the complex count that we get from our technique is the correct one.

\addtocontents{toc}{\SkipTocEntry}
\subsection*{Acknowledgments}
I would like to thank my advisor Prof. Eugenii Shustin for his support and valuable guidance,
Lev Blechman, Asaf Cohen, Rotem Gluzman and Shay Sadovsky for useful discussions.

The research was supported by Israel Science Foundation grant number 501/18 and by the Bauer-Neuman Chair in Real and Complex Geometry.

\section{Notation and basic facts}\label{sec:preliminaries}

\newcommand{\V}{\mathcal{V}}

\subsection{Tropical hypersurfaces}\label{sec:tropical_geometry}
We recall some basic facts on tropical geometry and set up notation.
For a detailed introduction to the subject see \cite{Mikhalkin_Shustin_Ittenberg, Maclagan2014}.
We assume basic familiarity with the geometry of convex polytopes, as presented for example in \cite{Grunbaum2003, Ziegler2007}.

We denote throughout the paper by 
$$ \Delta_d^{(n)} := \left\{\bs{x}\in \mbb{R}^{n+1}\; | \; \forall i, \bs{x}_i\ge 0, \sum_{i=0}^{n}\bs{x}_i=d \right\} $$
a realization of the $ n $-simplex with side length $ d $ in $ \mbb{R}^{n+1} $ (whenever $ n $ or $ d $ are obvious from the context we might omit them from the notation).
We denote points in $ \Delta $ by boldface Latin letters, such as $ \bs{u,v,w} $ etc.
For a vector $ \bs{v} $ we denote its $ i $'th coordinate as $ \bs{v}_i $.
To avoid confusion, if we have a sequence of vectors, we denote them using upper script in a bracket, for example $ \bs{v}^{(k)} $.
In this vein, we denote by $ \bs{e}^{(k)} $ the unit vector in $ \mbb{R}^{n+1} $ in direction $ k $.
To emphasize the homogeneity, we denote the coordinates of a point $ \bs{v}\in \Delta $ as $ [\bs{v}_0:\bs{v}_1:\ldots :\bs{v}_n] $.

One advantage of working in homogeneous coordinates is that if $ \sum_{\bs{v}\in\Delta}\alpha_{\bs{v}} \bs{v}=\bs{0} $ then $ \sum_{\bs{v}\in \Delta}\alpha_{\bs{v}} =0  $ so every linear dependency of points in $ \Delta $ is affine and every linear combination that lies in $ \Delta $ is an affine combination.
From now on, we will only work with linear dependencies and combinations, while remembering they are in fact affine.
We will denote the linear span of a set $ A $ by $ \Sp(A) $.

Denote by $ \V $ the space $ \mbb{R}^{n+1} $ in which $ \Delta  $ is embedded and by $ \V^* $ its dual (which is also isomorphic to $ \mbb{R}^{n+1} $).
Let $ \bs{1}\in\V^* $ be the functional that sums the coordinates of a point $ \bs{x}\in\V $, i.e. 
\[ \pair{\bs{1}}{\bs{x}} = \sum_{i=0}^n \bs{x}_i. \]
We will denote  $ \nicefrac{\V^*}{\Sp(\bs{1})} $ by $ \mbb{T}^n $ and its elements by boldface Greek letters such as $ \bs{\phi,\psi,\chi} $ etc. 
We view the elements of $ \mbb{T}^n $ as linear functionals on $ \Delta $, up to an additive constant.

Let $ \mbb{K}:=\bigcup_{m\ge 1}\mbb{C}((t^{1/m})) $ be the field of Puiseux series and $ {\mbb{K}_{\mbb{R}}:=\bigcup_{m\ge 1}\mbb{R}((t^{1/m}))} $, its real part.
The field $ \mbb{K} $ is algebraically closed of characteristic 0 and so, by Lefcshetz transfer principle \cite[Theorem 1.13]{Jensen_Lenzing_model_theory}, the solutions of the enumeration problems over it is equal to the solution over $ \mbb{C} $.
Similarly, since $ \mbb{K}_{\mbb{R}} $ is a real closed field, by Tarski's transfer principle \cite[Theorem 1.16]{Jensen_Lenzing_model_theory}, the enumeration problems over it are equivalent to those over $ \mbb{R} $.
Both those fields posses the non-Archimedean valuation $ \val(\sum_{r}a_r t^r) = \min \{ r\in \mbb{Q} \; | \; a_r\ne 0 \} $.

\begin{remark}\label{remark:puiseux_is_function}
	For an element $ a\in\mbb{K} $ with $ \val(a)\ge 0 $, exist a natural number $ m\in \mbb{N} $ s.t. after a parameter change $ t\mapsto t^m $, $ a $ becomes a power series, and thus an analytic function germ in the neighborhood of $ 0 $.
	In this vein, we denote by $ a(0) $ the coefficient of $ t^0 $ for $ a\in\mbb{K} $ with $ \val(a)\ge 0 $.
\end{remark}

For a point $ \bs{z}=[\bs{z}_0:\dots:\bs{z}_n]\in  (\mbb{K}^*)^n\subseteq \mbb{P}^n $ we can compute the negative of its coordinate-wise valuation, $ \trop(\bs{z}):=(-\val(\bs{z}_0),\dots,-\val(\bs{z}_n)) $, called the \emph{tropicalization} of $ \bs{z} $.
Since $ \bs{z} $ is defined up to multiplication by a non-zero constant, $ \trop(\bs{z}) $ is defined up to addition of $ (\alpha,\dots,\alpha) $, i.e. $ \trop(\bs{z})\in \mbb{T}^{n}  $.

We will interested only in homogeneous algebraic and tropical polynomials, so in the sequel a \emph{degree $ d $ polynomial} is a formal sum 
\[ \sum_{\bs{v}\in\Delta_d^{(n)}}a_{\bs{v}}X^{\bs{v}} \]
where $ X^{\bs{v}}=X_1^{\bs{v}_1}\cdot\ldots\cdot X_n^{\bs{v}_n} $, $ X_1,\ldots,X_n $ are variables and the coefficients $ a_{\bs{v}} $ are taken from $ \mbb{K} $ in the algebraic case and from $ \overline{\mbb{R}}:=\mbb{R}\cup\{-\infty\} $ in the tropical case.
The tropical addition and multiplication are denoted by $ \oplus $ and $ \odot $ respectively to distinguish them from the operations on $ \mbb{K} $.

Let $$ f = \sum_{\bs{v}\in\Delta^{(n)}_d\cap\mbb{Z}^{n+1}}a_{\bs{v}}X^{\bs{v}}\in \mbb{K}[X_0,\dots,X_n] $$ be a degree polynomial and $ V(f)\subseteq \mbb{P}^n $ be its vanishing locus.
Applying the tropicalization map pointwise and taking the euclidean closure, we get the \emph{tropicalization} of $ V(f) $:
\[ \trop(V(f)):=\overline{\{\trop(\bs{z}) | \bs{z}\in V(f)\cap (\mbb{K}^*)^n \}}\subseteq \mbb{T}^n. \]
By Kapranov's theorem \cite[Theorem 3.1.3]{Maclagan2014}, $ \trop(V(f)) $ equal to the non-differentiability of the function (defined up to global additive constant)
\begin{align*}
\mbb{T}^n & \to \mbb{R} \\
\bs{\psi} & \mapsto \max\{ \pair{\bs{\psi}}{\bs{v}}-\val(a_{\bs{v}}) | \bs{v}\in\Delta\cap\mbb{Z}^{n+1} \}. 
\end{align*}
It is called \emph{the tropical hypersurface defined by the tropical polynomial $ \bigoplus (-\val(a_{\bs{v}}))X^{\bs{v}} $}.

A central role is played by \emph{the dual subdivision corresponding to a tropical polynomial $ f $}, see \cite{Mikhalkin2003}.
We recall here the main property of the dual subdivision:

\begin{prop}
	Let $ f $ be a tropical polynomial and $ \sigma $ be the dual subdivision corresponding to $ f $.
	The tropical hypersurface $ V(f) $ is a polyhedral complex of dimension $ n-1 $ and the connected components of $ \mbb{T}^n\setminus V(f) $ are the interiors of $ n $ dimensional polytopes.
	Thus considering the closure of those components together with the cells of $ V(f) $ gives a structure of a polyhedral complex on $ \mbb{T}^n $.
	The cells of this complex are in an inclusion reversing bijection to cells of $ \sigma $.
	If $ \delta\in \sigma $, the corresponding cell in $ \mbb{T}^n $ is the set of points $ \bs{\chi} $ s.t. 
	\[ \pair{\bs{\chi}}{\bs{u}}+c_{\bs{u}} = \max \{ \pair{\bs{\chi}}{\bs{v}}+c_{\bs{v}} \; | \; \bs{v}\in\Delta\cap\mbb{Z}^{n+1} \} \]
	for all $ \bs{u}\in\delta\cap\mbb{Z}^{n+1} $.
	If $ \delta\in \sigma $ and $ Q\subseteq \mbb{T}^n $ is the corresponding cell, then $ \dim Q+\dim \delta = n $ and $ Q $ and $ \delta $ are perpendicular to each other, i.e. $ \pair{\bs{\chi}}{\bs{v}} $ remain constant as $ \bs{\chi}\in Q $ and $ \bs{v}\in \delta $ vary.
\end{prop}

\begin{proof}
	See \cite[Proposition 3.11]{Mikhalkin2003}.
\end{proof}

A related concept is that of \emph{Legendre dual function} of a tropical polynomial.
It is the function 
\begin{equation*}
\begin{split}
\nu : \Delta &\to \mbb{R} \\
\bs{x}	&\mapsto \min \{ - \sum_{\bs{v}\in A}\alpha_{\bs{v}} c_{\bs{v}} \; | \; A\subseteq \Delta \cap \mathcal{N}, \sum_{\bs{v}\in A}\alpha_{\bs{v}}=1, \sum_{\bs{v}\in A}\alpha_{\bs{v}}\bs{v}=\bs{x} \}.
\end{split}
\end{equation*}
where $ c_{\bs{v}} $ are the coefficients of the tropical polynomial we are considering.
It is a convex piecewise linear function with linearity domains cells of the dual subdivision.
The Legendre dual of a tropical hypersurface is defined to be the Legendre dual of some of its defining tropical polynomial.
A translation of the tropical hypersurface by a vector $ \bs{\psi} $ amounts to adding the linear functional $ \bs{\psi} $ to $ \nu $.

For a tropical hypersurface $ S $ and an algebraic hypersurface $ \mathcal{S} $, we say that $ \mathcal{S} $ is a \emph{lift} of $ S $, if $ \mathcal{S} $ tropicalize to $ S $.

\subsection{Singular tropical hypersurfaces}\label{sec:tropical_singularity}

Let $ S $ be a tropical hypersurface.
We follow \cite{Dickenstein2012, Markwig2012_trop_surf_sings} and define a point $ \bs{\psi}\in S $ to be \emph{singular} if exist a lift $ \mathcal{S} $ of $ S $ and a singular point $ \bs{q}\in \mathcal{S} $ with $ \trop(\bs{q})=\bs{\psi} $.
We denote the collection of singular tropical hypersurfaces supported on $ \Delta $ by $ \text{Sing}^{\text{tr}}(\Delta) $ and the collection of singular algebraic hypersurfaces with Newton polytope $ \Delta $ by $ \text{Sing}^{\text{alg}}(\Delta) $.
For an arrangement of point $ \xpoints\subseteq \mbb{T}^n $ we denote the collection of singular tropical hypersurfaces passing through them by $ \text{Sing}^{\text{tr}}(\Delta, \xpoints) $.

We begin by stating and proving the following generalization of \cite[Lemma 10]{Markwig2012_trop_surf_sings}, which gives an intrinsic characterization of a tropical singular point.
\begin{lemma}\label{lemma:origin_sing_condition}
	Let $ S $ be a tropical hypersurface, and let $ \nu $ be its Legendre dual.
	For $ c\in \mbb{R} $ denote by $ F_{c} := \{\bs{v}\in \Delta\cap\mbb{Z}^{n+1} \; | \; \nu(\bs{v})=c \} $.
	Then, $ S $ has singularity at $ \bs{0} $ if and only if for every $ c\in \mbb{R} $ exist $ \alpha: F_{c}\to \mbb{Q}^{*} $ such that $ {\sum_{\bs{v}\in F_{c} } \alpha_{\bs{v}} \bs{v} \in \Sp(\bigcup_{c'<c}F_{c'})} $.
	
	In particular, the integral points with lowest value of $ \nu $ are linearly dependent.
\end{lemma}

\begin{proof}
	Suppose that $ \bs{0}\in S $ is a singular tropical point.
	Pick a lift $ \mathcal{S}\in \text{Sing}^{\text{alg}}(\Delta) $ with singularity at a point $ \bs{q} $ tropicalizing to $ \bs{0} $.
	Let $ f(Z)=\sum_{\bs{v}\in\Delta\cap\mbb{Z}^{n+1}} a_{\bs{v}}Z^{\bs{v}} $ be the polynomial defining $ \mathcal{S} $.
	After multiplying the $ i $-th coordinate of $ z $ by $ \bs{q}_i^{-1} $ (which does not change the tropicalization since $ \val(\bs{q}_i)=0 $), we can assume that $ \mathcal{S} $ has singularity at $ {\bs{1}:=[1:\dots:1]} $. 
	Taking the derivative by $ Z_i $ and substituting $ Z_0=\dots=Z_n=1 $, we get that \linebreak $ {\sum_{\bs{v}\in \Delta\cap\mbb{Z}^{n+1}}\bs{v}_i a_{\bs{v}} = 0} $, meaning that $ \sum_{\bs{v}\in \Delta\cap\mbb{Z}^{n+1}}a_{\bs{v}} \bs{v} = \bs{0} $.
	Writing $ a_{\bs{v}} $ as the Puiseux series $ \sum_{r\in \mbb{Q}}\alpha_{\bs{v},r}t^r $ and taking the coefficient of $ t^{c} $ we get that
	\[ \sum_{\bs{v}\in F_{c}} \alpha_{\bs{v},c}\bs{v} = -\sum_{\bs{u}, \nu(\bs{u})<c} \alpha_{\bs{u},c}\bs{u}\in \Sp(\bigcup_{c'< c}F_{c'}) \]
	and $ \alpha_{\bs{v},c}\ne 0 $ for all $ \bs{v}\in F_{c} $ since $ \alpha_{\bs{v},c}t^c $ is the initial term of $ a_{\bs{v}} $ by the definition of $ F_c $.
	The change of field for the coefficients of the linear dependency follow from the usual argument in linear algebra.
	
	In the other direction, let us be given, for all $ \bs{v}\in\Delta\cap\mbb{Z}^{n+1}  $, $ \alpha_{\bs{v}} \in \mbb{Q}^* $ and for all $ \bs{u}\in\Delta\cap\mbb{Z}^{n+1} $ with $ \nu(\bs{u})<\nu(\bs{v}) $, $ \beta_{\bs{u}, \nu(\bs{v})}\in \mbb{Q} $, such that 
	\[ \sum_{\bs{v}\in F_c}\alpha_{\bs{v}} \bs{v} + \sum_{\bs{u},\nu({\bs{u}})<\nu({\bs{v}})} \beta_{\bs{u}, \nu({\bs{v}})} \bs{u} = \bs{0}. \]
	We can define the polynomial $ f(Z)=\sum_{\bs{v}\in\Delta\cap\mbb{Z}^{n+1}}a_{\bs{v}}Z^{\bs{v}} $ with coefficients $$ {a_{\bs{v}}=\alpha_{\bs{v}}t^{\nu({\bs{v}})}+\sum_{\bs{u}, \nu({\bs{u}})>\nu({\bs{v}})}\beta_{\bs{v}, \nu({\bs{u}})}t^{\nu({\bs{u}})}}. $$
	Since $ \val(a_{\bs{v}})=\nu({\bs{v}}) $, the tropicalization of $ V(f) $ is $ S $.
	We will show $ V(f) $ has a singularity at $ \bs{1} $, thus finishing the proof of the lemma.
	Indeed,
	\[ \der{Z_i}f|_{Z=\bs{1}} = \sum_{\bs{v}\in\Delta\cap\mbb{Z}^{n+1}}\bs{v}_i a_{\bs{v}}, \]
	the only possibly non-vanishing terms of the above expressions are terms with $ t^{\nu({\bs{v}})} $ for $ \bs{v}\in\Delta\cap\mbb{Z}^{n+1} $, and they vanish by the choice of $ \alpha_{\bs{v}} $ and $ \beta_{\bs{v},c} $.
\end{proof}

\begin{figure}[h]
	
	\begin{tikzpicture}
	\begin{scope}
	\filldraw [fill=black] (-1, -1) circle (0.03);
	\filldraw [fill=black] (0, 0) circle (0.03);
	\filldraw [fill=black] (0, 0) circle (0.03);
	\filldraw [fill=black] (0, 0) circle (0.03);
	\filldraw [fill=black] (2, 3) circle (0.03);
	\filldraw [fill=black] (0, 0) circle (0.03);
	\filldraw [fill=black] (2, 2) circle (0.03);
	\filldraw [fill=black] (2, 0) circle (0.03);
	\filldraw [fill=black] (3, 2) circle (0.03);
	\filldraw [fill=black] (3, 1) circle (0.03);
	\filldraw [fill=black] (4, 1) circle (0.03);
	
	\draw (-1, -1) -- (-2, -1);
	\draw (-1, -1) -- (-1, -2);
	\draw (-1, -1) -- (0, 0);
	\draw (0, 0) -- (-2, 0);
	\draw (0, 0) -- (2, 2);
	\draw (2, 2) -- (2, 3);
	\draw (2, 3) -- (-2, 3);
	\draw (2, 3) -- (3, 4);
	\draw (0, 0) -- (2, 0);
	\draw (2, 0) -- (2, -2);
	\draw (2, 0) -- (3, 1);
	\draw (3, 1) -- (3, 2);
	\draw (3, 2) -- (2, 2);
	\draw (3, 2) -- (5, 4);
	\draw (3, 1) -- (4, 1);
	\draw (4, 1) -- (4, -2);
	\draw (4, 1) -- (5, 2);
	\end{scope}
	
	\begin{scope}[xshift=250, scale=1.4]
	\draw (0,0) -- (3,0) -- (0,3) -- (0,0);
	\draw (0,1) -- (1,0) -- (1,2) -- (0,2) -- (2,0) -- (2,1) -- (1,1);
	
	\filldraw [fill=black] (0, 0) circle (0.05) node [below] {\footnotesize $2$};
	\filldraw [fill=black] (0, 1) circle (0.05) node [anchor=east] {\footnotesize $3$};
	\filldraw [fill=black] (0, 2) circle (0.05) node [anchor=east] {\footnotesize $3$};
	\filldraw [fill=black] (0, 3) circle (0.05) node [anchor=east] {\footnotesize $0$};
	\filldraw [fill=black] (1, 0) circle (0.05) node [below] {\footnotesize $3$};
	\filldraw [fill=black] (1, 1) circle (0.05) node [anchor=east] {\footnotesize $3$};
	\filldraw [fill=black] (1, 2) circle (0.05) node [above] {\footnotesize $1$};
	\filldraw [fill=black] (2, 0) circle (0.05) node [below] {\footnotesize $1$};
	\filldraw [fill=black] (3, 0) circle (0.05) node [below] {\footnotesize $-3$};
	
	\draw [->] (3,2) -- (3,2.5) node [above] {$ x_1 $ };
	\draw [->] (3,2) -- (3.5,2) node [anchor=west] {$ x_2 $ };
	\draw [->] (3,2) -- (2.65,1.65) node [anchor=east] {$ x_0 $ };
	\end{scope}
	\end{tikzpicture}
	
	\caption{The tropical hypersurface defined by the tropical polynomial \eqref{eq:tropical_polynomial_example} and the dual subdivision corresponding to it. 
		The vertical direction is $ x_1 $, the horizontal direction is $ x_2 $ and $ x_0 $ is the distance from the diagonal edge of $ \Delta $.
		The numbers near the points of $ \Delta $ indicate the coefficients of the tropical polynomial \eqref{eq:tropical_polynomial_example}.}\label{fig:tropical_hypersurface_example}
\end{figure}
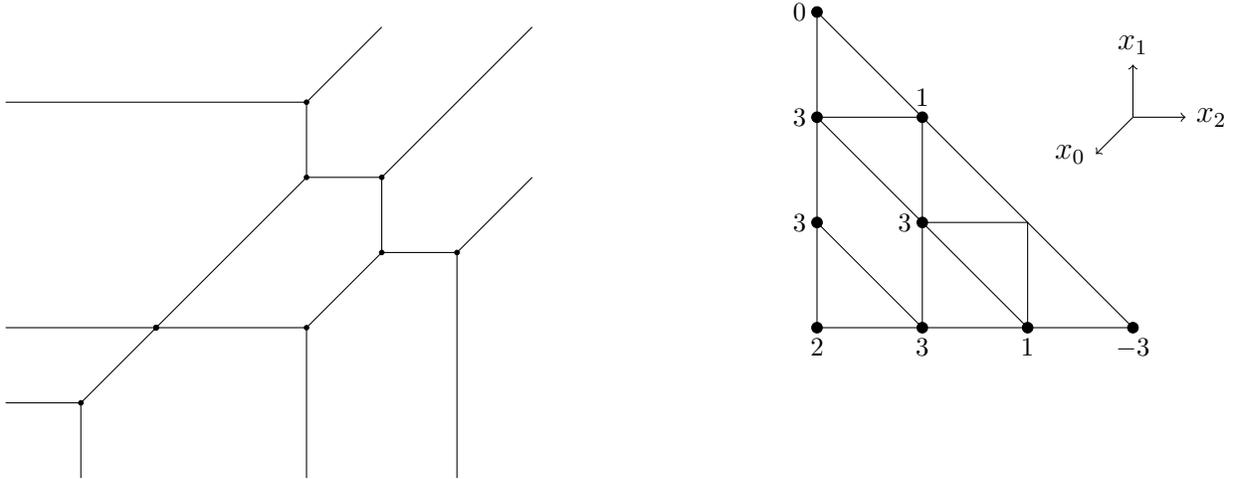

\begin{Example}
The tropical hypersurface defined by
\begin{equation}\label{eq:tropical_polynomial_example}
\begin{split}
f = 2\odot X_0^3 \oplus 3\odot X_0^2 X_1 \oplus 3\odot X_0 X_1^2 \oplus X_1^3 \oplus 3 \odot X_0^2 X_2 & \oplus \\ 
\oplus 3\odot X_0 X_1 X_2 \oplus 1\odot X_1^2 X_2 \oplus 1\odot X_0 X_2^2 \oplus X_1 X_2^2 \oplus (-3)\odot X_2^3 ,
\end{split}
\end{equation}
see Figure \ref{fig:tropical_hypersurface_example}, has a singularity at $ \bs{0} \in \mbb{T}^2 $.
Note the parallelogram in the dual subdivision.
\end{Example}

We get an immediate corollary of Lemma \ref{lemma:origin_sing_condition}.
\begin{cor}\label{cor:nu_image_size}
	Let $ S $ be a tropical hypersurface, with Legendre dual $ \nu $.
	\begin{enumerate}
		\item $ |\nu\left( \Delta\cap\mbb{Z}^{n+1} \right) | \le |\Delta\cap\mbb{Z}^{n+1}|-n-1 $.
		\item If there is equality in $ (1) $ then for every $ c\in im(\nu) $, the coefficients $ \alpha $ as in Lemma \ref{lemma:origin_sing_condition} are unique up to multiplication by a non zero constant.
	\end{enumerate}
\end{cor}

\begin{proof}
	\begin{enumerate}
		\item We use the notations of Lemma \ref{lemma:origin_sing_condition}.
		For $ c\in \left( \Delta\cap\mbb{Z}^{n+1} \right) $ define the linear transformation
		\begin{align*}
		T_c: \mbb{Q}^{F_c} & \mapsto \nicefrac{\Sp\left( \bigcup_{c'\le c}F_{c'} \right) }{\Sp(\bigcup_{c'< c}F_{c'})} \\
		\left( \alpha_{\bs{v}} \right)_{\bs{v}\in F_c} & \mapsto \left[\sum_{\bs{v}\in F_c}\alpha_{\bs{v}}\bs{v}\right] 
		\end{align*}
		where we denote by $ \mbb{Q}^{F_c} $ the free $ \mbb{Q} $-vector space generated by $ F_c $ and the square brackets in the image indicate equivalence class in the quotient by $ \Sp(\bigcup_{c'< c}F_{c'}) $.
		By definition, $ T_c $ is surjective, and by Lemma \ref{lemma:origin_sing_condition}, $ \dim \ker T_c\ge 1 $.
		So, we get
		\begin{equation}\label{eq:nu_image_size_dim_calcs}
		\begin{split}
		\dim \Sp\left( \bigcup_{c'\le c} F_{c'} \right) - \dim \Sp\left( \bigcup_{c'< c} F_{c'} \right) = & \dim \nicefrac{\Sp\left( \bigcup_{c'\le c}F_{c'} \right) }{\Sp(\bigcup_{c'< c}F_{c'})} = \\
		=  \dim \mbb{Q}^{F_c} - \dim \ker T_c \le & |F_c| - 1.
		\end{split}
		\end{equation}
		Summing this over all $ c\in \nu(\Delta\cap\mbb{Z}^{n+1}) $ we get that
		\[ n+1 \le |\Delta\cap\mbb{Z}^{n+1}| - |\nu(\Delta\cap\mbb{Z}^{n+1})| \]
		which yields the desired inequality.
		
		\item
		The given implies there is equality in \eqref{eq:nu_image_size_dim_calcs} for all $ c\in \nu(\Delta\cap\mbb{Z}^{n+1}) $ which means that $ \dim \ker T_c=1 $ for all $ c\in \nu(\Delta\cap\mbb{Z}^{n+1}) $ which is exactly what is needed to be shown.
	\end{enumerate}
\end{proof}

Let us be given a sequence of disjoint sets $ F_1,\dots,F_k\subseteq \Delta\cap\mbb{Z}^{n+1} $ that cover $ \Delta\cap\mbb{Z}^{n+1} $.
The collection of $ \nu\in\mbb{R}^{\Delta\cap\mbb{Z}^{n+1}} $ such that exist a sequence of real numbers $ c_1<\dots<c_k $ with $ \{ \bs{v}\in\Delta\cap\mbb{Z}^{n+1} \; | \; \nu(\bs{v})=c_i \} = F_i $ is either empty or is a relative interior of a polyhedral cone of dimension $ k $.
Since the Legendre dual of a tropical hypersurface is defined up to an additive constant, the collection of tropical hypersurfaces with Legendre dual satisfying the above condition (for a given sequence $ F_1,\dots,F_k $) is either of dimension $ k-1 $ or empty.

We do not restrict our attention to tropical hypersurfaces with singularity at $ \bs{0} $.
As mentioned above, translation of a tropical hypersurface correspond to addition of linear functional to the Legendre dual $ \nu $.
Let $ \bs{\psi} $ be a singular point of a tropical hypersurface $ S $.
If after translating $ \bs{\psi} $ to $ \bs{0} $ the Legendre dual of $ S $ satisfies condition $ (2) $ in Corollary \ref{cor:nu_image_size}, we say $ \bs{\psi} $ is \emph{of maximal-dimensional geometric type}.

Given sequence $ F_1,\dots,F_k $ as above, the space of tropical hypersurfaces that have some translate such that the Legendre dual of it satisfies $ \{ \bs{v}\in\Delta\cap\mbb{Z}^{n+1} \; | \; \nu(\bs{v})=c_i \} = F_i $ for some $ c_1<\dots<c_k $ is the relative interior of a cone of dimension at most $ n+k-1 $.
In particular, the dimension of tropical hypersurfaces that have a singularity of maximal-dimensional geometric type is at most $ |\Delta\cap\mbb{Z}^{n+1}|-2 $ and hence $ \dim \text{Sing}^{\text{tr}}(\Delta)\le|\Delta\cap\mbb{Z}^{n+1}|-2 $.

\section{Choice of initial data}\label{sec:initial_data}

In this section we describe the choice of initial data, namely the tropical point conditions $ \xpoints $.
We show that our choice gives rise to a lattice path similarly to \cite{Shustin2017_enum, Mikhalkin2003_lattice_path} and prove that a disconnected lattice path can not occur in our situation.

\subsection{Points in Mikhalkin position}

Recall that we defined in Section \ref{sec:tropical_geometry} $$ \Delta_d^{(n)} := \left \{ \bs{x}\in \mbb{R}^{n+1} \; | \; (\forall 0\le i\le n)\bs{x}_i\ge 0, \sum_{i=0}^n \bs{x}_i=d \right  \} $$ to be a realization of the simplex, and that we omit the upper and lower indices from the notation whenever the dimension and degree are clear from the context.
To apply the lattice path algorithm, as in \cite{Shustin2017_enum, Mikhalkin2003_lattice_path}, we will need to pick the points $ \overline{\bs{\chi}} $ in special position.
For that, define the following order on $ \Delta_d^{(n)}\cap\mbb{Z}^{n+1} $.

\begin{Def}\label{def:prec_order}
	Define the order $ \prec $ on $ \Delta\cap\mbb{Z}^{n+1} $ to be the reversed lexicographic order, i.e. 
	\begin{equation}\label{eq:def_prec_order}
	\bs{x}\prec \bs{y} \iff \left(\exists i\right)\left (\bs{x}_i < \bs{y}_i \land \left( \forall j>i \right) \bs{x}_j=\bs{y}_j \right ).
	\end{equation}
	
	Clearly $ \prec $ is defined by a linear functional, i.e. exist $ \bphi: \Delta\cap\mbb{Z}^{n+1}\to \mbb{Q} $ s.t. $ \pair{\bphi}{\bs{x}}=\pair{\bphi}{\bs{y}} $ if and only if $ \bs{x}=\bs{y} $ and $ \pair{\bphi}{\bs{x}}<\pair{\bphi}{\bs{y}} $ if and only if $ \bs{x}\prec\bs{y} $.
	We will fix such a functional throughout the paper.
\end{Def}

\begin{Construction}\label{con:Mikhalkin_position}
	Recall that we have identified, in Section \ref{sec:tropical_geometry}, points in $ \mbb{T}^n $ and linear functionals on $ \Delta $ up to additive constant.
	We thus can pick the point conditions to be $ \bs{\chi}^{(i)} = M_i \bphi $ for all $ i=1,\dots,N-1 $, where $ N:=|\Delta\cap\mbb{Z}^{n+1}|-1 $, $ \bphi $ is as in Definition $ \ref{def:prec_order} $ and $ M_1\ll M_2\ll \dots\ll M_{N-1} $ are positive rationals.
	
	The symbol $ \ll $ should be understood as follows. 
	If $ \sum_{i=1}^k \alpha_i M_i $ is a linear combination with $ 1\le k\le N-1 $ and $ \alpha_k> 0 $ we will require 
	\[ \sum_{i=1}^k\alpha_i M_i >0. \]
	Since for every $ d,n $ we will consider only finitely many such linear combinations, this choice is always possible.
	One can actually find a finite set of linear combinations that are sufficient to consider by examining the proofs of Lemma \ref{lemma:disconnected_path_illuminated}, Proposition \ref{prop:no_disconnected} and Lemma \ref{lemma:mikhalkin_pos_imply_cond}.
	
	We say that the points $ \overline{\bs{\chi}} $ are in \emph{Mikhalkin position}.
\end{Construction}

\begin{remark}\label{rem:mikhalkin_pos_is_generic}
	Note that we make the choice of $ N-1 $ parameters $ M_1,\dots,M_{N-1} $ in Construction \ref{con:Mikhalkin_position}.
	By the discussion following Corollary \ref{cor:nu_image_size}, this means that the singularities that can occur are of maximal-dimensional geometric type so, by Corollary \ref{cor:nu_image_size}, the dual subdivision of any singular surface $ S\in \text{Sing}^{\text{tr}}(\Delta, \xpoints) $ contains a unique cell whose integral points form a circuit.
	This also implies that $ \text{Sing}^{\text{tr}}(\Delta, \xpoints) $ is finite and that for every element $ S $ of $ \text{Sing}^{\text{tr}}(\Delta, \overline{\bs{\chi}}) $, each $ \bs{\chi}^{(i)} $ lies in the interior of an $ n-1 $ dimensional face $ F_i $ of $ S $, with $ F_i\ne F_j $ whenever $ i\ne j $.
\end{remark}

\subsection{Reformulation in terms of lattice path}

Let $ \overline{\bs{\chi}}\subseteq \mbb{T}^{n} $ be a collection of points in Mikhalkin position as in Construction \ref{con:Mikhalkin_position}, and let $ S\in \text{Sing}^{\text{tr}}(\Delta, \overline{\bs{\chi}}) $ be a singular tropical hypersurface passing through them.
By Remark \ref{rem:mikhalkin_pos_is_generic}, each $ \bs{\chi}^{(i)} $ is in the interior of a distinct $ (n-1) $-face $ F_i $ of $ S $.
If we denote by $ \nu:\Delta\to\mbb{R} $ the piecewise linear convex function dual to $ S $, the face $ F_i $ is dual to an edge $ E_i $ of the subdivision induced by $ \nu $.
We denote by $ P(S, \overline{\bs{\chi}}) $ the collection of all those edges and call it \emph{the lattice path corresponding to  $ (S, \overline{\bs{\chi}}) $}.

For a subset $ A\subseteq \Delta\cap\mbb{Z}^{n+1} $ consisting of the $ m $ points $ \bs{v}^{(1)}\prec \dots \prec \bs{v}^{(m)} $, denote by $ P(A) $ the \emph{complete lattice path supported on $ A $}
\[ P(A):=\left \{[\bs{v}^{(i)},\bs{v}^{(i+1)}] \; | \; i=1,\dots, m-1 \right  \}. \]
With this notation we can state the following analog of \cite[Lemma 3.2]{Shustin2017_enum} where it is stated for $ n=3 $.
We will not include a proof since the proof that appear in \cite{Shustin2017_enum} works in arbitrary dimension.

\begin{lemma}\label{lemma:lattice_path}
	For a singular tropical hypersurface $ S $ passing through $ \overline{\bs{\chi}} $, the lattice path $ P(S,\overline{\bs{\chi}}) $ defined above satisfies one of the following:
	\begin{enumerate}
		\item Either $ P(S,\overline{\bs{\chi}}) = P(\Delta\cap\mbb{Z}^{n+1}\setminus \{\bs{w}\}) $ for some $ \bs{w}\in \Delta\cap\mbb{Z}^{n+1} $, 
		
		\item or $ P(S, \overline{\bs{\chi}}) = P(\Delta\cap \mbb{Z}^{n+1})\setminus \{ [\bs{w}^{\text{pr}}, \bs{w}] \} $ for some $ \bs{w}^{\text{pr}},\bs{w}\in\Delta\cap\mbb{Z}^{n+1} $, where $ \bs{w}^{\text{pr}} $ is the first (w.r.t. $ \prec $) point of $ \Delta\cap\mbb{Z}^{n+1} $ preceding $ \bs{w} $\footnote{Here and in the rest of the Section, we denote by $ \bs{w}^{\text{pr}} $ the last point in $ \Delta\cap\mbb{Z}^{n+1} $ preceding $ \bs{w} $ and by $ \bs{w}^{\text{suc}} $ the first point in $ \Delta\cap\mbb{Z}^{n+1} $ succeeding $ \bs{w} $.}.
	\end{enumerate}
	We will call the first case {connected lattice path $ \Gamma_{\bs{w}} $} and the second case {disconnected lattice path missing the edge $ [\bs{w}^{\text{pr}},\bs{w}] $}.
\end{lemma}

\begin{lemma}\label{lemma:all_lattice_path_complex}
	Let $ \lambda\subseteq \mbb{R}^{\Delta\cap\mbb{Z}^{n+1}}/ (1,\dots,1) $ be the set of tropical polynomials that define a hypersurface which passes through $ \xpoints $.
	Then $ \lambda $ is a one dimensional polyhedral complex, its cells correspond to different lattice paths.
	The unbounded cells of $ \lambda $ are parallel to the axes and correspond to connected lattice paths (including the last two paths in the proof of Lemma \ref{lemma:lattice_path}).
	The bounded cells of $ \lambda $ correspond to disconnected lattice paths.
\end{lemma}

\begin{proof}
	Let $ \bs{w}\in\Delta\cap\mbb{Z}^{n+1} $, $ \Gamma_{\bs{w}}:=P(\Delta\cap\mbb{Z}^{n+1}\setminus \{\bs{w} \}) $ be the corresponding connected lattice path.
	Take a tropical polynomial $ F\in\lambda $ that defines a tropical hypersurface $ S $ with a lattice path $ P(S,\xpoints)=\Gamma_{\bs{w}} $ and write $ F=\bigoplus_{\bs{v}\in\Delta\cap\mbb{Z}^{n+1}}  c_{\bs{v}} \odot X^{\bs{v}} $.
	If $ [\bs{u},\bs{v}] \in \Gamma_{\bs{w}} $ is dual to the $ (n-1) $ dimensional cell of $ S $ containing $ \bchi{i} $, then 
	\[ \pair{\bchi{i}}{\bs{u}} + c_{\bs{u}} = \pair{\bchi{i}}{\bs{v}} + c_{\bs{v}} \]
	so $ c_{\bs{v}}=c_{\bs{u}}-M_i\pair{\bphi}{\bs{v}-\bs{u}} $.
	Since $ \Gamma_{\bs{w}} $ is connected, this defines the coefficients $ c_{\bs{v}} $ uniquely for all $ \bs{v}\ne \bs{w} $, up to additive constant, which is modded out in $ \mbb{R}^{\Delta\cap\mbb{Z}^{n+1}}/ (1,\dots,1) $.
	
	The only thing to make sure of for $ F $ to be in $ \lambda $, is that 
	\begin{equation}\label{eq:connected_collection_ineq}
	\pair{\bchi{i}}{\bs{w}}+c_{\bs{w}} \le \pair{\bchi{i}}{\bs{u}} + c_{\bs{u}}
	\end{equation} whenever $ [\bs{u}, \bs{v}] $ is an edge of $ \Gamma_{\bs{w}} $ that correspond to a cell passing through $ \bchi{i} $.
	We assert that the collection \eqref{eq:connected_collection_ineq} of inequalities is equivalent to 
	\begin{equation}\label{eq:connected_path_ineq}
	c_{\bs{w}} \le M_k\pair{\bphi}{\bs{w}^{\text{pr}}-\bs{w}^{\text{suc}}}+c_{\bs{w}^{\text{pr}}}
	\end{equation}
	where $ \bchi{k} $ is contained in the cell dual to $ [\bs{w}^{\text{pr}},\bs{w}^{\text{suc}}] $.
	Indeed, \eqref{eq:connected_path_ineq} is one of the inequalities appearing in \eqref{eq:connected_collection_ineq} so it must hold.
	Conversely, suppose that \eqref{eq:connected_path_ineq} holds and let $ [\bs{u},\bs{v}] \in \Gamma_{\bs{w}} $.
	If $ \bs{w}\prec \bs{u} $ then the coefficient of $ M_i $ in the right-hand side of \eqref{eq:connected_collection_ineq} is greater than its coefficient on the left-hand side, so the inequality holds.
	If, on the other hand, $ \bs{v}\prec\bs{w} $ then $ M_k $ is the last of the $ M_j $'s appearing in \eqref{eq:connected_collection_ineq} and it appears only on the left-hand side, with negative coefficient.
	
	Finally, the cell of $ \lambda $ that corresponds to disconnected path $ P(\Delta\cap\mbb{Z}^{n+1})\setminus \{[\bs{w}^{\text{pr}},\bs{w}]\} $ is homeomorphic to the segment $ [\bchi{k}, \bchi{k+1}] $ for some $ 1\le k< N-2 $, and thus bounded.
\end{proof}

\begin{remark}
	It is easy to show that the polyhedral complex $ \lambda $ in Lemma \ref{lemma:all_lattice_path_complex} is in fact a \emph{tropical line}.
	We will not need neither this fact nor the definition of tropical line, so we omit them from this work.
	We refer the curious reader to \cite{Maclagan2014} for the definition of tropical line.
\end{remark}

\subsection{Impossibility of singular tropical hypersurfaces with disconnected path}
In this subsection, we will show that disconnected path could not happen for singular tropical hypersurfaces passing through $ \xpoints $.

\begin{Def}\label{def:illuminated}
	Let $ \V $ be a real vector space, $ \bs{w}\in \V $ and $ X\subseteq \V $ convex.
	Denote the \emph{$ \bs{w} $-illuminated set of $ X $} \[ \partial_{\bs{w}} X = \left \{\bs{x}\in X\; | \; \left[\bs{x}, \bs{w}\right]\cap X= \{\bs{x}\}  \right \} \]
	where $ \left[\bs{a},\bs{b}\right] $ denotes the straight line segment connecting $ \bs{a} $ and $ \bs{b} $. 
\end{Def}

\begin{lemma}\label{lemma:illuminated_calc}
	Let $ \bs{w}\in\Delta\cap\mbb{Z}^{n+1} $.
	Set $ \Omega:=\conv \{ \bs{v}\in\Delta\cap\mbb{Z}^{n+1} \; | \; \bs{v}\prec \bs{w} \} $, then $ \partial_{\bs{w}}\Omega $ is a face of $ \Omega $.
\end{lemma}

\begin{proof}
	For $ k=1,\dots,n $, denote 
	\[ \Omega^{(k)} = \conv \{ \bs{v}\in\Delta\cap\mbb{Z}^{n+1} \; | \; \bs{v}\prec\bs{w}, (\forall i>k) \bs{v}_i=\bs{w}_i \}. \]
	We will show, by induction on $ k $, that $ \partial_{\bs{w}}\Omega^{(k)} $ is a face of $ \Omega^{(k)} $ of dimension at most $ k-1 $.
	
	If $ \bs{w}_i=0 $ for all $ 0<i\le k $, then $ \Omega^{(k)} $ is empty, so there is nothing to prove.
	If $ \bs{w}_i=0 $ for all $ 0<i<k $ and $ \bs{w}_k>0 $, $ \partial_{\bs{w}}\Omega^{(k)} $ is contained in the supporting hyperplane $ \bs{x}_{k}=\bs{w}_k-1 $.
	If $ \bs{w}_k=0 $, then $ \Omega^{(k)}=\Omega^{(k-1)} $, so the assertion follows by induction.
	
	Now, suppose that $ \bs{w}_k>0 $ and there exist $ 0<i<k $ s.t. $ \bs{w}_i>0 $.
	Let $ \ell<k $ be the maximal index with $ \bs{w}_{\ell}>0 $. 
	By induction, $ \partial_{\bs{w}}\Omega^{(\ell)} $ is a polytope of dimension at most $ \ell-1 $, so 
	\begin{equation*}
	F:=\conv\left( \partial_{\bs{w}}\Omega^{(\ell)} \cup \{ \bs{v}\in\Delta \; | \; (\forall j<\ell)\bs{v}_j=0, \bs{v}_k=\bs{w}_k-1, (\forall i>k)\bs{v}_i=\bs{w}_i \} \right)
	\end{equation*}
	is a polytope of dimension at most $ k-1 $, we assert that $ F=\partial_{\bs{w}}\Omega^{(k)} $.
	Indeed, by considering the $ k $'th coordinate, it is immediate that for every point $ \bs{v} $ of $ F $ the segment $ [\bs{v},\bs{w}] $ intersect $ \Omega^{(k)} $ in $ \bs{v} $ only.
	For the other direction, it is enough to show that all the points of $ \Omega^{(k)}\cap\mbb{Z}^{n+1}\setminus F $  are not contained in $ \partial_{\bs{w}}\Omega^{(k)} $.
	Let $ \bs{v}\in \Omega^{(k)}\cap\mbb{Z}^{n+1}\setminus F $.
	If $ \bs{v}_k<\bs{w}_k-1 $ then $ [\bs{v},\bs{w}] $ intersects $ \Omega^{(k)} $ in a point $ \bs{u} $ with $ \bs{u}_k=\bs{w}_k-1 $.
	Otherwise, by construction of $ F $, $ \bs{v}_j>0 $ for some $ j<\ell $.
	Then $ \bs{v}+\bs{e}^{(\ell)}-\bs{e}^{(j)}\in \Omega^{(k)} $ and $ \bs{w}+\bs{e}^{(j)}-\bs{e}^{(\ell)}\in \Omega^{(k)} $, so $$ \frac{1}{2}(\bs{v}+\bs{w}) = \frac{1}{2}(\bs{v}+\bs{e}^{(\ell)}-\bs{e}^{(j)} + \bs{w}+\bs{e}^{(j)}-\bs{e}^{(\ell)}) \in \Omega^{(k)} $$
	what was needed to be shown.
\end{proof}

\begin{lemma}\label{lemma:disconnected_path_illuminated}
	\begin{sloppypar}
	Let $ S $ be a singular hypersurface passing through $ \xpoints $ with $ {P(S,\xpoints)=P(\Delta\cap\mbb{Z}^{n+1})\setminus \{[\bs{w}^{\text{pr}},\bs{w}]\}} $ and let $ C_S $ be the unique cell of the dual subdivision that contains a circuit of the subdivision corresponding to $ S $ (see Remark \ref{rem:mikhalkin_pos_is_generic}).
	Denote by $ \Omega:=\conv(\{ \bs{v}\in\Delta\cap\mbb{Z}^{n+1} \; | \; \bs{v}\prec \bs{w} \}) $.
	Then $ C_S\cap\Omega\subseteq \partial_{\bs{w}}\Omega $ and $ C_S\cap \mbb{Z}^{n+1}\setminus \Omega = \{\bs{w}, \bs{w}^{\text{suc}} \} $.
	\end{sloppypar}
\end{lemma}

\begin{proof}
	Let $ \nu:\Delta\to\mbb{R} $ be the Legendre dual piecewise linear convex function, so the defining tropical polynomial of $ S $ is 
	\[ F_S := \bigoplus_{\bs{v}\in\Delta\cap\mbb{Z}^{n+1}}  (- \nu(\bs{v}))\odot X^{\bs{v}} . \]
	If an edge $ [\bs{u}, \bs{v}] $ is dual to to the $ (n-1) $ cell of $ S $ that passes through $ \bs{\chi}_i $ then 
	\begin{equation*} 
	\langle \bs{\chi}_i, \bs{u} \rangle - \nu(\bs{u}) = \langle \bs{\chi}_i, \bs{v} \rangle - \nu (\bs{v})
	\end{equation*}
	meaning that,
	\begin{equation*}
	\nu(\bs{v}) = \nu(\bs{u}) + M_i\cdot \langle \bphi, \bs{v}-\bs{u} \rangle
	\end{equation*}
	and similarly, 
	\begin{equation*}
	\nu(\bs{w}) = \nu(\bs{w}^{\text{pr}}) + M_{\bs{\eta}}\cdot \pair{ \bphi}{ \bs{w}-\bs{w}^{\text{pr}}}.
	\end{equation*}
	where $ \bs{\eta} = M_{\bs{\eta}}\cdot\bphi $ is the additional point on $ L $ as in the proof of Lemma \ref{lemma:lattice_path}. We thus conclude, since $ M_k<M_{\bs{\eta}}<M_{k+1} $, that for all $ \bs{u},\bs{v}\in\Delta\cap\mbb{Z}^{n+1} $ with $ \bs{u}\prec \bs{v} $, $ \nu(\bs{u})\ll \nu(\bs{v}) $ (where $ \ll $ is understood as in Construction \ref{con:Mikhalkin_position}), unless either $ \{ \bs{u},\bs{v} \} = \{ \bs{w}^{\text{pr}}, \bs{w} \} $ or $ \{ \bs{u},\bs{v} \} = \{ \bs{w}, \bs{w}^{\text{suc}} \} $.
	
	Denote by $ \sum_{\bs{x}\in C_S\cap\mbb{Z}^{n+1}}\alpha_{\bs{x}}\bs{x}=\bs{0} $ the unique (up to multiplication by a scalar) linear relation of $ C_S\cap\mbb{Z}^{n+1} $ and by $ \bs{u}_{\max} $ the last (w.r.t. $ \prec $) point of $ C_S $.
	Since $ C_S $ is a linearity domain of $ \nu $, we have $ \sum_{\bs{x}\in C_S\cap\mbb{Z}^{n+1}} \alpha_{\bs{x}}\nu(\bs{x})=0 $.
	Thus $ \bs{u}_{\max} $ is either $ \bs{w} $ or $ \bs{w}^{\text{suc}} $ since otherwise, we would have $ \nu(\bs{u}_{\max})\gg \nu(\bs{v}) $ for all $ \bs{v}\in C_S\cap\mbb{Z}^{n+1}\setminus \{ \bs{u}_{\max} \} $.
	If $ \bs{u}_{\max}=\bs{w}^{\text{suc}} $ then, since  $ \nu(\bs{w}^{\text{suc}})\gg \nu(\bs{v}) $ for all $ \bs{v}\prec \bs{w} $, we get that $ \bs{w}\in C_S $ and of course this holds if $ \bs{u}_{\max}=\bs{w} $.
	
	Now, let $ \bs{v}\in C_S\cap\mbb{Z}^{n+1} $ with $ \bs{v}\prec \bs{w} $ and $ \bs{x}\in [\bs{v}, \bs{w}] $.
	Since $ C_S $ is convex, $ \bs{x}\in C_S $ meaning that $ \nu(\bs{x}) $ is a convex combination of $ \nu(\bs{v}) $ and $ \nu(\bs{w}) $. 
	Thus, if $ \bs{x}\ne\bs{v} $, we get $ \nu(\bs{x})\gg \nu(\bs{y}) $ for all $ \bs{y}\prec \bs{w}^{\text{pr}} $ meaning that $ \bs{x}=\bs{w}^{\text{pr}} $.
	This can only happen if $ C_S=\conv\{ \bs{v},\bs{w}^{\text{pr}},\bs{w} \} $ but then $ \bs{\eta}=\bs{\chi}_k $.
	We conclude that $ \bs{v}\in\partial_{\bs{w}}\Omega $, i.e. $ C_S\cap\Omega\subseteq \partial_{\bs{w}}\Omega $.
	
	Finally, we show that $ \bs{w}\ne \bs{u}_{\max} $.
	Indeed, since $ \partial_{\bs{w}}\Omega $ contained in an hyperplane by Lemma \ref{lemma:illuminated_calc}, and every proper subset of $ C_S\cap\mbb{Z}^{n+1} $ is linearly independent, $ |C_S\cap\mbb{Z}^{n+1}\setminus \Omega|>1 $ meaning that $ \bs{u}_{\max}=\bs{w}^{\text{suc}} $.
\end{proof}

\begin{prop}\label{prop:no_disconnected}
	\begin{sloppypar}
	There are no singular tropical hypersurfaces passing through $ \xpoints $ and corresponding to disconnected path as in option $ (b) $ of Lemma \ref{lemma:lattice_path}.
	\end{sloppypar}
\end{prop}

\begin{proof}
	\begin{sloppypar}
	Let $ \Gamma $ be a disconnected lattice path missing the edge $ [\bs{w}^{\text{pr}},\bs{w}] $.
	Assume, towards a contradiction, that $ S\in \text{Sing}^{\text{tr}}(\Delta, \xpoints) $ with $ \Gamma=P(S, \xpoints) $.
	Let $ \nu:\Delta\to\mbb{R} $ be the Legendre dual piecewise linear function, and $ C_S $ be the unique linearity domain of $ \nu $ that contains a circuit.
	Following the proof of Lemma \ref{lemma:lattice_path}, denote $ L:=\Sp(\bphi) $ and $ L\cap S = \{\bchi{1},\dots,\bchi{k} , \bs{\eta}, \bchi{k+1},\dots, \bchi{N-1} \} $.
	Write $ \bs{\eta} = M_{\bs{\eta}} \cdot \bphi $ for some $ M_k < M_{\bs{\eta}} < M_{k+1} $.
	\end{sloppypar}
	
	Let $ \delta = C_S\cap\mbb{Z}^{n+1}\setminus \{ \bs{w},\bs{w}^{\text{suc}} \} $.
	As in the proof of Lemma \ref{lemma:disconnected_path_illuminated}, we get that
	\begin{equation}\label{eq:no_disconnected_nu_w}
	\nu(\bs{w}) = \nu(\bs{w}^{\text{pr}}) + \pair{\bphi}{\bs{w}-\bs{w}^{\text{pr}}}M_{\bs{\eta}}
	\end{equation}
	and
	\begin{equation}\label{eq:no_disconnected_nu_v_s}
	\nu(\bs{w}^{\text{suc}}) = \nu(\bs{w}) + \pair{ \bphi}{\bs{w}^{\text{suc}}-\bs{w}}M_{k+1} = \nu(\bs{w}^{\text{pr}}) + \pair{ \bphi}{\bs{w}-\bs{w}^{\text{pr}}}M_{\bs{\eta}} + \pair{ \bphi}{\bs{w}^{\text{suc}}-\bs{w}}M_{k+1}
	\end{equation}
	Denote by $ \sum_{\bs{x}\in C_S\cap \mbb{Z}^{n+1}}\alpha_{\bs{x}} \bs{x}=\bs{0} $ the circuit relation, and assume, without loss of generality, that $ \alpha_{\bs{w}}+\alpha_{\bs{w}^{\text{suc}}}\ge 0 $. 
	Since $ \nu $ is linear on $ C_S $, $ \sum_{\bs{x}\in C_S}\alpha_{\bs{x}} \nu(\bs{x})=0 $.
	Substituting \eqref{eq:no_disconnected_nu_w} and  \eqref{eq:no_disconnected_nu_v_s} into this and using the fact that $ \nu(\bs{v})\ll \nu(\bs{w}) $ (in the sense of Construction \ref{con:Mikhalkin_position}) for all $ \bs{v}\prec \bs{w} $, we get that 
	\begin{equation*}
	(\alpha_{\bs{w}}+\alpha_{\bs{w}^{\text{suc}}})\pair{\bphi}{\bs{w}-\bs{w}^{\text{pr}}}M_{\bs{\eta}} + \alpha_{\bs{w}^{\text{suc}}}\pair{\bphi}{\bs{w}^{\text{suc}}-\bs{w}}M_{k+1} = 0.
	\end{equation*}
	Since $ 0<M_{\bs{\eta}}<M_{k+1} $,
	\begin{equation}\label{eq:disconected_path_alpha_ineq}
	0 < -\alpha_{\bs{w}^{\text{suc}}} < \alpha_{\bs{w}}\frac{\pair{\bphi}{\bs{w}-\bs{w}^{\text{pr}}}}{\pair{\bphi}{\bs{w}^{\text{suc}}-\bs{w}^{\text{pr}}}}
	\end{equation}
	in particular, $$ \alpha_{\bs{w}}+\alpha_{\bs{w}^{\text{suc}}}>\alpha_{\bs{w}}\left(1-\frac{\pair{\bphi}{\bs{w}-\bs{w}^{\text{suc}}}}{\pair{\bphi}{\bs{w}^{\text{suc}}-\bs{w}^{\text{pr}}}} \right)=\alpha_{\bs{w}}\frac{\pair{\bphi}{\bs{w}^{\text{suc}}-\bs{w}}}{\pair{\bphi}{\bs{w}^{\text{suc}}-\bs{w}^{\text{pr}}}}> 0 $$ so we can assume, without loss of generality, that $ \alpha_{\bs{w}}+\alpha_{\bs{w}^{\text{suc}}}=1 $.
	Substituting $ \alpha_{\bs{w}}=1-\alpha_{\bs{w}^{\text{suc}}} $ into \eqref{eq:disconected_path_alpha_ineq} we get
	\begin{equation*}
	0 < -\alpha_{\bs{w}^{\text{suc}}} < (1-\alpha_{\bs{w}^{\text{suc}}})\frac{\pair{\bphi}{\bs{w}-\bs{w}^{\text{pr}}}}{\pair{\bphi}{\bs{w}^{\text{suc}}-\bs{w}^{\text{pr}}}} \Rightarrow 0<-\alpha_{\bs{w}^{\text{suc}}} < \frac{\pair{\bphi}{\bs{w}-\bs{w}^{\text{pr}}}}{\pair{\bphi}{\bs{w^{(s)} - \bs{w}}}}
	\end{equation*}
	We can deduce that the point  $ -\sum_{\bs{x}\in \delta}\alpha_{\bs{x}} \bs{x}=\alpha_{\bs{w}}\bs{w}+\alpha_{\bs{w}^{\text{suc}}}\bs{w}^{\text{suc}} = \bs{w}-\alpha_{\bs{w}^{\text{suc}}}(\bs{w}-\bs{w}^{\text{suc}}) $ is contained in the open interval $$ I:=(\bs{w}, \bs{w}+\frac{\pair{\bphi}{\bs{w}-\bs{w}^{\text{pr}}}}{\pair{\bphi}{\bs{w}^{\text{suc}}-\bs{w}}}(\bs{w}-\bs{w}^{\text{suc}})) $$.
	
	Now, suppose that $ \ell $ is the first index with $ \bs{w}_{\ell} > 0 $, so $$ \bs{w}^{\text{suc}} = \bs{w} + (\bs{w}_{\ell} -1)\bs{e}^{(0)} - \bs{w}_{\ell}\bs{e}^{(\ell)} + \bs{e}^{(\ell+1)}. $$
	Let $ H $ be a supporting hyperplane of $ \Omega $ passing through $ \partial_{\bs{w}}\Omega $ (See Lemma \ref{lemma:illuminated_calc}).
	Since $ -\sum_{\bs{x}\in \delta}\alpha_{\bs{x}}\bs{x} = \alpha_{\bs{w}}\bs{w}+\alpha_{\bs{w}^{\text{suc}}}\bs{w}^{\text{suc}} $ is contained both in $ H $ and in the interval $ I $ it is enough to show this interval is disjoint from $ H $ to arrive at a contradiction.
	Since $ I $ is contained in the space where $ \bs{x}_i=\bs{w}_i $ for all $ i>\ell+1 $, we can consider only this space.
	On it, $ H $ is defined by the equation
	\begin{equation*}
	\theta(\bs{x}):=(1-\bs{w}_{\ell}-2\bs{w}_{\ell+1})\sum_{i=0}^{\ell-1}\bs{x}_i + (1-\bs{w}_{\ell+1})\bs{x}_{\ell}+(\bs{w}_{\ell}+1)\bs{x}_{\ell+1} = 0
	\end{equation*}
	which attain the positive value $ \bs{w}_{\ell}+\bs{w}_{\ell+1} $ on $ \bs{w} $.
	We are thus left with proving that $ \theta $ is positive on the other end of the interval $ I $ as well. 
	\begin{equation*}
	\begin{split}
	\theta\left (\bs{w}^{\text{suc}} - \bs{w}\right ) 
	=   \theta((\bs{w}_{\ell} -1)\bs{e}^{(0)} - \bs{w}_{\ell}\bs{e}^{(\ell)} + \bs{e}^{(\ell+1)}) & = \\
	=  \left( (\bs{w}_{\ell}-1)(1-\bs{w}_{\ell}-2\bs{w}_{\ell+1})-\bs{w}_{\ell}(1-\bs{w}_{\ell+1})+\bs{w}_{\ell}+1 \right) & = \\
	=(\bs{w}_{\ell} -1 -\bs{w}_{\ell}^2 +\bs{w}_{\ell} -2\bs{w}_{\ell}\bs{w}_{\ell+1} +2\bs{w}_{\ell+1} -\bs{w}_{\ell} +\bs{w}_{\ell}\bs{w}_{\ell+1} +\bs{w}_{\ell}+1 ) &= \\
	=  (2\bs{w}_{\ell}  -\bs{w}_{\ell}^2  -\bs{w}_{\ell}\bs{w}_{\ell+1} +2\bs{w}_{\ell+1} ) 
	=   (\bs{w}_{\ell}+\bs{w}_{\ell+1})(2 - \bs{w}_{\ell})
	\end{split}
	\end{equation*}
	
	If $ \bs{w}_{\ell}>1 $, we have $ \theta\left (\bs{w}+\frac{\pair{\bphi}{\bs{w}-\bs{w}^{\text{pr}})}}{\pair{\bphi}{\bs{w}^{\text{suc}}-\bs{w}}}(\bs{w}-\bs{w}^{\text{suc}}) \right ) \ge \theta(\bs{w}) > 0 $.
	If $ \bs{w}_{\ell+1}=0 $ then $ \Omega^{(\ell+1)} $ (see the notation in the beginning of the proof of Lemma \ref{lemma:illuminated_calc}) does not contain a point with $ \bs{x}_{\ell+1}>0 $.
	This means, since $ \bs{v}^{\bs{(s)}}_{\ell+1}>0 $, that $ I $ can not intersect with $ \Sp(\Omega) $.
	If $ \ell=0 $ and $ \bs{w}_1>0 $, we have $ \bs{w}^{\text{pr}} = \bs{w}+\bs{e}^{(0)}-\bs{e}^{(1)} $ and $ \bs{w}^{\text{suc}}=\bs{w}-\bs{e}^{(0)}+\bs{e}^{(1)} $ meaning that $ {\pair{\bphi}{\bs{w}-\bs{w}^{\text{pr}})}}={\pair{\bphi}{\bs{w}^{\text{suc}}-\bs{w}}} $ and so $ I $ and $ H $ are disjoint ($ H $ does contain an endpoint of $ I $, namely $ \bs{w}^{\text{pr}} $, but not an interior point).
		
	Finally, if $ \ell>0 $ and $ \bs{w}_{l+1}>0 $, we can compute $ \bs{w}^{\text{pr}} = \bs{w}+\bs{e}^{(\ell-1)}-\bs{e}^{(\ell)} $ and \linebreak $ {\bs{w}^{\text{suc}} = \bs{w}-\bs{e}^{(\ell)}+\bs{e}^{(\ell+1)}} $.
	This means that 
	$$ \pair{\bphi}{\bs{w}^{\text{suc}} - \bs{w}} - \pair{\bphi}{\bs{w} - \bs{w}^{\text{pr}}} = \pair{\bphi}{\bs{w}^{\text{pr}}} - \pair{\bphi}{\bs{w}+\bs{e}^{(\ell)}-\bs{e}^{(\ell+1)}} > 0 $$
	since $ \bs{w}-\bs{e}^{(\ell)}+\bs{e}^{(\ell+1)}\prec \bs{w}^{\text{pr}} $.
	Thus 
	\[ \theta\left (\bs{w}+\frac{\pair{\bphi}{\bs{w}-\bs{w}^{\text{pr}})}}{\pair{\bphi}{\bs{w}^{\text{suc}}-\bs{w}}}(\bs{w}-\bs{w}^{\text{suc}})\right ) = (\bs{w}_{\ell}+\bs{w}_{\ell+1})\left (1 - \frac{\pair{\bphi}{\bs{w}-\bs{w}^{\text{pr}})}}{\pair{\bphi}{\bs{w}^{\text{suc}}-\bs{w}}} \right )> 0. \]
\end{proof}

\subsection{Connected lattice path and Mikhalkin condition}
In this Section we begin the investigation of singular hypersurfaces that pass through $ \xpoints $ and correspond to connected lattice path, which will continue throughout Section \ref{sec:construction_graded_circuits}.
In view of Lemma \ref{lemma:origin_sing_condition}, we translate the tropical hypersurface $ S $ such that its singularity is at the origin. 
Since such a translation corresponds to adding a linear functional to the coefficients of the defining tropical polynomial of $ S $, we would like to replace the notion of Mikhalkin position of $ \xpoints $ by one that is invariant to this operation.
This is the reason for the following definition.

\begin{Def}\label{def:mikhalkin_condition}
	Let $ \bs{w}\in \Delta\cap\mbb{Z}^{n+1} $ and $ \nu $ is a piecewise linear convex function on $ \Delta $. 
	We will say that $ \nu $ \emph{satisfies Mikhalkin condition away from $ \bs{w} $} if for any circuit $ Q\subseteq \Delta\cap\mbb{Z}^{n+1}\setminus\{\bs{w}\} $, given the unique (up to multiplication by a scalar) linear dependency $ \sum_{\bs{v}\in Q}\alpha_{\bs{v}} \bs{v} =\bs{0} $ we have
	\[ \sum_{\bs{v}\in Q}\alpha_{\bs{v}} \nu(\bs{v}) > 0 \iff \alpha_{\bs{u}_{\max}}>0 \]
	where $ \bs{u}_{\max} $ is the last (w.r.t. $ \prec $) vector of $ Q $.
\end{Def}

It is immediate that this definition satisfies the desired invariance under adding a linear functional to $ \nu $.
\begin{lemma}\label{lemma:mikh_cond_add_functional}
	If $ \nu:\Delta\to\mbb{R} $ satisfies Mikhalkin condition away from $ \bs{w} $ and $ \Lambda:\mbb{R}^{n+1}\to\mbb{R} $ is a linear functional, then $ \nu+\Lambda $ satisfies Mikhalkin condition away from $ \bs{w} $.
\end{lemma}

\begin{proof}
	If $ Q\subseteq \Delta\cap\mbb{Z}^{n+1} \setminus \{ \bs{w}\} $ is a circuit with relation $ \sum_{\bs{v}\in Q}\alpha_{\bs{v}}\bs{v}=\bs{0} $ s.t $ \alpha_{\bs{u}_{\max}}>0 $ then 
	\[ \sum_{\bs{v}\in Q}\alpha_{\bs{v}}(\nu(\bs{v})+\Lambda(\bs{v})) = \sum_{\bs{v}\in Q}\alpha_{\bs{v}}\nu(\bs{v})) + \Lambda\left( \sum_{\bs{v}\in Q}\alpha_{\bs{v}}\bs{v}\right) = \sum_{\bs{v}\in Q}\alpha_{\bs{v}}\nu(\bs{v})>0. \]
\end{proof}

Definition \ref{def:mikhalkin_condition} in fact extends the definition of Mikhalkin position, as the following lemma shows.

\begin{lemma}\label{lemma:mikhalkin_pos_imply_cond}
	Let $ \bs{w}\in \Delta\cap\mbb{Z}^{n+1} $ and $ \Gamma_{\bs{w}}:=P(\Delta\cap\mbb{Z}^{n+1}\setminus \{ \bs{w} \}) $ the corresponding connected lattice path.
	Suppose that $ S $ is a tropical hypersurface passing through $ \xpoints $ whose corresponding lattice path is $ \Gamma_{\bs{w}} $ and let $ \nu:\Delta\to \mbb{R} $ its Legendre dual piecewise linear convex function.
	Then $ \nu $ satisfies Mikhalkin condition away from $ \bs{w} $.
\end{lemma}

\begin{proof}
	Let $ Q\subseteq \Delta\cap\mbb{Z}^{n+1}\setminus \{ \bs{w}\} $ be a circuit with relation $ \sum_{\bs{v}\in Q}\alpha_{\bs{v}}\bs{v} = \bs{0}  $ and suppose that $ \alpha_{\bs{u}_{\max}}>0 $ where $ \bs{u}_{\max}:=\max_{\prec}Q $.
	For $ \bs{y}\in \Delta\cap\mbb{Z}^{n+1} $, if $ [\bs{y},\bs{u}_{\max}]\subseteq \Gamma_{\bs{w}} $ then $ \nu(\bs{u}_{\max}) = \nu(\bs{y})+M_k\pair{\bphi}{\bs{u}_{\max}-\bs{y}} $ for some $ 1\le k \le N-1 $, and $ M_i $ with $ i\ge k $ can not affect $ \nu(\bs{v}) $ for any $ \bs{v}\in Q\setminus \{ \bs{u}_{\max} \} $.
	Thus the maximal $ M_i $ appearing in $ \sum_{\bs{v}\in Q}\alpha_{\bs{v}}\nu(\bs{v}) $ is $ M_k $ and it appear with a positive coefficient, so by Construction \ref{con:Mikhalkin_position}, $ \sum_{\bs{v}\in Q}\alpha_{\bs{v}}\nu(\bs{v})>0 $.
\end{proof}

In particular, if we pick any $ S $ that correspond to a connected lattice path, by choosing the coefficient of $ \bs{w} $ to satisfy \eqref{eq:connected_path_ineq},  we get the following corollary.

\begin{cor}\label{cor:exist_mikh_condition}
	For every $ \bs{w}\in\Delta\cap\mbb{Z}^{n+1} $ there exist a piecewise linear function $ \nu:\Delta \to \mbb{R} $ that satisfies Mikhalkin condition away from $ \bs{w} $.
\end{cor}

\section{Construction of tropical singular hypersurfaces} \label{sec:construction_graded_circuits}
In this section we will construct the examples of tropical singular hypersurfaces that will allow us to exhaust assymptotically all possible singular hypersurfaces (see Section \ref{sec:discriminant_degree}).

In Section \ref{subsec:comb1} we will define graded circuits (Definition \ref{def:graded_circuit}) which is the combinatorial object that corresponds to singular points of maximal-dimensional type on tropical \linebreak hypersurfaces and recall low dimensional examples as classified in \cite{Markwig2012_trop_curve_sing} and \cite{Markwig2012_trop_surf_sings}.
In Section \ref{subsec:comb2} we will give straightforward combinatorial condition for a graded circuit to correspond to a singular point of a tropical hypersurface that passes through points satisfying Mikhalkin condition.
In Section \ref{subsec:comb3} we will construct the circuits that will appear in our enumeration.
Finally, in Section \ref{subsec:comb4} we will show how the circuits considered in Section \ref{subsec:comb3} can be glued together into graded circuits that admit Mikhalkin condition.

\subsection{Definitions and first examples}\label{subsec:comb1}

In view of Lemma \ref{lemma:origin_sing_condition} and Corollary \ref{cor:nu_image_size} it is natural to define the following:

\begin{Def}\label{def:graded_circuit} 
	Let $ \bs{w}\in \Delta\cap\mbb{Z}^{n+1}d $. 
	A sequence of subsets $ C_i\subseteq  \Delta_d^{(n)} \cap\mbb{Z}^{n+1} (i\in \mbb{Z}_{\ge 0}) $ is called \emph{a graded circuit centered at $ \bs{w} $} if it satisfies $ C_0 = \{ \bs{w}\} $, each $ C_i $ is linearly independent and for every $ i>0 $ there exist a unique, up to multiplication by a scalar, linear combination
\begin{equation}\label{eq:def_GC_unique_linear_comb}
	\sum_{\bs{x}\in C_i}\alpha_{\bs{x}} \bs{x} \in \Sp(\bigcup_{j<i}C_j)
\end{equation}
	with $ \alpha_{\bs{x}}\ne 0 $ for all $ \bs{x}\in C_i $. 
	
	Each $ C_i $ is called a \emph{level} of the graded circuit. 
	The \emph{span} of a graded circuit $ C $ is 
	\[ \Sp(C) := \Sp(\bigcup_{i} C_i). \]
	We set $ \dim C := \dim \Sp(C)-1 $. 
	
	A level $ C_i $ is \emph{trivial} if it consist of a single point. 
	Given a graded circuit $ C $, we can construct the \emph{simplification of $ C $} by throwing out all trivial levels, it has the same span as the original circuit. 
	Two graded circuits will be called \emph{equivalent} if they have the same simplification. 
	The number of non trivial levels of $ C $ is called the \emph{height of $ C $} and denoted $ ht(C) $.
	
	For a level $ C_i $, if the unique linear combination \eqref{eq:def_GC_unique_linear_comb} has $ \sum_{\bs{x}\in C_i}\alpha_{\bs{x}} = 0 $, we will say $ C_i $ \emph{is parallel to previous levels}. 
	If this is not the case for all $ C_i\in C $, we will say \emph{$ C $ has no parallel levels}.

\end{Def}

\begin{Example}\label{ex:graded_circuit_height_1}
If $ C $ is a graded circuit of height 1 centered at $ \bs{w} $, then $ C_0\cup C_1 = \{\bs{w} \}\cup C_1 $ is a circuit.
Conversely, if $ C $ is a circuit containing $ \bs{w} $ then $ \left\langle \{\bs{w}\}, C\setminus \{\bs{w}\} \right\rangle $ is a graded circuit of height 1 centered at $ \bs{w} $.
We thus refer to graded circuits of height 1 as circuits.
\end{Example}

\begin{Def}\label{def:circuit_admit_nu}
	Given a function $ \nu:\Delta\cap\mbb{Z}^{n+1}\to \mbb{R} $, the \emph{levels} of $ \nu $ is the subsets $ {\{\nu^{-1}(a) \; | \; a\in\mbb{R} \}} $.
	Let $ C $ be a graded circuit centered at $ \bs{w} $.
	We will say $ C $ \emph{admits} $ \nu $ if the levels of $ \nu|_{\Sp(C)} $ form a graded circuit centered at $ w $ and equivalent to $ C $.
	We will say that $ C $ \emph{admits Mikhalkin condition} if it admits some function $ \nu $ that satisfies Mikhalkin condition away from $ \bs{w} $ (see Definition \ref{def:mikhalkin_condition}). 
\end{Def}

\begin{Example}\label{ex:all_graded_circuits_of_dim2}
	The graded circuits of dimension 2 and 3 admitting some convex, piecewise linear function with linearity domains lattice polytopes, were classified in \cite{Markwig2012_trop_curve_sing} (for dimension 2) and \cite{Markwig2012_trop_surf_sings} (for dimension 3).
	
	In dimension 2 their classification give us the following:
	\begin{enumerate}
		\item 
		The height $ ht(C)=1 $ (i.e. $ C $ is a circuit) and $ \conv(C_0\cup C_1) $ is a parallelogram (see Figure \ref{fig:graded_circuits_dim2} (a)).
		
		\item
		The height $ ht(C)=1 $ (i.e. $ C $ is a circuit) and $ \conv(C_1) $ is a triangle containing $ \bs{w} $ in its interior (see Figure \ref{fig:graded_circuits_dim2} (b)).
		 
		\item
		The height $ ht(C)=2 $, $ \conv(C_0\cup C_1) $ is a segment of lattice length 2 and $ C_2 $ contains a point on either side of $ \Sp(C_1) $ (see Figure \ref{fig:graded_circuits_dim2} (c)).
		
		\item
		The height $ ht(C)=2 $, $ \conv(C_0\cup C_1) $ is a segment of lattice length 2 and $ C_2 $ contains 2 points on the same side of $ \Sp(C_1) $ (see Figure \ref{fig:graded_circuits_dim2} (d)).
	\end{enumerate}

	Note that the highest level of the last example is parallel to previous levels and the rest of the examples have no parallel levels.
\end{Example}

\begin{figure}
	
	\begin{tikzpicture}[scale=1.2]
	
	\begin{scope}[scale = 0.25]
	\draw (0,0) -- (0,10) -- (10,0) -- (0,0);
	\filldraw (3,3) circle (0.3);
	\filldraw [fill=gray] (3,2) circle (0.3);
	\filldraw [fill=gray] (4,3) circle (0.3);
	\filldraw [fill=gray] (4,2) circle (0.3);
	
	\draw (5, -2) node {$ (a) $};
	\end{scope}
	
	\begin{scope}[scale = 0.25, xshift=400]
	\draw (0,0) -- (0,10) -- (10,0) -- (0,0);
	\filldraw (3,3) circle (0.3);
	\filldraw [fill=gray] (2,2) circle (0.3);
	\filldraw [fill=gray] (4,3) circle (0.3);
	\filldraw [fill=gray] (3,4) circle (0.3);
	
	\draw (5, -2) node {$ (b) $};
	\end{scope}
	
	\begin{scope}[scale = 0.25, xshift=800]
	\draw (0,0) -- (0,10) -- (10,0) -- (0,0);
	\filldraw (3,3) circle (0.3);
	\filldraw [fill=gray] (3,2) circle (0.3);
	\filldraw [fill=gray] (3,4) circle (0.3);
	\filldraw [fill=white] (4,2) circle (0.3);
	\filldraw [fill=white] (2,4) circle (0.3);
	
	\draw (5, -2) node {$ (c) $};
	\end{scope}
	
	\begin{scope}[scale = 0.25, xshift=1200]
	\draw (0,0) -- (0,10) -- (10,0) -- (0,0);
	\filldraw (3,3) circle (0.3);
	\filldraw [fill=gray] (3,2) circle (0.3);
	\filldraw [fill=gray] (3,4) circle (0.3);
	\filldraw [fill=white] (4,2) circle (0.3);
	\filldraw [fill=white] (4,3) circle (0.3);
	
	\draw (5, -2) node {$ (d) $};
	\end{scope}
	
	\end{tikzpicture}
	
	\caption{Examples of graded circuits that admit some convex, piecewise linear function with lattice polytopes as linearity domains.
	(a) Parallelogram; (b) Triangle with inner point; (c) Segment with the the points of $ C_2 $ being on either side; (d) Segment with the points of $ C_2 $ on the same side. The black circle in each example denotes $ \bs{w} $, gray circles denote points of $ C_1 $ and white circles denote points of $ C_2 $.}
	\label{fig:graded_circuits_dim2}
\end{figure}
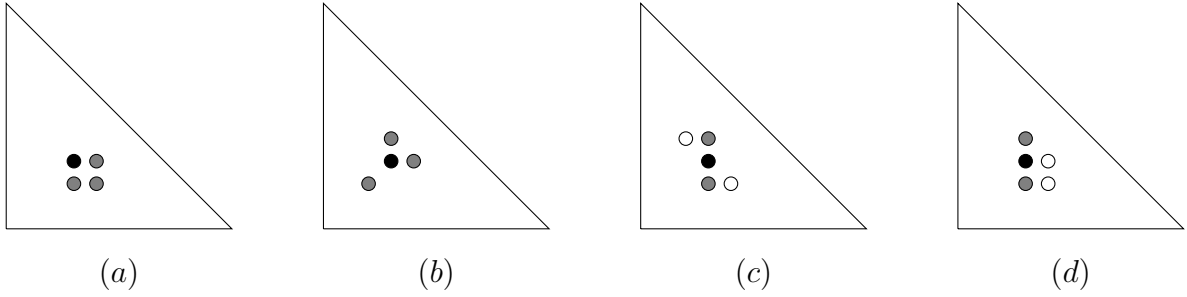
\begin{remark}\label{rem:no_circ_classification}
	Example \ref{ex:all_graded_circuits_of_dim2} relies on the classification, up to unimodular transformation, of lattice circuits. 
	This approach is impossible in higher dimensions, since even the empty lattice simplices in dimension 4 were classified only recently in \cite{Santos_Valino_class_dim4_simpl}, and the classification of lattice circuits seems out of reach with current methods.
\end{remark}

\subsection{Graded circuits admitting Mikhalkin condition}\label{subsec:comb2}

In this subsection we will study the properties of graded circuits that admit Mikhalkin condition and will give a combinatorial characterization of them.

We will need the following notation.
\begin{Def}
	Denote 
	\[ H_0=\{ \bs{x}\in \mbb{R}^{n+1} \; | \; \sum \bs{x}_i=0 \} \].
	For a graded circuit $ C $ and $ i>0 $ we denote
	\[ \partial C_i = \left \{ \sum_{\bs{x}\in C_i}\alpha_{\bs{x}} \bs{x} | \sum_{\bs{x}\in C_i}\alpha_{\bs{x}}=0 \right \}\subseteq H_0 \] 
	and
	\[ \partial C = \sum_{i=1}^{ht(C)} \partial C_i\subseteq H_0 \]
	where the sum denotes sum of the linear spaces $ \partial C_i $.
	
	We denote by $ \overline{C}_k $ the truncated graded circuit with $ ht(\overline{C}_k)=k $ and $ (\overline{C}_k)_i = C_i $ for $ i\le k $. 
\end{Def}

A useful property of graded circuits that admits Mikhalkin condition is that they do not have parallel levels (Lemma \ref{lemma:graded_circ_admit_mikhalkin_no_parallel}).
This property will be used through the following lemma and its corollary \ref{cor:graded_circuit_linear_combination_with_sum_1}.

\begin{lemma}\label{lemma:difference_set_no_parallel_circuit}
	Let $ C $ be a graded circuit without parallel levels, then for all $ i>0 $,
	\begin{enumerate}
		\item \[ \partial \overline{C}_i = \bigoplus_{j\le i} \partial C_j = \Sp(\overline{C}_i)\cap H_0. \]
		\item \[ \Sp(\overline{C}_i) = \Sp(C_i)\oplus \partial\overline{C}_{i-1}. \]
	\end{enumerate}
\end{lemma}

\begin{proof}
	\begin{enumerate}
		\item 
		Assume that $ \bs{\eta} \in \left( \sum_{j<i}\partial C_j \right)\cap \partial C_i   $.
		Then $ \bs{\eta} = \sum_{\bs{x}\in C_i}\alpha_{\bs{x}} \bs{x} $ with $ \sum \alpha_{\bs{x}}=0 $.
		But since $ \bs{\eta}\in \Sp(\bigcup_{j<i}C_j) $ the coefficients $ \alpha_{\bs{x}} $ are proportional to the coefficients in \eqref{eq:def_GC_unique_linear_comb} and because $ C_i $ is not parallel to previous levels, $ \alpha_{\bs{x}}=0 $ for all $ \bs{x}\in C_i $, meaning that $ \bs{\eta} = 0 $.
		
		Now, let $ \bs{\eta} \in H_0\cap \Sp(\overline{C}_i) $, so $ \bs{\eta }= \sum_{j\le i}\sum_{\bs{x}\in C_j} \beta_{\bs{x}} \bs{x} $. For every $ j\le i $ exist $ \sum_{k\le j}\sum_{\bs{x}\in C_k} \alpha_{\bs{x}} \bs{x}=0 $ and since $ C $ has no parallel levels, $ \sum_{\bs{x}\in C_j}\alpha_{\bs{x}}\ne 0 $. So by adding correct multiple of this combination to $ \bs{\eta} $ we can assume that $ \sum_{\bs{x}\in C_j}\beta_{\bs{x}} =0 $ for all $ 2\le j\le i $. Since $ \bs{\eta} \in H_0 $, $ \sum_{\bs{x}\in C_1}\beta_{\bs{x}}=0 $ as well, meaning that $ \bs{\eta}\in\partial \overline{C}_i $.
		
		\item 
		By $ (1) $, 
		\[ \Sp(C_i)\cap \partial \overline{C}_{i-1} = (\Sp(C_i)\cap H_0) \cap \bigoplus_{j< i} \partial C_j = \partial C_i\cap \bigoplus_{j< i} \partial C_j = {0}.   \]
		Now, let $ \bs{\eta}\in \Sp(\overline{C}_i) $, so $ \bs{\eta}=\sum_{j\le i}\sum_{\bs{x}\in C_j}\beta_{\bs{x}} \bs{x} $. As in the proof of $ (1) $, we can assume, without loss of generality, that $ \sum_{\bs{x}\in C_j}\beta_{\bs{x}}=0 $ for all $ j<i $. Then $$ \sum_{j<i}\sum_{\bs{x}\in C_j}\beta_{\bs{x}} \bs{x}\in \partial \overline{C}_{i-1} $$ and $ \sum_{\bs{x}\in C_i}\beta_{\bs{x}} \bs{x}\in \Sp(C_i) $, so the assertion follows.
	\end{enumerate}
\end{proof}

\begin{cor}\label{cor:graded_circuit_linear_combination_with_sum_1}
	Let $ C $ be a graded circuit with no parallel levels. Then for all $ h\le ht(C) $ there exist a linear combination $ \sum_{\bs{x}\in C_h} m_{\bs{x}} \bs{x}\in \Sp(\overline{C}_{h-1} ) $ s.t. $ \sum m_{\bs{x}}=1 $.
\end{cor}

\begin{proof}
	Let $ \sum _{\bs{x}\in C_h}\alpha_{\bs{x}} \bs{x}\in \Sp\left (\overline{C}_{h-1}\right ) $. By Lemma \ref{lemma:difference_set_no_parallel_circuit} (1), $ \sum_{\bs{x}\in C_h}\alpha_{\bs{x}} \ne 0 $ so after re-scaling we get the desired result.
\end{proof}

\begin{lemma}\label{lemma:graded_circ_admit_mikhalkin_no_parallel}
	A graded circuit that admits Mikhalkin condition has no parallel levels.
\end{lemma}

\begin{proof}
	Let $ \nu: \Delta^{(n)} \to \mbb{R} $ satisfying Mikhalkin condition and admitted by $ C $.
	Assume, towards a contradiction, that $ C_h $ is the first level which is parallel to previous levels. 
	Then there exist a linear combination $ \sum_{\bs{x}\in C_h}\alpha_{\bs{x}} \bs{x}\in \Sp\left (\overline{C}_{h-1}\right ) $ with $ \sum \alpha_{\bs{x}}=0 $.
	Since $ \overline{C}_{h-1} $ has no parallel levels, by Lemma \ref{lemma:difference_set_no_parallel_circuit} 
	\[ \sum_{\bs{x}\in C_h} \alpha_{\bs{x}} \bs{x} = \sum_{0< i< h}\sum_{\bs{x}\in C_i} \beta_{\bs{x}} \bs{x} \]
	with $ \sum_{\bs{x}\in C_i}\beta_{\bs{x}}=0 $ for all $ 0<i<h $.
	But then 
	\[ \sum_{\bs{x}\in C_h}\alpha_{\bs{x}} \nu(\bs{x}) - \sum_{0<i<h}\sum_{\bs{x}\in C_i}\beta_{\bs{x}} \nu(\bs{x}) = \left( \sum_{\bs{x}\in C_h}\alpha_{\bs{x}} \right) \nu(C_h) - \sum_{0<i <h}\left( \sum_{\bs{x}\in C_i} \beta_{\bs{x}} \right) \nu(C_i) = 0 \]
	in contradiction to Mikhalkin condition on $ \nu $.  
\end{proof}

Using this, we can now express the fact that a graded circuit $ C $ admits Mikhalkin condition without quantifying on $ \nu $.

\begin{lemma}\label{lemma:graded_circ_admit_mikhalkin_comb_cond}
	A graded circuit $ C $ centered at $ \bs{w} $ admits Mikhalkin condition if and only if it has no parallel levels and for all levels $ C_h\in C $ and all $ \bs{v}\in \Sp(\overline{C}_h)\cap\Delta\cap\mbb{Z}^{n+1}\setminus C_h $ the coefficient of the last (w.r.t. $ \prec $) point appearing in the unique (see Lemma \ref{lemma:difference_set_no_parallel_circuit}) linear combination 
\begin{equation}\label{eq:lemma_mikh_GC_comb_cond}
	\bs{v} -  \sum_{0 < i\le h}\sum_{\bs{x}\in C_i}\alpha_{\bs{x}} \bs{x} = 0
\end{equation}
	with $ \sum_{\bs{x}\in C_i}\alpha_{\bs{x}}=0 $ for all $ i<h  $, is positive.
\end{lemma}

\begin{proof}
	The fact that a graded circuit that admits Mikhalkin condition has no parallel levels was proved in Lemma \ref{lemma:graded_circ_admit_mikhalkin_no_parallel}.
	Suppose $ C $ admits a function $ \nu:\Delta \to \mbb{R} $ satisfying Mikhalkin condition and $ \bs{v}\in \Sp(\overline{C}_h)\cap\Delta\cap\mbb{Z}^{n+1}\setminus C_h $ for some level $ C_h\in C $.
	Then 
	\[ \nu(\bs{v}) - \sum_{0< i\le h}\sum_{\bs{x}\in C_i}\alpha_{\bs{x}} \nu(\bs{x}) = \nu(\bs{v}) - \sum_{0<i\le h}\left( \sum_{\bs{x}\in C_i}\alpha_{\bs{x}} \right)\nu(C_i)=\nu(\bs{v})-\nu(C_h)>0, \]
	so the coefficient of the last point has to be positive.
	
	Now, suppose $ C $ is a graded circuit centered at $ \bs{w} $ with no parallel levels and \eqref{eq:lemma_mikh_GC_comb_cond} satisfied for all $ \bs{v},h $.
	Pick arbitrary piecewise linear convex function $ \nu_0:\Delta\to\mbb{R} $ satisfying Mikhalkin condition away from $ \bs{w} $, see Corollary \ref{cor:exist_mikh_condition}. 
	Since $ \partial C = \bigoplus \partial C_i $ by Lemma \ref{lemma:difference_set_no_parallel_circuit}, there exist a linear function $ \Lambda:\Delta\to\mbb{R} $ s.t. $ \nu:=\nu_0+\Lambda $ is constant on each $ C_i $. 
	By Lemma \ref{lemma:mikh_cond_add_functional}, $ \nu $ satisfies Mikhalkin condition. 
	Finally, \eqref{eq:lemma_mikh_GC_comb_cond} ensures that the levels of $ \nu $ form a graded circuit which is equivalent to $ C $.
\end{proof}

\begin{remark}\label{remark:v_last_in_comb_cond}
	The coefficient of $ \bs{v} $ in \eqref{eq:lemma_mikh_GC_comb_cond} is 1 which is positive, meaning that whenever $ \bs{v} $ is the last point in \eqref{eq:lemma_mikh_GC_comb_cond}, the condition is satisfied.
\end{remark}

\subsection{Elementary Mikhalkin circuits}\label{subsec:comb3}

In this subsection we will describe the circuits we will use in our construction of graded circuits and show that they admit Mikhalkin condition.

\begin{Example}\label{ex:type_I_circ}
	For an index $ 0\le \ell< n $, let $ \bs{w}_i=0 $ for all $ i<\ell $, and $ \bs{w}_{\ell},\bs{w}_{\ell+1}\ne 0 $. 
	The graded circuit $ C $ of height 1 centered at $ \bs{w} $ with 
	\[C_1=\{\bs{u}^{(-1)}:= \bs{w}-\bs{e}^{(\ell+1)}+\bs{e}^{(\ell)}, \bs{u}^{(1)}:=\bs{w}+\bs{e}^{(\ell+1)}-\bs{e}^{(\ell)} \},\]
	see Figure \ref{fig:elementary_mikhalkin_circuits_dim2} (a,b), clearly satisfy the conditions of Lemma \ref{lemma:graded_circ_admit_mikhalkin_comb_cond}.
	We will call it \emph{elementary Mikhalkin circuit of type I and index $ \ell $}.
	The linear combination \eqref{eq:def_GC_unique_linear_comb} is $ \bs{u}^{(1)}+\bs{u}^{(-1)} $.
	Note that the circuit depicted in Figure \ref{fig:elementary_mikhalkin_circuits_dim2} (b) does not give rise to a singularity, as shown in \cite{Markwig2012_trop_curve_sing}.
\end{Example}

\begin{Example}\label{ex:type_II_circ}
	For an index $ 0<\ell<n $, let $ \bs{w}_i=0 $ for all $ i<\ell $ and $ \bs{w}_{\ell},\bs{w}_{\ell+1}\ne 0 $.
	The graded circuit $ C $ of height 1 centered at $ \bs{w} $ with 
	\[C_1=\{ \bs{u}^{(-2)}:= \bs{w}-\bs{e}^{(\ell+1)}+\bs{e}^{(\ell)}, \bs{u}^{(-1)}:= \bs{w}-\bs{e}^{(\ell)}+\bs{e}^{(\ell-1)}, \bs{u}^{(1)}:=\bs{w}+\bs{e}^{(\ell+1)}-2\bs{e}^{(\ell)}+\bs{e}^{(\ell-1)} \}\]
	(see Figure \ref{fig:elementary_mikhalkin_circuits_dim2} (c)) admits Mikhalkin condition, which we will show using Lemma \ref{lemma:graded_circ_admit_mikhalkin_comb_cond}.
	Let $ \bs{v}\in \Sp(C_1)\cap\Delta\cap\mbb{Z}^{n+1} \setminus \{\bs{w} \} $ with 
	\begin{equation}\label{eq:example_type_II_Mikh_cond} 
	\bs{v} - \alpha_{-2}\bs{u}^{(-2)}-\alpha_{-1}\bs{u}^{(-1)}-\alpha_1\bs{u}^{(1)} =0.
	\end{equation}
	If $ \alpha_{1}>0 $ then $ \bs{v}\succeq \bs{u}^{(1)} $ meaning that the condition of Lemma \ref{lemma:graded_circ_admit_mikhalkin_comb_cond} is satisfied.
	If $ \alpha_1 = 0 $, $ \bs{v}\in \Sp\{\bs{u}^{(-2)}, \bs{u}^{(-1)}\} $ and in any case $ \bs{v}\succeq \bs{u}^{(-1)} $ meaning that the coefficient of the last vector in  \eqref{eq:example_type_II_Mikh_cond} is positive.
	Finally, if $ \alpha_1<0 $, the coefficient of the last vector in \eqref{eq:example_type_II_Mikh_cond} is $ -\alpha_1>0 $. 
	
	We will call such a graded circuit \emph{elementary Mikhalkin circuit of type II and index $ \ell $}. The linear combination \eqref{eq:def_GC_unique_linear_comb} is $ \bs{u}^{(-2)}+\bs{u}^{(1)}-\bs{u}^{(-1)} $.
\end{Example}

The next example will require the notion of Mikhalkin triangulation which is the construction of \cite[Lemma 3.6]{Shustin2017_enum}.

\begin{lemma}\label{lemma:mikhalkin_triangulation}
	Let $ Q\subseteq \Delta $ be a lattice polytope not containing $ \bs{w} $.
	There exist a unique subdivision $ \sigma  $ of $ Q $ s.t every $ \nu:\Delta \to \mbb{R}$ satisfying Mikhalkin condition away from $ \bs{w} $ has cells of $ \sigma  $ as linearity domains. 
	The subdivision $ \sigma $ can be found applying successive applications of smooth extension as in \cite[Example 3.5]{Shustin2017_enum}.
\end{lemma}

\begin{proof}
	See \cite[Lemma 3.6]{Shustin2017_enum}.
\end{proof}

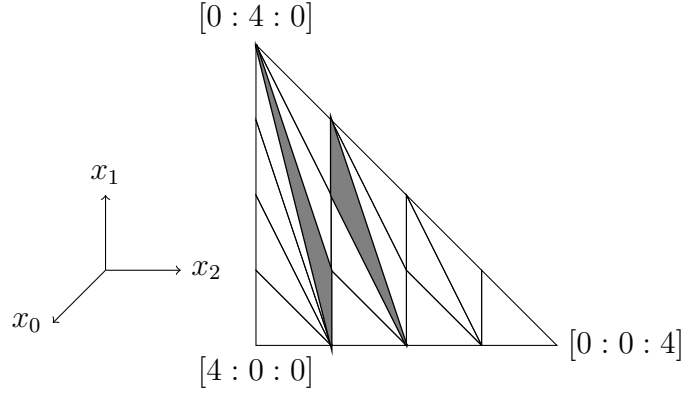
\begin{figure}[h]
	
	\begin{tikzpicture}
	\draw (1, 0) -- (0, 1) -- (0, 0) -- (1, 0);
	\draw (1, 0) -- (0, 1) -- (0, 2) -- (1, 0);
	\draw (1, 0) -- (0, 3) -- (0, 2) -- (1, 0);
	\draw (1, 0) -- (0, 3) -- (0, 4) -- (1, 0);
	\draw (1, 0) -- (1, 1) -- (0, 4) -- (1, 0);
	\draw (1, 2) -- (1, 1) -- (0, 4) -- (1, 2);
	\draw (1, 2) -- (1, 3) -- (0, 4) -- (1, 2);
	\draw (1, 0) -- (2, 0) -- (1, 1) -- (1, 0);
	\draw (1, 2) -- (2, 0) -- (1, 1) -- (1, 2);
	\draw (1, 2) -- (2, 0) -- (1, 3) -- (1, 2);
	\draw (2, 1) -- (2, 0) -- (1, 3) -- (2, 1);
	\draw (2, 2) -- (1, 3) -- (2, 1) -- (2, 2);
	\draw (3, 0) -- (2, 0) -- (2, 1) -- (3, 0);
	\draw (3, 0) -- (2, 2) -- (2, 1) -- (3, 0);
	\draw (3, 0) -- (3, 1) -- (2, 2) -- (3, 0);
	\draw (3, 0) -- (4, 0) -- (3, 1) -- (3, 0);
	
	\draw [line width=1pt] (0,4) -- (1,0) -- (1,1) -- (0,4);
	\draw [line width=1pt] (2,0) -- (1,2) -- (1,3) -- (2,0);
	\fill [fill = gray] (0,4) -- (1,0) -- (1,1) -- (0,4);
	\fill [fill = gray] (2,0) -- (1,2) -- (1,3) -- (2,0);
	
	\draw (0,0) node [below] {$ [4:0:0] $};
	\draw (4,0) node [anchor=west] {$ [0:0:4] $};
	\draw (0,4) node [above] {$ [0:4:0] $};
	
	\draw [->] (-2,1) -- (-2,2) node [above] {$ x_1 $ };
	\draw [->] (-2,1) -- (-1,1) node [anchor=west] {$ x_2 $ };
	\draw [->] (-2,1) -- (-2.7,0.3) node [anchor=east] {$ x_0 $ };
	\end{tikzpicture}
	
	\caption{Mikhalkin triangulation of $ \Delta_4^{(2)} $. The vertical direction is $ x_1 $, the horizontal direction is $ x_2 $ and $ x_0 $ is the distance from the diagonal edge of $ \Delta $. The grayed triangles are the 2 options referenced in Example \ref{ex:all_elementary_circuits_dim3} }
	\label{fig:Mikhalkin_triangulation}
\end{figure}

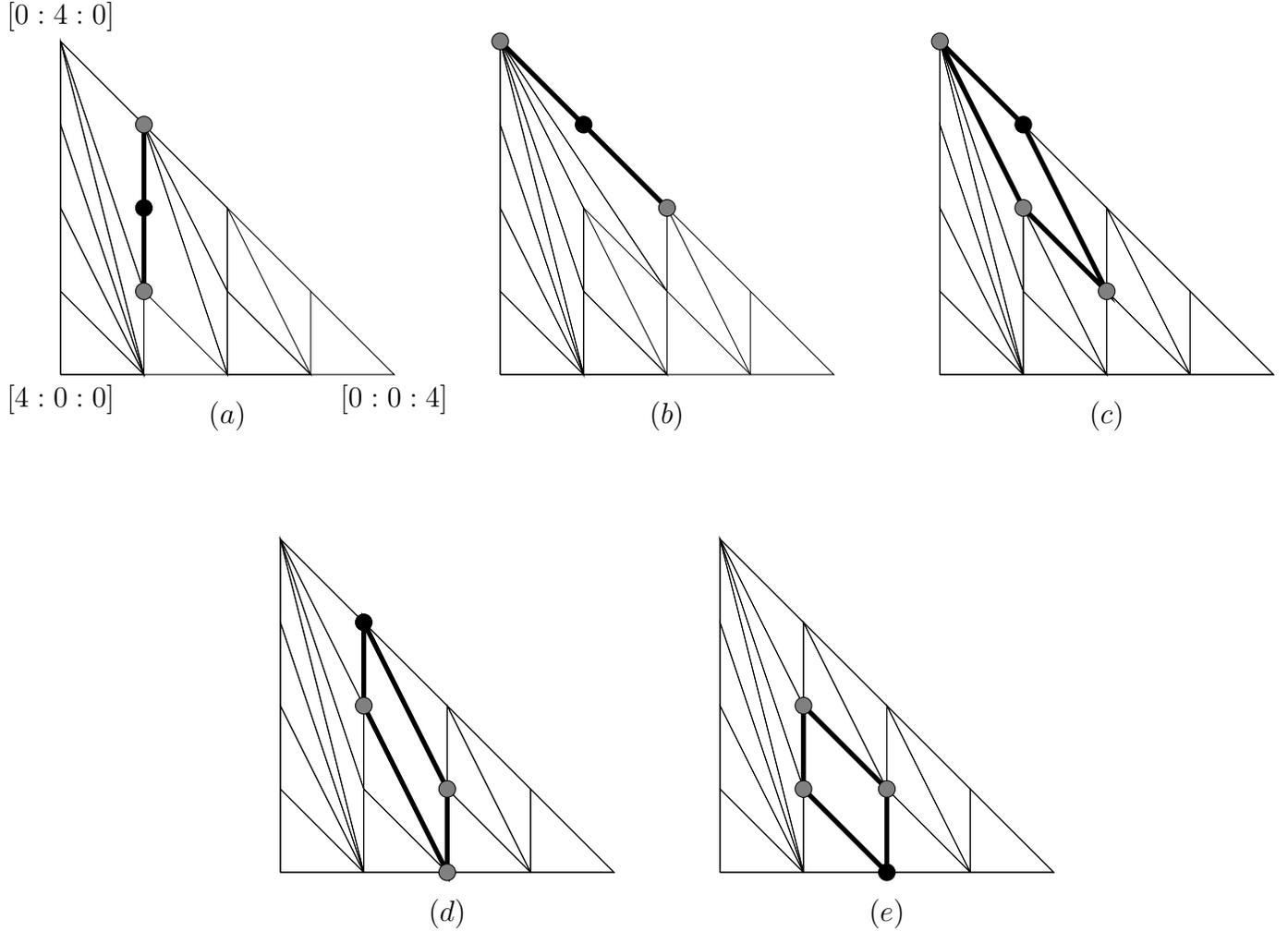
\begin{figure}[h]
	
	\begin{tikzpicture}[scale=1.2]
	
	\begin{scope}
	
	\draw (0,0) -- (0,4) -- (4,0) -- (0,0);
	
	\draw (0, 0) -- (0, 1) -- (1, 0) -- (0, 0);
	\draw (0, 1) -- (0, 2) -- (1, 0) -- (0, 1);
	\draw (0, 2) -- (0, 3) -- (1, 0) -- (0, 2);
	\draw (0, 3) -- (0, 4) -- (1, 0) -- (0, 3);
	\draw (0, 4) -- (1, 0) -- (1, 1) -- (0, 4);
	\draw (0, 4) -- (1, 1) -- (1, 3) -- (0, 4);
	\draw (1, 1) -- (1, 3) -- (2, 0) -- (1, 1);
	\draw (1,3) -- (2,0) -- (2,1) -- (1,3);
	\draw (1,3) -- (2,1) -- (2,2) -- (1,3);
	\draw (2,0) -- (2,1) -- (3,0) -- (2,0);
	\draw (2,1) -- (2,2) -- (3,0) -- (2,1);
	\draw (3,0) -- (3,1);
	
	\draw [line width=2pt] (1,1)--(1,3);
	\filldraw [fill=black] (1,2) circle (0.1);
	\filldraw [fill=gray] (1,1) circle (0.1);
	\filldraw [fill=gray] (1,3) circle (0.1);
	
	\draw (0,0) node [below] {$ [4:0:0] $};
	\draw (4,0) node [below] {$ [0:0:4] $};
	\draw (0,4) node [above] {$ [0:4:0] $};
	\draw (2,-0.5) node {$ (a) $};
	\end{scope}
	
	\begin{scope}[xshift=150]
	
	\draw (0,0) -- (0,4) -- (4,0) -- (0,0);
	
	\draw (0, 0) -- (0, 1) -- (1, 0) -- (0, 0);
	\draw (0, 1) -- (0, 2) -- (1, 0) -- (0, 1);
	\draw (0, 2) -- (0, 3) -- (1, 0) -- (0, 2);
	\draw (0, 3) -- (0, 4) -- (1, 0) -- (0, 3);
	\draw (0, 4) -- (1, 0) -- (1, 1) -- (0, 4);
	\draw (0, 4) -- (1, 1) -- (1, 2) -- (0, 4);
	\draw (1, 0) -- (1, 1) -- (2, 0) -- (1, 0);
	\draw (1, 1) -- (1, 2) -- (2, 0) -- (1, 1);
	\draw (0, 4) -- (1, 2) -- (2, 1) -- (0, 4);
	\draw (0, 4) -- (2, 1) -- (2, 2) -- (0, 4);
	\draw (2,0) -- (2,1);
	\draw (3,0) -- (2,1);
	\draw (3,0) -- (2,2);
	\draw (3,0) -- (3,1);
	
	\draw [line width=2pt] (2,2)--(0,4);
	
	\filldraw [fill=black] (1,3) circle (0.1);
	\filldraw [fill=gray] (0,4) circle (0.1);
	\filldraw [fill=gray] (2,2) circle (0.1);
	
	\draw (2,-0.5) node {$ (b) $};
	\end{scope}
	
	\begin{scope}[xshift=300]
	
	\draw (0,0) -- (0,4) -- (4,0) -- (0,0);
	\draw (1, 0) -- (0, 1) -- (0, 0) -- (1, 0);
	\draw (1, 0) -- (0, 1) -- (0, 2) -- (1, 0);
	\draw (1, 0) -- (0, 3) -- (0, 2) -- (1, 0);
	\draw (1, 0) -- (0, 3) -- (0, 4) -- (1, 0);
	\draw (1, 0) -- (1, 1) -- (0, 4) -- (1, 0);
	\draw (1, 2) -- (1, 1) -- (0, 4) -- (1, 2);
	\draw (1, 0) -- (2, 0) -- (1, 1) -- (1, 0);
	\draw (1, 2) -- (2, 0) -- (1, 1) -- (1, 2);
	\draw (1, 2) -- (2, 0) -- (2, 1) -- (1, 2);
	\draw (2, 2) -- (1, 3) -- (2, 1) -- (2, 2);
	\draw (3, 0) -- (2, 0) -- (2, 1) -- (3, 0);
	\draw (3, 0) -- (2, 2) -- (2, 1) -- (3, 0);
	\draw (3, 0) -- (3, 1) -- (2, 2) -- (3, 0);
	\draw (3, 0) -- (4, 0) -- (3, 1) -- (3, 0);
	\draw (0,4) -- (1,3) -- (2,1) -- (1,2) -- (0,4);
	
	\draw [line width=2pt] (0,4) -- (1,2) -- (2,1) -- (1,3) -- (0,4);
	
	\filldraw [fill=black] (1,3) circle (0.1);
	\filldraw [fill=gray] (0,4) circle (0.1);
	\filldraw [fill=gray] (2,1) circle (0.1);
	\filldraw [fill=gray] (1,2) circle (0.1);
	
	\draw (2,-0.5) node {$ (c) $};
	\end{scope}
	
	\begin{scope}[xshift=75, yshift=-170]
	
	\draw (0,0) -- (0,4) -- (4,0) -- (0,0);
	
	\draw (1, 0) -- (0, 1) -- (0, 0) -- (1, 0);
	\draw (1, 0) -- (0, 1) -- (0, 2) -- (1, 0);
	\draw (1, 0) -- (0, 3) -- (0, 2) -- (1, 0);
	\draw (1, 0) -- (0, 3) -- (0, 4) -- (1, 0);
	\draw (1, 0) -- (1, 1) -- (0, 4) -- (1, 0);
	\draw (1, 2) -- (1, 1) -- (0, 4) -- (1, 2);
	\draw (1, 0) -- (2, 0) -- (1, 1) -- (1, 0);
	\draw (1, 2) -- (2, 0) -- (1, 1) -- (1, 2);
	\draw (1, 2) -- (1, 3) -- (0, 4) -- (1, 2);
	\draw (2, 2) -- (1, 3) -- (2, 1) -- (2, 2);
	\draw (3, 0) -- (2, 0) -- (2, 1) -- (3, 0);
	\draw (3, 0) -- (2, 2) -- (2, 1) -- (3, 0);
	\draw (3, 0) -- (3, 1) -- (2, 2) -- (3, 0);
	\draw (3, 0) -- (4, 0) -- (3, 1) -- (3, 0);
	
	\draw [line width=2pt] (1,2) -- (1,3) -- (2,1) -- (2,0)--(1,2);
	
	\filldraw [fill=black] (1,3) circle (0.1);
	\filldraw [fill=gray] (1,2) circle (0.1);
	\filldraw [fill=gray] (2,0) circle (0.1);
	\filldraw [fill=gray] (2,1) circle (0.1);
	
	\draw (2,-0.5) node {$ (d) $};
	\end{scope}
	
	\begin{scope}[xshift=225, yshift=-170]
	
	\draw (0,0) -- (0,4) -- (4,0) -- (0,0);
	\draw (1, 0) -- (0, 1) -- (0, 0) -- (1, 0);
	\draw (1, 0) -- (0, 1) -- (0, 2) -- (1, 0);
	\draw (1, 0) -- (0, 3) -- (0, 2) -- (1, 0);
	\draw (1, 0) -- (0, 3) -- (0, 4) -- (1, 0);
	\draw (1, 0) -- (1, 1) -- (0, 4) -- (1, 0);
	\draw (1, 2) -- (1, 1) -- (0, 4) -- (1, 2);
	\draw (1, 2) -- (1, 3) -- (0, 4) -- (1, 2);
	\draw (1, 0) -- (2, 0) -- (1, 1) -- (1, 0);
	\draw (1, 2) -- (2, 1) -- (1, 3) -- (1, 2);
	\draw (2, 2) -- (1, 3) -- (2, 1) -- (2, 2);
	\draw (3, 0) -- (2, 0) -- (2, 1) -- (3, 0);
	\draw (3, 0) -- (2, 2) -- (2, 1) -- (3, 0);
	\draw (3, 0) -- (3, 1) -- (2, 2) -- (3, 0);
	\draw (3, 0) -- (4, 0) -- (3, 1) -- (3, 0);
	
	\draw [line width=2pt] (1,1)--(1,2)--(2,1)--(2,0)--(1,1);
	
	\filldraw [fill=black] (2,0) circle (0.1);
	\filldraw [fill=gray] (1,1) circle (0.1);
	\filldraw [fill=gray] (1,2) circle (0.1);
	\filldraw [fill=gray] (2,1) circle (0.1);
	
	\draw (2,-0.5) node {$ (e) $};
	\end{scope}
	
	\end{tikzpicture}
	\caption{Some elementary Mikhalkin circuits in dimension $ n=2 $ and degree $ d=4 $.
		(a) Type I circuit with index 0; (b) Type I circuit with index 1; (c) Type II circuit (index 0); (d) Type III circuit; (e) Type IV circuit.
		In each example, the black point is $ \bs{w} $ and the gray points are elements of $ C_1 $.
	The coordinates are the same as in Figure \ref{fig:Mikhalkin_triangulation}. The coordinates of the vertices of $ \Delta $ are indicated in $ (a) $.}
	\label{fig:elementary_mikhalkin_circuits_dim2}
\end{figure}

\begin{Def}\label{def:mikhalkin_triangulation}
	A subdivision of $ Q $ as in Lemma \ref{lemma:mikhalkin_triangulation} is called \emph{Mikhalkin triangulation of $ Q $}, see Figure \ref{fig:Mikhalkin_triangulation}.
\end{Def}

By \cite[Example 3.5]{Shustin2017_enum}, given $ \bs{v}\in \Delta\cap\mbb{Z}^{n+1} $ and $ {\nu_0:\conv(\{\bs{u}\prec \bs{v}\; | \; \bs{u}\in \Delta\cap\mbb{Z}^{n+1}\})\to \mbb{R}} $ satisfying Mikhalkin condition (away some $ \bs{w}\ne \bs{v} $) we can extend $ \nu_0 $ to a function \linebreak $ {\nu:\conv(\{\bs{u}\preceq \bs{v}\; | \; \bs{u}\in \Delta\cap\mbb{Z}^{n+1}  \}) \to \mbb{R}} $ satisfying Mikhalkin condition away from $ \bs{w} $.

\begin{Example}\label{ex:type_III_circ}
	For indices $ 0<r<\ell\le n $, let $ \bs{w}_i=0 $ for all $ i<\ell $ and $ \bs{w}_{\ell-1} > 0 $.
	Let $ \delta $ be a maximal dimensional cell in the Mikhalkin triangulation of 
	\[ F  := \{ \bs{x}\in \Delta \; | \; \bs{x}_{\ell}=\bs{w}_{\ell}+1, (\forall i>\ell)\bs{x}_i=\bs{w}_i, \bs{x}_{0}=\bs{x}_{1}=\dots=\bs{x}_{r-1}=0  \}  \]
	s.t. $ \bs{x}_{\ell-1}<\bs{w}_{\ell-1}-2 $ for all $ \bs{x}\in\delta $, except possibly the last vertex (w.r.t. $ \prec $).
	If $ {\bs{e}^{(\ell-1)}-\bs{e}^{(\ell-2)}\notin \Sp(\rho)} $ for any proper face $ \rho $ of $ \delta $, the graded circuit $ C $ of height 1 centered at $ \bs{w} $ with $ C_1=\{ \bs{u}^{(-1)}:=\bs{w}-\bs{e}^{(\ell-1)}+\bs{e}^{(\ell)}\}\cup \delta $ admits Mikhalkin condition.
	
	Indeed, we can construct $ \nu:\Delta\to\mbb{R} $ admitted by $ C $ and satisfying Mikhalkin condition which we will do in steps.
	Denote the last vertex of $ \delta $ by $ \bs{u}^{(\ell-r)} $.
	First, we will define $ \nu|_{\conv(\{\bs{v}\prec \bs{w} \}\cup \{ \bs{w}\prec \bs{v} \prec \bs{u}^{(\ell-r)} \}) } $ so that it would satisfy Mikhalkin condition, which can be done by successive applications of \cite[Example 3.5]{Shustin2017_enum}.
	Next, note that since $ \bs{x}^{(\ell-1)}<\bs{w}^{(\ell-1)}-2 $ for all $ \bs{x}\in \delta\setminus \{ \bs{u}^{(\ell-r)} \} $, 
	\[ \conv(C_1\cup C_0)\cap \conv(\{\bs{v}\prec \bs{w} \}\cup \{ \bs{w}\prec \bs{v} \prec \bs{u}^{(\ell-r)} \}) = \conv(\{ \bs{u}^{(-1)} \}\cup\delta\setminus\{\bs{u}^{(\ell-r)}\}) \]
	which is a codimension 1 face of $ \conv(C_0\cup C_1) $, so we can extend $ \nu $ to $ \conv\left \{ \bs{v}\preceq \bs{u}^{(\ell-r)} \right \} $ by \cite[Lemma 3.4]{Shustin2017_enum}.
	Finally, we can apply \cite[Example 3.5]{Shustin2017_enum} successively as in Lemma \ref{lemma:mikhalkin_triangulation}, to extend $ \nu $ to the whole $ \Delta $.
	
	We will call such a graded circuit \emph{elementary Mikhalkin circuit of type III, lower index $ r $ and upper index $ l $}, see Figure \ref{fig:elementary_mikhalkin_circuits_dim2} (d) for a two dimensional example.
\end{Example}

\begin{figure}[h]
	
	\begin{tikzpicture}
	
	\draw (0,0) -- (4,0) -- (0,4) -- (0,0);
	\draw (1, 0) -- (0, 1) -- (0, 0) -- (1, 0);
	\draw (1, 0) -- (0, 1) -- (0, 2) -- (1, 0);
	\draw (1, 0) -- (0, 3) -- (0, 2) -- (1, 0);
	\draw (1, 0) -- (0, 3) -- (0, 4) -- (1, 0);
	\draw (1, 0) -- (1, 1) -- (0, 4) -- (1, 0);
	\draw (1, 2) -- (1, 1) -- (0, 4) -- (1, 2);
	\draw (1, 0) -- (2, 0) -- (1, 1) -- (1, 0);
	\draw (1, 2) -- (2, 0) -- (1, 1) -- (1, 2);
	\draw (1,2)--(2,1) --(0,4);
	\draw ( 2,0) -- (2,1);
	\draw [line width=2pt] (2,1) -- node [anchor=west] {$ \delta $} (2,2);
	
	\filldraw [fill=black] (1,3) circle (0.1) node [above] {$ \bs{w} $};
	\filldraw [fill=gray] (1,2) circle (0.1);

	\draw (0,0) node [below] {$ [4:0:0] $};
	\draw (4,0) node [anchor=west] {$ [0:0:4] $};
	\draw (0,4) node [above] {$ [0:4:0] $};
	
	\draw [->] (-2,1) -- (-2,2) node [above] {$ x_1 $ };
	\draw [->] (-2,1) -- (-1,1) node [anchor=west] {$ x_2 $ };
	\draw [->] (-2,1) -- (-2.7,0.3) node [anchor=east] {$ x_0 $ };
	\end{tikzpicture}
	
	\caption{The resulting subdivision of beginning of the process described in Example \ref{ex:type_III_circ} with $ \delta=\{ [1:1:2], [0:2:2] \} $.
	The black circle is $ \bs{w} $, the gray circle is $ \b{u}^{(-1)} $ and the bold line is $ \conv(\delta) $. The coordinates are as in Figure \ref{fig:Mikhalkin_triangulation}.}
	\label{fig:type_III_obstruction}
\end{figure}
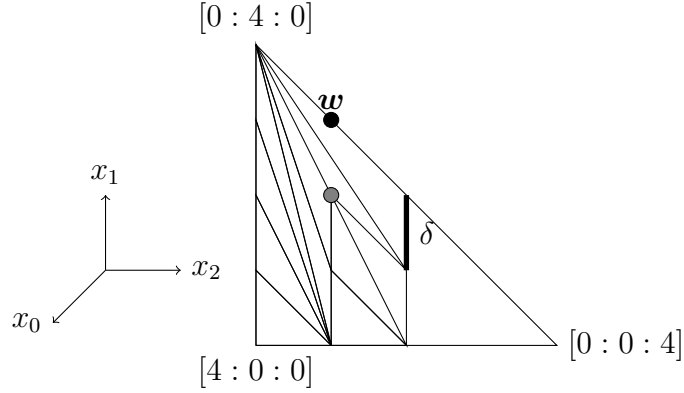

\begin{remark}
	The condition $ \bs{x}_{\ell-1}<\bs{w}_{\ell-1}-2  $ can not be omitted entirely.
	Indeed, consider an elementary Mikhalkin circuit of type III in dimension $ n=2 $ with lower index $ r=0 $ and upper index $ l=2 $.
	Then, if we perform the beginning of the construction of $ \nu $ as in Example \ref{ex:type_III_circ}, we get the subdivision depicted (for $ d=4 $) in Figure \ref{fig:type_III_obstruction} for which we can not fit $ \conv(\{\bs{u}^{(-1)}, \bs{w}\}\cup \delta) $.
	
	It is possible that the condition $ \bs{x}_{\ell-1}<\bs{w}_{\ell-1}-2 $ we imposed for the construction is not necessary and a weaker condition can be found.
	However, asymptotically as $ d\to \infty $ almost all of the circuits that admit Mikhalkin condition satisfy it (see Proposition \ref{prop:discriminant_degree}), so it wouldn't affect our count.
	
	The situation where $ \delta $ is of maximal dimension is also not the most general example of a circuit constructed this way and admitting Mikhalkin condition.
	For an example of a more general construction, consider $ \bs{w}=[0:0:\bs{w}_2:\bs{w}_3] $ with $ \bs{w}_2>0 $ and 
	\[ C_1=\{ \bs{u}^{(-1)}:=[0:1:\bs{w}_2-1:\bs{w}_3], \bs{u}^{(1)}:=[a:1:b:\bs{w}_3+1], \bs{u}^{(2)}:=[a:0:b+1:\bs{w}_3+1] \} \]
	for some $ a,b\in\mbb{Z}, a>0, b\ge 0 $ and $ a+b=\bs{w}_2-2 $, then $ C $ is a circuit admitting Mikhalkin condition (see \cite[Lemma 3.9]{Shustin2017_enum}).
	As before, the fact we are ignoring those circuits wouldn't affect the asymptotic of our result.
\end{remark}

\begin{Example}\label{ex:type_IV_circ}
	Again, pick two indices $ 0\le r,\ell \le n $ with $ r+1<\ell $, let $ \bs{w}_i=0 $ for all $ i<\ell $ except $ i=r $ and $ \bs{w}_r,\bs{w}_{\ell}>0 $.
	Let $ \bs{u}\in \Delta\cap\mbb{Z}^{n+1} $ with $ \bs{u}_{\ell}=\bs{w}_{\ell}-1 $, $ \bs{u}_j=0 $ for all $ j<r $, $ \bs{u}_r>0 $ and $ \bs{u}_i=\bs{w}_i $ for all $ i>\ell $.
	The graded circuit $ C $ of height 1 centered at $ \bs{w} $ with $ C_1=\{ \bs{u}^{(-2)}:=\bs{u}, \bs{u}^{(-1)}:=\bs{u}+\bs{e}^{(r+1)}-\bs{e}^{(r)}, \bs{u}^{(1)}:=\bs{w}+e^{(r+1)}-\bs{e}^{(r)}  \}  $ admits Mikhalkin condition by Lemma \ref{lemma:graded_circ_admit_mikhalkin_comb_cond}.
	Indeed, let $ \bs{v}\in \Sp(C)=\Sp(C_1) $, in particular, $ \bs{v}_i=\bs{w}_i $ for all $ i>\ell $.
	The point $ \bs{u}^{(1)} $ is the first (w.r.t. $ \prec $) point of $ \Sp(C_1) $ with $ \bs{x}_{\ell}\ge \bs{w}_{\ell} $ meaning that if $ \bs{v}_{\ell}\ge \bs{w}_{\ell} $ then the last coefficient of \eqref{eq:lemma_mikh_GC_comb_cond} is positive by Remark \ref{remark:v_last_in_comb_cond}.
	If $ \bs{v}_{\ell}<\bs{w}_{\ell}-1 $ then, by computing the $ \ell $'th coordinate of \eqref{eq:lemma_mikh_GC_comb_cond} we get that the coefficient of $ \bs{u}^{(1)} $ is positive and obviously it is the last point.
	Finally, if $ \bs{v}_{\ell}=\bs{w}_{\ell}-1 $ then the coefficient of $ \bs{u}^{(1)} $ is 0 meaning that 
	\[ \bs{v} - (1-\alpha) \bs{u} - \alpha\cdot \bs{u}^{(-1)} = 0 \]
	where $ \alpha<0 $ whenever $ \bs{v}\prec \bs{u} $ which finishes the proof that $ C $ admits Mikhalkin condition.
	We will call such a graded circuit \emph{elementary Mikhalkin circuit of type IV, lower index $ r $ and upper index $ \ell $}.
	In Figure \ref{fig:elementary_mikhalkin_circuits_dim2} (e), one can find an example in dimension $ n=2 $.
	The linear combination \eqref{eq:def_GC_unique_linear_comb} is $ \bs{u}^{(-2)}+\bs{u}^{(1)}-\bs{u}^{(-1)} $.
\end{Example}

\begin{remark}\label{remark:C_1_is_regular}
	Note that in Examples \ref{ex:type_II_circ}-\ref{ex:type_IV_circ}, $ C_1 $ is a regular simplex (i.e. \linebreak $ {\Sp_\mbb{Z}(C_1)=\mbb{Z}^{n+1}\cap \Sp_\mbb{R}(C_1)} $).
	This means that there exist a linear combination with integer coefficients $$ \bs{w}+\sum_{\bs{v}\in C_1}m_{\bs{v}}\bs{v}=0 $$.
\end{remark}

\begin{Example}\label{ex:all_elementary_circuits_dim3}
	We consider examples of elementary Mikhalkin circuits in dimension $ n=3 $.
	For type I, the index can be $ 0,1 $ or $ 2 $, the latter can not produce a singularity while the 2 former can.
	For circuit of type II, the index can be $ 1 $, which can correspond to a singularity, or $ 2 $, which can not.
	
	For type III, the pair of indices can take the following values: $ (0,2),(1,3),(0,3) $.
	The first give a parallelogram which can correspond to a singularity.
	The second give a parallelogram in the boundary of $ \Delta $ and thus can not give a singularity \cite[Theorem 2]{Markwig2012_trop_surf_sings}.
	Next, an elementary Mikhalkin circuit of type III, lower index 0 and upper index 3 is \linebreak $ {\{\bs{w}-\bs{e}^{(2)}+\bs{e}^{(1)},\bs{w}\}\cup \delta} $ where $ \delta $ is either
	\[ \left \{ [\bs{v}_0:\bs{v}_1:\bs{v}_2:\bs{w}_3+1], [\bs{v}_0-1:\bs{v}_1+1:\bs{v}_2:\bs{w}_3-1], [\bs{v}_0+\bs{v}_1-1:0:\bs{v}_2+1:\bs{w}_3+1] \right \} \] 
	or
	\[ \left \{ [0:\bs{v}_0+\bs{v}_1+1:\bs{v}_2-1:\bs{w}_3+1], [\bs{v}_0:\bs{v}_1:\bs{v}_2:\bs{w}_3+1], [\bs{v}_0-1:\bs{v}_1+1:\bs{v}_2:\bs{w}_3+1]  \right \}, \]
	see Figure \ref{fig:Mikhalkin_triangulation}, where in both those examples $ \bs{w}=[0:0:\bs{w}_2:\bs{w}_3] $, $ \bs{v}_0+\bs{v}_1+\bs{v}_2=\bs{w}_2-1 $ and in the first example $ \bs{v}_1>1 $ while in the second $ \bs{v}_0>1 $ (see \cite[Lemma 5.5]{Shustin2017_enum}).
	
	Finally, a type IV circuit of lower index 1 (and upper index 3) is a parallelogram contained in $ \partial\Delta $ meaning that it does not correspond to a singularity.
	We are left with type IV circuit of lower index 0 and upper index either 2 or 3, both of which are parallelograms in the interior of $ \Delta $ and thus correspond to a singularity. 
\end{Example}

\newpage 
\subsection{Gluing graded circuits}\label{subsec:comb4}

\begin{Def}
	Let $ C $ be a graded circuit centered at $ \bs{w} $. A linear surjective map  
	\[ \pr:\mbb{Z}^{n+1} \to \mbb{Z}^{k+1} \]
	is called \emph{a projection along of $ C $} if $ \Sp(C)\cap\Delta=\pr^{-1}(\bs{w}') $ for some $ \bs{w}'\in \mbb{Z}^{k+1} $. 
\end{Def}

\begin{Def}
	Let $ C $ be a graded circuit centered at $ \bs{w} $ and $ \pr:\mbb{Z}^{n+1}\to\mbb{Z}^{k+1} $ a projection along $ C $. A function $ \sigma:\pr(\Delta^{(n)}\cap \mbb{Z}^{n+1})\setminus\{\pr(\bs{w}) \} \to \Delta^{(n)} $ is called \emph{min-section of $ \pr $ (or of $ C $)} if $ \pr\circ \sigma = Id $ and for every $ \nu:\Delta^{(n)}\to\mbb{R} $ satisfying Mikhalkin condition that $ C $ admits, 
	\[ \nu(\sigma(\bs{y})) = \min \{ \nu(\bs{x}) \; | \; \pr(\bs{x})=\bs{y} \} \]
	for all $ \bs{y}\in \pr(\Delta^{(n)}) \setminus \{ \pr(\bs{w})\} $.
\end{Def}

\begin{lemma}
	If $ C $ admits Mikhalkin condition it has a unique min-section.
\end{lemma}

\begin{proof}
	Let $ \nu $ be a function satisfying Mikhalkin condition admitted by $ C $, $ \pr $ a projection along $ C $ and $ \bs{y}\in \pr(\Delta\cap\mbb{Z}^{n+1}) $. By Lemma \ref{lemma:difference_set_no_parallel_circuit}, for any distinct $ \bs{x}^{(1)},\bs{x}^{(2)}\in \pr^{-1}(\bs{y})  $ there exist a linear combination
	\[ \bs{x}^{(2)} - \bs{x}^{(1)} - \sum_{i=1}^{ht(C)} \sum_{\bs{x}\in C_i} \alpha_{\bs{x}} \bs{x} = 0 \]
	with $ \sum_{\bs{x}\in C_i}\alpha_{\bs{x}} = 0 $ for all $ i>0 $. Since not all the coefficients in this dependency are 0 and $ \nu $ satisfies Mikhalkin condition, 
	\[ \nu(\bs{x}^{(2)}) - \nu(\bs{x}^{(1)}) = \nu(\bs{x}^{(2)}) - \nu(\bs{x}^{(1)}) - \sum_{i=1}^{ht(C)} \sum_{\bs{x}\in C_i} \alpha_{\bs{x}} \nu(\bs{x}) \ne 0 \]
	Meaning that the minimum of $ \nu|_{\pr^{-1}(\bs{y})} $ occur only once, so the min-section is unique.
	
	To prove existence, note that the sign of $ \nu(\bs{x}^{(2)})-\nu(\bs{x}^{(1)}) $ depends only on the coefficients $ \alpha_{\bs{x}} $, meaning that the sign would be the same for every function satisfying Mikhalkin condition and admitted by $ C $. Thus the minimum would be attained at the same point for every such function.
\end{proof}

\begin{Construction}\label{con:graded_circuit_glue}
	Let $ \widetilde{C} $ be a graded circuit centered at $ \bs{w} $ and admitting Mikhalkin condition, let $ \pr:\mbb{Z}^{n+1}\to\mbb{Z}^{k+1} $ and $ \sigma:\pr(\Delta^{(n)}\cap\mbb{Z}^{n+1})\to \Delta^{(n)} $ be the corresponding projection and min-section.
	Let $ C'\subseteq \pr(\Delta^{(n)}\cap\mbb{Z}^{n+1}) $ be a graded circuit centered at $ \pr(\bs{w}) $.
	We denote by $ \widetilde{C}\Vert C' $ the graded circuit $ C\subseteq \Delta^{(n)} $ that have $ C_i=\widetilde{C}_i $ for all $ i\le ht(\widetilde{C} ) $ and $ C_{i}=\sigma\left(C'_{i - ht(\widetilde{C})}\right) $ for $ i>ht(\widetilde{C}) $.
	We will call such a graded circuit the \emph{gluing} of $ \widetilde{C} $ and $ C' $.
\end{Construction}

We will show that $ \widetilde{C}\Vert C' $ admits Mikhalkin condition whenever $ \widetilde{C} $ and $ C' $ satisfy certain restrictions. 
The requirement from $ C' $ encapsulated in the following definition.

\begin{Def}\label{def:elementary_Mikhalkin_graded_circuit}
	We will call a graded circuit $ C $ \emph{elementary graded Mikhalkin circuit} if it is one of the following:
	\begin{enumerate}
		\item 
		Elementary Mikhalkin circuit (see Examples \ref{ex:type_I_circ}-\ref{ex:type_III_circ} ) of types I, II or III.
		
		\item
		A graded circuit of any height, with $ C_1 $ being elementary Mikhalkin circuit of type IV, lower index $ r $ and upper index $ \ell $ (Example \ref{ex:type_IV_circ}), and $$ {\Sp(C)=\Sp\{ \bs{w}, \bs{e}^{(i)}-\bs{e}^{(j)} \; | \; r\le i,j \le \ell \}}. $$
	\end{enumerate}
\end{Def}

\begin{lemma}\label{lemma:parrallel_projection}
	Let $ \widetilde{C} $ be a graded circuit admitting Mikhalkin condition and let 
	\begin{align*}
	\pr : \Delta^{(n)} &\to \Delta^{(n-\ell_2+\ell_1)} \\ 
	[\bs{x}_0:\dots:\bs{x}_n] &\mapsto [\bs{x}_0:\dots:\bs{x}_{\ell_1-1}:\sum_{i=\ell_1}^{\ell_2}\bs{x}_i:\bs{x}_{\ell_2+1}:\dots:\bs{x}_n]
	\end{align*} 
	be the projection along $ \widetilde{C} $ (for some $ 0\le \ell_1< \ell_2 \le n $). 
	Suppose additionally, that $ \bs{w}_i=0 $ for all $ i<\ell_1 $.
	Let $ C' \subseteq \Delta^{(n-\ell_2+\ell_1)} $ be an elementary Mikhalkin graded circuit (see Definition \ref{def:elementary_Mikhalkin_graded_circuit}) with $ \Sp(C')=\mbb{R}^{n-\ell_2+\ell_1+1} $ .
 	Then $ \widetilde{C}\Vert C' $ admits Mikhalkin condition.
\end{lemma}

\begin{proof}
	Denote by $ \sigma $ the min-section of $ \widetilde{C} $.
	First, we will show that $ C $ has no parallel levels. 
	The first $ ht(\widetilde{C}) $ levels are not parallel to previous ones since they are the levels of $ \widetilde{C} $.
	For $ i>ht(\widetilde{C}) $, if 
	\[ \sum_{\bs{x}\in C_{i}}\alpha_{\bs{x}} \bs{x}\in \Sp(\overline{C}_{i-1}) \] 
	then, by projecting along $ \widetilde{C} $, we get that 
	\[ \sum_{\bs{x}\in C_{i-ht(\widetilde{C})}'}\alpha_{\bs{x}} \bs{x} \in \Sp(\overline{C}_{i-ht(\widetilde{C})-1}), \] 
	meaning that $ \sum \alpha_{\bs{x}}\ne 0 $, since $ C' $ has no parallel levels.

	Assume, towards a contradiction, that $ C $ does not admits Mikhalkin condition. Then, by Lemma \ref{lemma:graded_circ_admit_mikhalkin_comb_cond}, there exist $ C_h \in C $, $ \bs{v}\in \Sp(\overline{C}_h)\cap\Delta\cap\mbb{Z}^{n+1}\setminus C_h $ and a linear combination
	\[ \bs{v}-\sum_{0< i \le h}\sum_{\bs{x}\in C_i} \alpha_{\bs{x}} \bs{x} =0 \]
	with $ \sum_{\bs{x}\in C_i}\alpha_{\bs{x}} =0 $ for all $ i<h $ and $ \alpha_{\bs{u}_{\max}}<0 $ where $ \bs{u}_{\max} $ is the last (w.r.t. $ \prec $) vector with $ \alpha_{\bs{u}_{\max}}\ne 0 $.
	Projecting along $ \widetilde{C} $ we get 
	\begin{equation}\label{eq:lemma_proj_parallel_decomp_of_pr_v}
	\pr(\bs{v})=\sum_{i=1}^{h-ht(\widetilde{C})}\sum_{\bs{y}\in C_i'}\alpha_{\sigma(\bs{y})}\bs{y}. 
	\end{equation}
	Since $ C' $ admits Mikhalkin condition, we are left to show that the maximal vector appearing in \eqref{eq:lemma_proj_parallel_decomp_of_pr_v} is $ pr(\bs{u}_{\max}) $.
	We will do it separately for each type of $ C'_1 $:
	
	\begin{description}
		\item[Type I] 
		Since we assumed $ \Sp(C')=\mbb{R}^{(n-\ell_2+\ell_1+1)} $, $ \ell_1=0 $ and $ \ell_2=n-1 $ meaning that $ n-\ell_2+\ell_1=1 $.
		If $ \bs{y},\bs{z}\in \Delta^{(1)}\cap\mbb{Z}^{2} $ with $ \bs{y}\prec \bs{z} $ then $ \bs{y}_1 < \bs{z}_1 $ meaning that 
		\[ \sigma(\bs{y})_n=\bs{y}_1<\bs{z}_1=\sigma(\bs{z})_n \]
		so $ \sigma(\bs{y})\prec \sigma(\bs{z}) $ and thus $ \sigma $ is monotone.
		Now, \eqref{eq:lemma_proj_parallel_decomp_of_pr_v} reads as $ pr(\bs{v})=\alpha_1\bs{u}^{(1)}+\alpha_{-1}\bs{u}^{(-1)} $ where $ \bs{u}^{(-1)}:=\bs{w}'-\bs{e}^{(1)}+\bs{e}^{(0)} $ and $ \bs{u}^{(1)}:=\bs{w}'+\bs{e}^{(1)}-\bs{e}^{(0)} $ and the conclusion can be obtained by considering the options $ \pr(\bs{v})\succ \bs{u}^{(1)} $ (in which $ \bs{u}_{\max}=\bs{v} $) and $ \bs{v}\prec \bs{u}^{(-1)} $ ($ \pr(\bs{u}_{\max})=\bs{u}^{(1)} $).
		
		\item [Type II] 
		In this case, since $ \Sp(C')=\mbb{R}^{(n-\ell_2+\ell_1+1)} $, we have $ \ell_1=1 $, $ \ell_2=n-1 $ and $ C'_1=\{\bs{u}^{(-2)}:=\bs{w}'-\bs{e}^{(2)}+\bs{e}^{(1)}, \bs{u}^{(-1)}:=\bs{w}'-\bs{e}^{(1)}+\bs{e}^{(0)}, \bs{u}^{(1)}:=\bs{w}+\bs{e}^{(2)}-2\bs{e}^{(1)}+\bs{e}^{(0)} \} $.
		Note that for every $ \bs{v}'\in \Delta^{(2)}\cap\mbb{Z}^3 $ with $ \bs{v}'_2>\bs{w}'_2=\bs{w}_n $, we have 
		\[ \sigma(\bs{v}')=[\bs{v}'_0:\bs{v}'_1:0:\dots:0:\bs{v}'_2] \]
		meaning that $ \sigma $ is monotone for vectors succeeding $ \bs{w}' $.
		Thus, if $ \pr(\bs{v})\succ \bs{u}^{(1)} $ then $ \bs{u}_{\max}=\bs{v}\succeq\sigma(\pr(\bs{v})) $ and if $ \pr(\bs{v}) \prec \bs{u}^{(1)} $ and $ \bs{u}^{(1)} $ appear in $ \eqref{eq:lemma_proj_parallel_decomp_of_pr_v} $, then $ pr(\bs{u}_{\max})=\bs{u}^{(1)} $.
		If $ \bs{u}^{(1)} $ does not appear in \eqref{eq:lemma_proj_parallel_decomp_of_pr_v} then $ \pr(\bs{v})\in \Sp(\bs{u}^{(-2)}, \bs{u}^{(-1)}) $ and thus either $ \pr(\bs{v})\in \{ \bs{u}^{(-1)},\bs{u}^{(-2)} \} $ or $ \bs{v}_n>\bs{w}_n $ meaning that $ \bs{v} $ is the last vector appearing in \eqref{eq:lemma_proj_parallel_decomp_of_pr_v} and $ \bs{v}=\bs{u}_{\max} $.
		
		\item [Type III]
		Since $ \Sp(C')=\mbb{R}^{n-\ell_2+\ell_1+1} $, $ \ell_2=n-1 $. 
		By the same reason as in the case of type I, $ \sigma $ is monotone on $ \{ \bs{x}\in \Delta^{(n-\ell_2+\ell_1)}\cap\mbb{Z}^{n-\ell_2+\ell_1+1} \; | \; \bs{x}_{n-\ell_2+\ell_1}>\bs{w}_n \} $ meaning that if $ (\bs{u}_{\max})_n>\bs{w}_n $, the maximal vector appearing in \eqref{eq:lemma_proj_parallel_decomp_of_pr_v} is $ pr(\bs{u}_{\max}) $.
		Otherwise, since the only vector of $ C'_1 $ with $ \bs{x}_n\le \bs{w}_n $ is $ \bs{u}^{(-1)} $ (see the notation of Example \ref{ex:type_III_circ}), $ \bs{v}=\sigma(\bs{u}^{(-1)}) $ and there is nothing to prove. 	
		
		\item [Type IV] 
		Similar to the case of type I, $ \ell_1=0 $ and $ \sigma $ is monotone meaning that if the last vector in \eqref{eq:lemma_proj_parallel_decomp_of_pr_v} succeeds $ \bs{w}' $ we are done.
		In the case where all the summands of \eqref{eq:lemma_proj_parallel_decomp_of_pr_v} precedes $ \bs{w}' $ it is sufficient to show that $ \alpha_{\bs{x}}=0 $ for all $ \bs{x}\in\bigcup \widetilde{C} $, which we will do by showing that $ \sigma $ is linear in the region $ \{ \bs{y}\prec \bs{w}' \} $.
		More precisely, we will find an index $ 0\le j \le \ell_2 $ s.t. $\sigma(\bs{y})_i=0 $ for all $ j\ne i\le \ell_2 $ and $ \bs{y}\prec \bs{w}' $.
		
		For every $ 0\le i \le \ell_2 $ write $ \bs{e}^{(i)} = \sum_{\bs{x}\in \bigcup \widetilde{C} }\beta_{i,\bs{x}}\bs{x} $.
		We pick $ j $ to be such that for every $ j \ne i\le \ell_2 $ there exists $ \bs{x}^{(i)}\in \bigcup\widetilde{C} $ s.t. $ \beta_{j,\bs{x}^{(i)}} > \beta_{i,\bs{x}^{(i)}} $ and for all $ \bs{z}\succ \bs{x}^{(i)} $, $ \beta_{j,\bs{z}} = \beta_{i,\bs{z}} $.
		Indeed, with this choice, if $ \sigma(\bs{y})_i> 0  $ for $ j\ne i\le \ell_2 $ and $ \bs{y}\prec \bs{w}' $, then 
		\[ \bs{e}^{(j)}-\bs{e}^{(i)} = \sum_{\bs{x}\in \widetilde{C}}\left( \beta_{j,\bs{x}}-\beta_{i,\bs{x}} \right)\bs{x} = (\beta_{j,\bs{x}^{(i)}}-\beta_{i,\bs{x}^{(i)}})\bs{x}^{(i)} + \sum_{\bs{x}^{(i)}\succ \bs{x}\in \widetilde{C}}\left( \beta_{j,\bs{x}}-\beta_{i,\bs{x}} \right)\bs{x} \]
		and since $ \beta_{j,\bs{x}^{(i)}}>\beta_{i,\bs{x}^{(i)}} $, from Mikhalkin condition on $ \nu $,
		\[ \nu(\sigma(\bs{y})-\bs{e}^{(i)}+\bs{e}^{(j)})<\nu(\sigma(\bs{y})) \]
		which is a contradiction.
	\end{description}
\end{proof}

Lemma \ref{lemma:parrallel_projection} allows us to glue graded circuits when $ \Sp(\widetilde{C}) $ is parallel to the axes. 
It is not always the case for type IV circuits.
We will deal with the other option in Lemma \ref{lemma:gluing_type_IV}, but first we state the restrictions we impose on $ \widetilde{C} $.

\begin{Def}\label{def:prefix_type_IV}
	Let $ \widetilde{C}\subseteq \Delta^{(n)} $ be a graded circuit centered at $ \bs{w}\in\Delta\cap\mbb{Z}^{n+1} $ with $ \widetilde{C}_1 \cup \widetilde{C}_0 $ being an elementary Mikhalkin graded circuit of type IV, lower index $ 0 $ and upper index $ n $, see Example \ref{ex:type_IV_circ}.
	We call $ \widetilde{C} $ a \emph{dimension $ \ell $ prefix of a type IV elementary Mikhalkin graded circuit} if 
	\begin{enumerate}
		\item
		$ \widetilde{C} $ admits Mikhalkin condition.
		
		\item
		$ \Sp(\widetilde{C}) = \Sp\{ \bs{w},\bs{u},\bs{e}^{(i)}-\bs{e}^{(j)} \; | \; i,j\le \ell \} $. 
		
		\item
		All the elements of $ \bigcup \widetilde{C} $ except $ \bs{w} $ and $ \bs{w}-\bs{e}^{(1)}+\bs{e}^{(0)} $ have $ n $'th coordinate equal to $ \bs{w}_n-1 $.
	\end{enumerate}
\end{Def}

\begin{lemma}\label{lemma:gluing_type_IV}
	Let $ \widetilde{C} $ be a dimension $ \ell $ prefix of a type IV elementary Mikhalkin graded circuit, see Definition \ref{def:prefix_type_IV}.
	Let $ \bs{u}\in\widetilde{C}_1 $ be the first point of $ C_1 $ as in Example \ref{ex:type_IV_circ}.
	Consider the projection 
	\begin{align*}
	\pr :\Delta &\to \mbb{Z}^{n-\ell}\cap\{\sum \bs{x}_i = \bs{w}_0+1\} \\
	[\bs{x}_0:\dots:\bs{x}_n] &\mapsto [\sum_{i=0}^\ell \bs{x}_i : \bs{x}_{\ell+1}: \dots: \bs{x}_{n-1}] + (\bs{x}_n-\bs{w}_n+1)\cdot [\sum_{i=0}^{\ell} \bs{u}_i-\bs{w}_0: \bs{u}_{\ell+1}:\dots : \bs{u}_{n-1}]
	\end{align*}
	along $ \widetilde{C} $.
	
	Let $ C'\subseteq \Delta_{\bs{w}_0+1}^{(n-\ell-1)} $ be an elementary Mikhalkin graded circuit (see Definition \ref{def:elementary_Mikhalkin_graded_circuit}) centered at $ \bs{w}'=\pr(\bs{w}) $.
	Then $ C:=\widetilde{C}\Vert C' $ is a prefix of type IV elementary Mikhalkin graded circuit.
\end{lemma}

\begin{proof}
	Since $ \widetilde{C}_1 \cup \widetilde{C}_0 $ is of type IV and lower index $ 0 $, we have $ \bs{u}_0>0 $ and so $ \bs{w}'_0 = \sum_{i=0}^{\ell}\bs{u}_i>0 $.
	Thus $ C'_1\cup C'_0 $ has to be of type I or IV.
	The first option implies that 
	\[ C'_1 = \{ \bs{w}'-\bs{e}^{(1)}+\bs{e}^{(0)}, \bs{w}'+\bs{e}^{(1)}-\bs{e}^{(0)} \} \]
	and $ C' $ has no other levels.
	The latter option is, since we assumed that $ C' $ is an elementary Mikhalkin graded circuit, that $ C'_1 $ is of type IV, lower index $ 0 $ and upper index $ \ell' $,
	\[ \bs{u}_{\ell+1}=\textellipsis=\bs{u}_{\ell+\ell'-1}=0 \]
	and $ \Sp(C')=\Sp\{ \bs{w}', \bs{e}^{(i)}-\bs{e}^{(j)} \; | \; i,j\le \ell' \} $.
	
	Throughout the proof we denote the maximal point of $ \widetilde{C} $ by $ \bs{u}_{\max} $ which, by condition (3) of Definition \ref{def:prefix_type_IV}, is equal to $ \bs{w}-\bs{e}^{(1)}+\bs{e}^{(0)} $.
	We also denote by $ \sigma:\pr(\Delta^{(n)}\cap\mbb{Z}^{n+1})\to \Delta^{(n)} $ the min-section of $ \pr $.
	
	We begin by proving that $ C $ satisfies condition (2) in Definition \ref{def:prefix_type_IV}.
	We have 
	\[ \Sp(C)=\pr^{-1}(\Sp(C'))=
	\Sp\{ \bs{w},\bs{u},\bs{e}^{(i)}-\bs{e}^{(j)} \; | \; i,j\le \ell+\dim C' \}. \]
	
	To show that $ C $ satisfies (3) in Definition \ref{def:prefix_type_IV}, we need to show that $ \sigma(\bs{y})_n=\bs{w}_n-1 $ for all $ \bs{y}\in \bigcup C'\setminus \{\bs{w}'\} $. 
	First, note that 
	\[ [0: \textellipsis  :0: \bs{y}_0 : \bs{y}_1:\textellipsis:\bs{y}_{n-\ell-2}:\bs{y}_{n-\ell-1}:\bs{w}_n-1]\in \pr^{-1}(\bs{y}) \]
	so the fiber of $ \bs{y} $ contains elements with $ \bs{x}_n=\bs{w}_n-1 $. 
	Next, if $ \bs{x}'\in \pr^{-1}(\bs{y})\cap\Delta^{(n)}\cap \mbb{Z}^{n+1} $ with $ \bs{x}'_n\ne  \bs{w}_n -1 $, write
	\begin{equation}\label{eq:lemma_glue_IV_cond_b_proof}
	\bs{x}'-\bs{x}+\sum_{\bs{v}\in \bigcup \widetilde{C}\setminus \{\bs{w}\}}\alpha_{\bs{v}} \bs{v} = 0
	\end{equation}
	with $ \sum_{\bs{v}\in \widetilde{C}_i} \alpha_{\bs{v}} =0 $ for all $ i\le ht(\widetilde{C}) $.
	
	If $ \bs{x}_n\ge \bs{w}_n $ then $ \bs{x}' \succ \bs{u}_{\max} $ so it is the last vector in \eqref{eq:lemma_glue_IV_cond_b_proof}	meaning that $ \nu(\bs{x}')>\nu(\bs{x}) $ for any $ \nu $ satisfying Mikhalkin condition, and so $ \sigma(\bs{y})\ne \bs{x}' $.
	If, on the other side, $ \bs{x}'_n < \bs{w}_n-1 $, computing the $ n $'th coordinate of \eqref{eq:lemma_glue_IV_cond_b_proof} gives us 
	\[ \alpha_{\bs{u}_{\max}}= \alpha_{\bs{u}_{\max}} + \sum_{1\le i\le ht(\widetilde{C})} (\bs{w}_n-1)\sum_{\bs{v}\in \widetilde{C}_i} \alpha_{\bs{v}}= \sum_{\bs{v}\in \bigcup\widetilde{C}\setminus\{\bs{w}\}}\alpha_{\bs{v}} \bs{v}_n=\bs{x}_n-\bs{x}'_n=\bs{w}_n-1-\bs{x}'_n>0, \]
	and since $ \bs{u}_{\max} $ is the last vector appearing in \eqref{eq:lemma_glue_IV_cond_b_proof}, $ \nu(\bs{x}')>\nu(\bs{x}) $ for all $ \nu $ satisfying Mikhalkin condition and so $ \sigma(\bs{y})\ne \bs{x}' $.
	
	We now show that $ C $ admits Mikhalkin condition.
	We argue as in the proof of Lemma \ref{lemma:parrallel_projection}.
	Assume, towards a contradiction, that there exists $ \bs{v}\in \Delta^{(n)} $ s.t. the coefficient of the last vector of 
	\[ \bs{v}+\sum_{\bs{x}\in \bigcup C}\alpha_{\bs{x}} \bs{x}=0 \]
	is negative (as always, $ \sum_{\bs{x}\in C_i}\alpha_{\bs{x}} $ equals 0 for all $ i< ht(\widetilde{C}  ) $).
	Since $ \bs{u}_{\max}=\max_\prec \bigcup C $,  $ \bs{v}\prec \bs{u}_{\max} $ meaning that $ \bs{v}_n\le \bs{w}_n-1 $.
	If $ \bs{v}_n< \bs{w}_n-1 $ then $ \alpha_{\bs{u}_{\max}}>0 $, so $ \bs{v}_n=\bs{w}_n-1 $.
	Similarly to the proof of Lemma \ref{lemma:parrallel_projection}, $ \sigma $ is monotone on $ \Delta_{\bs{w}_0+1}^{n-\ell-1} $ meaning that if $ \pr(\bs{v})\succ \bs{w}' $ we are done.
	If $ \pr(\bs{v})\prec \bs{w}' $ and $ C' $ is of type I then $ \pr(\bs{v})=\bs{w}'-\bs{e}^{(1)}+\bs{e}^{(0)} $ and there is nothing to prove.
	Finally, suppose that $ C'_1 $ is of type IV and $ \pr(\bs{v})\prec \bs{w}' $.
	As in the corresponding case of Lemma \ref{lemma:parrallel_projection}, there exists $ 0\le j\le \ell $ s.t.  $ \sigma(\bs{y})_i=0 $ for all $ j\ne i\le \ell $ and $ \bs{y}\prec \bs{w}' $ meaning that $ \alpha_{\bs{x}}=0 $ for all $ \bs{x}\in \bigcup \widetilde{C} $ which finishes the proof.
\end{proof}

\newcommand{\mgc}{\operatorname{MGC}}

\begin{Def}
	We will call a graded circuit that admits Mikhalkin condition and can be constructed by gluing (see Construction \ref{con:graded_circuit_glue}) of elementary Mikhalkin circuits, a \emph{Mikhalkin graded circuit}.
	
	We will denote by $ \mgc(\Delta) $ the set of all Mikhalkin graded circuits in $ \Delta $.
\end{Def}

\begin{table}
	\begin{tabular}{| c | c | c |}
		\hline
		Type & $ \bs{w} $ & number of singular hyper-surfaces \\ \hline \hline
		$ \text{ I} \;\Vert \text{ I}  $ & $  [\bs{w}_{0} : \bs{w}_{1} : \bs{w}_{2}] $ & $ 2 d^{2} + O(d) $ \\ \hline
		
		$ \text{ II}  $ & $  [0 : \bs{w}_{1} : \bs{w}_{2}] $ & $ O(d) $ \\ \hline
		
		$ \text{ III}^{(2)}  $ & $  [0 : \bs{w}_{1} : \bs{w}_{2}] $ & $ \frac{1}{2} d^{2} + O(d) $ \\ \hline
		
		$ \text{ IV}_{2}  $ & $  [\bs{w}_{0} : 0 : \bs{w}_{2}] $ & $ \frac{1}{2} d^{2} + O(d) $ \\ \hline
		
	\end{tabular}
	
	\caption{A table of Mikhalkin graded circuits in dimension 2.}
	\label{table:GC_table_dim2}
\end{table}

\begin{table}
	
	\begin{tabular}{| c | c | c |}
		\hline
		Type & $ \bs{w} $ & number of singular hyper-surfaces \\ \hline \hline
		$ \text{ I} \;\Vert \text{ I} \;\Vert \text{ I}  $ & $  [\bs{w}_{0} : \bs{w}_{1} : \bs{w}_{2} : \bs{w}_{3}] $ & $ \frac{4}{3} d^{3} + O\left(d^{2}\right) $ \\ \hline
		
		$ \text{ II} \;\Vert \text{ I}  $ & $  [0 : \bs{w}_{1} : \bs{w}_{2} : \bs{w}_{3}] $ & $ O\left(d^{2}\right) $ \\ \hline
		
		$ \text{ III}^{(2)} \;\Vert \text{ I}  $ & $  [0 : \bs{w}_{1} : \bs{w}_{2} : \bs{w}_{3}] $ & $ \frac{1}{3} d^{3} + O\left(d^{2}\right) $ \\ \hline
		
		$ \text{ IV}_{2} \;\Vert \text{ I}  $ & $  [\bs{w}_{0} : 0 : \bs{w}_{2} : \bs{w}_{3}] $ & $ \frac{1}{3} d^{3} + O\left(d^{2}\right) $ \\ \hline
		
		$ \text{ I} \;\Vert \text{ II}  $ & $  [0 : \bs{w}_{1} : \bs{w}_{2} : \bs{w}_{3}] $ & $ O\left(d^{2}\right) $ \\ \hline
		
		$ \text{ I} \;\Vert \text{ III}^{(2)}  $ & $  [0 : \bs{w}_{1} : \bs{w}_{2} : \bs{w}_{3}] $ & $ \frac{2}{3} d^{3} + O\left(d^{2}\right) $ \\ \hline
		
		$ \text{ III}^{(3)}  $ & $  [0 : 0 : \bs{w}_{2} : \bs{w}_{3}] $ & $ \frac{1}{3} d^{3} + O\left(d^{2}\right) $ \\ \hline
		
		$ \text{ I} \;\Vert \text{ IV}_{2}  $ & $  [\bs{w}_{0} : \bs{w}_{1} : 0 : \bs{w}_{3}] $ & $ \frac{2}{3} d^{3} + O\left(d^{2}\right) $ \\ \hline
		
		$ \text{ IV}_{3} \;\Vert \text{ I}  $ & $  [\bs{w}_{0} : 0 : 0 : \bs{w}_{3}] $ & $ \frac{1}{3} d^{3} + O\left(d^{2}\right) $ \\ \hline
		
	\end{tabular}
	
	\caption{A table of Mikhalkin graded circuits in dimension 3.}
	\label{table:GC_table_dim3}
\end{table}

\begin{Example}
	We present a complete list of Mikhalkin graded circuits in dimensions $ n=2,3 $ and $ 4 $ in Tables \ref{table:GC_table_dim2}, \ref{table:GC_table_dim3} and \ref{table:GC_table_dim4} respectively.
	The first column of those tables indicate the types of circuits we need to glue in order to get the desired graded circuit. 
	For type III circuits we indicate their dimension by a superscript and for type IV circuits we indicate the difference between their upper and lower indices by a subscript.
	The second column shows for which $ \bs{w} $ we can get the specified graded circuit, by indicating which coordinates should equal to zero.
	Setting all the other coordinates to non-zero values would yield a graded circuit (at least if the degree $ d $ is high enough).
	For type III circuits with upper index $ r $, $ \bs{w}_r $ can be arbitrary and for those coordinates that occur as $ \bs{w}_r $ of type III circuit we can assign zero values as well, but for other types, the coordinates that do not have to be equal to zero should be strictly positive.
	
	Given a description of a graded circuit as in the first column of Tables \ref{table:GC_table_dim2} - \ref{table:GC_table_dim4}, the indices of the levels can be inferred as follows:
	We read the description from left to right; For levels of type I or II the index is the smallest non zero coordinate of $ \bs{w} $; For types III and IV, the lower index is the least non-zero coordinate and the upper index can be inferred either from the dimension (for type III) or the subscript (for type IV); We then can project along this circuit and continue as in Lemmas \ref{lemma:parrallel_projection} or \ref{lemma:gluing_type_IV}.
	
	The last column shows the (asymptotic) number of singular hypersurfaces corresponding to a particular graded circuit, see Sections \ref{sec:patchworking} and \ref{sec:discriminant_degree} for the proof.
	Note that type II circuits are always yielding an asymptotically small number of hypersurfaces.
	This follows from the fact that its dimension is 2 but it has only 1 parameter ($ \bs{w}_{\ell} $) given the value of the projection of $ \bs{w} $, see the proof of Proposition \ref{prop:discriminant_degree}.
\end{Example}

\begin{table}[h!]
	\centering
	
	\begin{tabular}{| c | c | c |}
		\hline
		Type & $ \bs{w} $ & number of singular hyper-surfaces \\ \hline \hline
		$ \text{ I} \;\Vert \text{ I} \;\Vert \text{ I} \;\Vert \text{ I}  $ & $  [\bs{w}_{0} : \bs{w}_{1} : \bs{w}_{2} : \bs{w}_{3} : \bs{w}_{4}] $ & $ \frac{2}{3} d^{4} + O\left(d^{3}\right) $ \\ \hline
		
		$ \text{ II} \;\Vert \text{ I} \;\Vert \text{ I}  $ & $  [0 : \bs{w}_{1} : \bs{w}_{2} : \bs{w}_{3} : \bs{w}_{4}] $ & $ O\left(d^{3}\right) $ \\ \hline
		
		$ \text{ III}^{(2)} \;\Vert \text{ I} \;\Vert \text{ I}  $ & $  [0 : \bs{w}_{1} : \bs{w}_{2} : \bs{w}_{3} : \bs{w}_{4}] $ & $ \frac{1}{6} d^{4} + O\left(d^{3}\right) $ \\ \hline
		
		$ \text{ IV}_{2} \;\Vert \text{ I} \;\Vert \text{ I}  $ & $  [\bs{w}_{0} : 0 : \bs{w}_{2} : \bs{w}_{3} : \bs{w}_{4}] $ & $ \frac{1}{6} d^{4} + O\left(d^{3}\right) $ \\ \hline
		
		$ \text{ I} \;\Vert \text{ II} \;\Vert \text{ I}  $ & $  [0 : \bs{w}_{1} : \bs{w}_{2} : \bs{w}_{3} : \bs{w}_{4}] $ & $ O\left(d^{3}\right) $ \\ \hline
		
		$ \text{ I} \;\Vert \text{ III}^{(2)} \;\Vert \text{ I}  $ & $  [0 : \bs{w}_{1} : \bs{w}_{2} : \bs{w}_{3} : \bs{w}_{4}] $ & $ \frac{1}{3} d^{4} + O\left(d^{3}\right) $ \\ \hline
		
		$ \text{ III}^{(3)} \;\Vert \text{ I}  $ & $  [0 : 0 : \bs{w}_{2} : \bs{w}_{3} : \bs{w}_{4}] $ & $ \frac{1}{6} d^{4} + O\left(d^{3}\right) $ \\ \hline
		
		$ \text{ I} \;\Vert \text{ IV}_{2} \;\Vert \text{ I}  $ & $  [\bs{w}_{0} : \bs{w}_{1} : 0 : \bs{w}_{3} : \bs{w}_{4}] $ & $ \frac{1}{3} d^{4} + O\left(d^{3}\right) $ \\ \hline
		
		$ \text{ IV}_{3} \;\Vert \text{ I} \;\Vert \text{ I}  $ & $  [\bs{w}_{0} : 0 : 0 : \bs{w}_{3} : \bs{w}_{4}] $ & $ \frac{1}{6} d^{4} + O\left(d^{3}\right) $ \\ \hline
		
		$ \text{ I} \;\Vert \text{ I} \;\Vert \text{ II}  $ & $  [0 : \bs{w}_{1} : \bs{w}_{2} : \bs{w}_{3} : \bs{w}_{4}] $ & $ O\left(d^{3}\right) $ \\ \hline
		
		$ \text{ II} \;\Vert \text{ II}  $ & $  [0 : 0 : \bs{w}_{2} : \bs{w}_{3} : \bs{w}_{4}] $ & $ O\left(d^{3}\right) $ \\ \hline
		
		$ \text{ III}^{(2)} \;\Vert \text{ II}  $ & $  [0 : 0 : \bs{w}_{2} : \bs{w}_{3} : \bs{w}_{4}] $ & $ O\left(d^{3}\right) $ \\ \hline
		
		$ \text{ IV}_{2} \;\Vert \text{ II}  $ & $  [0 : \bs{w}_{1} : 0 : \bs{w}_{3} : \bs{w}_{4}] $ & $ O\left(d^{3}\right) $ \\ \hline
		
		$ \text{ I} \;\Vert \text{ I} \;\Vert \text{ III}^{(2)}  $ & $  [0 : \bs{w}_{1} : \bs{w}_{2} : \bs{w}_{3} : \bs{w}_{4}] $ & $ \frac{1}{2} d^{4} + O\left(d^{3}\right) $ \\ \hline
		
		$ \text{ II} \;\Vert \text{ III}^{(2)}  $ & $  [0 : 0 : \bs{w}_{2} : \bs{w}_{3} : \bs{w}_{4}] $ & $ O\left(d^{3}\right) $ \\ \hline
		
		$ \text{ III}^{(2)} \;\Vert \text{ III}^{(2)}  $ & $  [0 : 0 : \bs{w}_{2} : \bs{w}_{3} : \bs{w}_{4}] $ & $ \frac{1}{8} d^{4} + O\left(d^{3}\right) $ \\ \hline
		
		$ \text{ IV}_{2} \;\Vert \text{ III}^{(2)}  $ & $  [0 : \bs{w}_{1} : 0 : \bs{w}_{3} : \bs{w}_{4}] $ & $ \frac{1}{8} d^{4} + O\left(d^{3}\right) $ \\ \hline
		
		$ \text{ I} \;\Vert \text{ III}^{(3)}  $ & $  [0 : 0 : \bs{w}_{2} : \bs{w}_{3} : \bs{w}_{4}] $ & $ \frac{1}{2} d^{4} + O\left(d^{3}\right) $ \\ \hline
		
		$ \text{ III}^{(4)}  $ & $  [0 : 0 : 0 : \bs{w}_{3} : \bs{w}_{4}] $ & $ \frac{1}{4} d^{4} + O\left(d^{3}\right) $ \\ \hline
		
		$ \text{ I} \;\Vert \text{ I} \;\Vert \text{ IV}_{2}  $ & $  [\bs{w}_{0} : \bs{w}_{1} : \bs{w}_{2} : 0 : \bs{w}_{4}] $ & $ \frac{1}{2} d^{4} + O\left(d^{3}\right) $ \\ \hline
		
		$ \text{ II} \;\Vert \text{ IV}_{2}  $ & $  [0 : \bs{w}_{1} : \bs{w}_{2} : 0 : \bs{w}_{4}] $ & $ O\left(d^{3}\right) $ \\ \hline
		
		$ \text{ III}^{(2)} \;\Vert \text{ IV}_{2}  $ & $  [0 : \bs{w}_{1} : \bs{w}_{2} : 0 : \bs{w}_{4}] $ & $ \frac{1}{8} d^{4} + O\left(d^{3}\right) $ \\ \hline
		
		$ \text{ IV}_{2} \;\Vert \text{ IV}_{2}  $ & $  [\bs{w}_{0} : 0 : \bs{w}_{2} : 0 : \bs{w}_{4}] $ & $ \frac{1}{8} d^{4} + O\left(d^{3}\right) $ \\ \hline
		
		$ \text{ I} \;\Vert \text{ IV}_{3} \;\Vert \text{ I}  $ & $  [\bs{w}_{0} : \bs{w}_{1} : 0 : 0 : \bs{w}_{4}] $ & $ \frac{1}{2} d^{4} + O\left(d^{3}\right) $ \\ \hline
		
		$ \text{ IV}_{4} \;\Vert \text{ I} \;\Vert \text{ I}  $ & $  [\bs{w}_{0} : 0 : 0 : 0 : \bs{w}_{4}] $ & $ \frac{1}{6} d^{4} + O\left(d^{3}\right) $ \\ \hline
		
		$ \text{ IV}_{4} \;\Vert \text{ IV}_{2}  $ & $  [\bs{w}_{0} : 0 : 0 : 0 : \bs{w}_{4}] $ & $ \frac{1}{12} d^{4} + O\left(d^{3}\right) $ \\ \hline
		
	\end{tabular}
\caption{A table of Mikhalkin graded circuits in dimension 4.}
\label{table:GC_table_dim4}
\end{table}

\newpage

\section{Patchworking}\label{sec:patchworking}

It follows from Lemmas \ref{lemma:origin_sing_condition} and \ref{lemma:mikhalkin_pos_imply_cond} that graded circuits admitting Mikhalkin condition correspond bijectively to singular hypersurfaces $ S\in \text{Sing}^{tr}(\Delta, \xpoints) $ together with a singular point $ \bs{\psi}\in \text{Sing}(S) $. 
Fix a graded circuit $ C\subseteq \Delta $ admitting Mikhalkin condition and let $ \bs{\psi}\in S $ be the corresponding tropical singular point and hypersurface.
In this section we will construct all algebraic singular hypersurfaces in $ (\mbb{K}^{*})^n $ of degree $ d $, passing through generic configuration of points whose tropicalization is $ \xpoints $, that tropicalize to $ S $ and have a singular point tropicalizing to $ \bs{\psi} $.
In particular, we will compute their number $ \text{mt}(S, \xpoints) $.

\begin{Def}\label{def:graded_circuit_multiplicity}
	If $ C $ is an elementary Mikhalkin circuit of type I as in Example \ref{ex:type_I_circ} we would set its multiplicity to $ m(C)=2 $.
	If $ C $ is an elementary Mikhalkin circuit of types II-IV (see Examples \ref{ex:type_II_circ}-\ref{ex:type_IV_circ}) we would set its multiplicity to $ m(C)=1 $.
	For a general Mikhalkin graded circuit $ C $ define its multiplicity $ m(C):=\prod_{1\le h\le ht(C)}m(C_h) $ to be the product of the multiplicities of its levels.
\end{Def}

We will show in this section that if $ \bs{\psi} $ is the singular point of $ S $ corresponding to the graded circuit $ C $, the number of algebraic lifts in $ (\mbb{K}^{*})^n $ of $ S $ with singularity tropicalizing to $ \bs{\psi} $, is equal to $ m(C) $.
We follow the patchworking technique as described, for example, in \cite[Chapter 2]{Mikhalkin_Shustin_Ittenberg}, \cite[Section 4.1]{Shustin2017_enum}, \cite{Shustin2005}.
To that end, write the defining polynomial $ f $ of an algebraic lift of $ S $ as
\begin{equation}\label{eq:patchworking_def_polynomial}
f(z) = \sum_{\bs{v}\in \Delta\cap\mbb{Z}^{n+1}} a_{\bs{v}} t^{\nu(\bs{v})} z^{\bs{v}}
\end{equation}
where $ \nu:\Delta\to\mbb{R} $ is the convex piece-wise linear function Legendre dual to the defining tropical polynomial of $ S $ and $ a_{\bs{v}}\in\mbb{K}^{*} $ with $ \val(a_{\bs{v}})=0 $ for all $ \bs{v} $.
We start by the first term of the coefficients $ a_{\bs{v}} $, namely $ a_{\bs{v}}(0) $ (see Remark \ref{remark:puiseux_is_function}), in Lemma \ref{lemma:patchwork_first_term}.
We next solve for higher terms of $ a_{\bs{v}} $, in Lemma \ref{lemma:patch_not_part_of_path} and Lemma \ref{lemma:patch_part_of_path}.
Note that some solutions $ a_{\bs{v}}(0) $ can be extended to a singular hypersurface that passes through the prescribed collection of points and tropicalize to $ S $ in more than one way.
To get a unique solution we will introduce additional variables, namely the coordinates of the singular point, which can be thought as blowing up the equations to get a transversal conditions.
In the case of $ \conv(C_1)\nsubseteq \Gamma_{\bs{w}} $ (Lemma \ref{lemma:patch_not_part_of_path}), this will be enough to get a system of equations with non-degenerate linearization, which then give us, by the implicit function theorem, solutions of the system in $ (\mbb{K}^{*})^N $.
In the harder case of $ \conv(C_1)\subseteq \Gamma_{\bs{w}} $, we would need to explore additional terms in the Puiseux series of the conditions to have singular point in order to get a non-degenerate system of equations, which can be thought of us additional blowing up.

Denote by $ r_h:=\dim \overline{C}_{h} $ for $ 1\le h\le ht(C) $.
We will work in affine coordinates with the additional assumption that  $ \Sp(\overline{C}_h)=\Sp\{ \bs{w}, \bs{e}^{(i)}-\bs{e}^{(j)} \; | \; i,j\le r_h \} $ for all $ 0< h\le ht(C) $. 
Note that this change of coordinates does not preserve the order $ \prec $, so we would not be able to invoke the fact that $ C $ admits Mikhalkin condition in those coordinates.
By translating tropically by $ -\bs{\psi} $, which algebraically amounts to multiplying each coordinate by $ t^{-\bs{\psi}_i} $, we can put the tropical singularity at $ \boldsymbol{0} $.
Fix points $ \bs{p}^{(k)}\in (\mbb{K}^{*})^n $ for $ k=1,\dots,N-1 $ with $ \text{Val}(\bs{p}^{(k)})=-M_k\cdot \bphi+\bs{\psi} $ (cf. Construction \ref{con:Mikhalkin_position}) and denote $ \mu_{k,i}:=-M_k\cdot \bphi_i + \bs{\psi}_i $, so we can write $ \bs{p}^{(k)}_i = \xi_{k,i}t^{\mu_{k,i}} $ where $ \xi_{k,i}(0)\ne 0 $. 
Denote the translated hypersurface by $ S_0:=\{ \bs{\chi}-\bs{\psi} \; | \; \bs{\chi}\in S \} $.

We will use the following notation throughout this section:
\begin{Def}
	The notation $ o(t^s) $ (here $ s\in \mbb{R} $) would mean "some Puiseux series $ a\in\mbb{K} $ with $ \val(a)>s $".
	So, for example, $ a_{\bs{v}} = \alpha_{\bs{v}}+o(1) $ would mean that $ a_{\bs{v}} $ has constant term $ \alpha_{\bs{v}} $.
	This coincide with the usual little $ o $ notation (as $ t\to 0 $) if elements of $ \mbb{K} $ are viewed as germs (in a punctured neighborhood of 0) of functions $ \mbb{C}\to \mbb{C} $ after parameter change $ t\mapsto t^m $, see Remark \ref{remark:puiseux_is_function}.
\end{Def}

Denote by $ Ini(f)(z)=\sum_{\bs{v}\in C_0\cup C_1}a_{\bs{v}}(0)z^{\bs{v}}\in \mbb{C}[z] $.
Recall that by \cite[Lemma 4.3]{Shustin2017_enum}, if $ f $ has a singularity in $ (\mbb{K}^*)^n $ and $ \bs{\psi} $ is of maximal-dimensional geometric type, $ V(Ini(f)) $ has a singularity in $ (\mbb{C}^*)^n $.

\begin{lemma}\label{lemma:patchwork_first_term}
	If $ f(z) $ of the form \eqref{eq:patchworking_def_polynomial} defines a hypersurface in $ (\mbb{K}^*)^n $ passing through the configuration $ (\bs{p}^{(k)})_{k\le N-2} $, and the polynomial $ Ini(f) $ has a singularity in $ (\mbb{C}^*)^n$, then:
	\begin{enumerate}
		\item 
		If $ C_1 $ is an elementary Mikhalkin circuit of type II, III or IV, there exists a unique, up to multiplication by non-zero scalar, possible value for $ (a_{\bs{v}}(0))_{\bs{v}\in \Delta\cap\mbb{Z}^n} $.
		
		\item
		If $ C_1 $ is an elementary Mikhalkin circuit of type I and $ \conv(C_1)\nsubseteq \Gamma_{\bs{w}} $ there exists, up to multiplication by non-zero scalar, 2 possible values for $ (a_{\bs{v}}(0))_{\bs{v}\in \Delta\cap\mbb{Z}^n} $.
		
		\item
		If $ C_1 $ is an elementary Mikhalkin circuit of type I and $ \conv(C_1)\subseteq \Gamma_{\bs{w}} $ there exists a unique, up to multiplication by non-zero scalar, possible value for $ (a_{\bs{v}}(0))_{\bs{v}\in \Delta\cap\mbb{Z}^n} $.
		
	\end{enumerate}
\end{lemma}

\begin{proof}
	We begin by investigating the condition imposed by the marked points $ \bs{p}^{(k)} $.
	Write as before $ \bs{p}^{(k)}_i:=\xi_{k,i}t^{\mu_{k,i}} $, then the tropical point $ \bchi{k}-\bs{\psi} $ with $ i $'th coordinate $ \mu_{k,i} $ lies on a $ (n-1) $-face dual to an edge $ E_k\subseteq \Gamma_{\bs{w}} $.
	Thus,
	\begin{equation}\label{eq:lemma_mult_point_cond}
	0=t^{-s'_k}f(\bs{p}^{(k)})=t^{-s'_k}\sum_{\bs{v}\in \Delta\cap\mbb{Z}^n}a_{\bs{v}} t^{\nu(\bs{v})} \left(\bs{p}^{(k)}\right)^{\bs{v}} = \sum_{\bs{v}\in E_k\cap\mbb{Z}^n}a_{\bs{v}}\xi_{k}^{\bs{v}} + o(1)
	\end{equation}
	for some $ s'_k\in\mbb{R} $.
	
	We first deal with the case that $ \conv(C_1)\nsubseteq \Gamma_{\bs{w}} $ (which is always the case if $ C_1 $ is of type II, III or IV).
	Then all the edges of $\Gamma_{\bs{w}}$ have lattice length 1 meaning that the equations $ t^{-s'_k}f(\bs{p}^{(k)})=0 $ consist of only 2 terms at $ t=0 $,
	\[ a_{\bs{v}_1}\xi_k^{\bs{v}_1}+a_{\bs{v}_2}\xi_k^{\bs{v}_2}=o(1) \]
	which gives as a unique, up to multiplication by a non-zero scalar, solution to $ (a_{\bs{v}}(0))_{{\bs{v}}\in\Delta\cap\mbb{Z}^n\setminus\{w\}} $.
	
	Since we assumed $ Ini(f) $ is singular, $ (a_{\bs{v}}(0))_{\bs{v}\in C_1\cup C_0} $ satisfy the discriminantal equation 
	\begin{equation}\label{eq:lemma_mult_discriminant}
	\prod_{\bs{v}\in C_1\cup C_0} a_{\bs{v}}(0)^{m_{\bs{v}}} = \prod_{\bs{v}\in C_1\cup C_0}m_{\bs{v}}^{m_{\bs{v}}} 
	\end{equation}
	where $ \sum m_{\bs{v}} \bs{v}=0 $ is the circuit relation of $ C_1\cup C_0 $ (see \cite[Chapter 9, Proposition 1.8]{Gelfand_kapranov}).
	If $ C_1 $ is of type II, III or IV, the coefficient $ m_{\bs{w}} $ in the circuit relation is 1 by Remark \ref{remark:C_1_is_regular}, meaning that the system $ t^{-s'_k}f(\bs{p}^{(k)})=0 $ together with \eqref{eq:lemma_mult_discriminant} have a unique solution at $ t=0 $.
	If, on the other hand, 	$ C_1 $ is of type I, \eqref{eq:lemma_mult_discriminant} becomes the familiar quadratic discriminant $ a_{\bs{w}}^2=4a_{\bs{u}^{(-1)}}a_{\bs{u}^{(1)}} $ (here $ C_1=\{ u^{(-1)}, u^{1} \} $ as in Example \ref{ex:type_I_circ}) which yields 2 solutions at $ t=0 $ for the above system, $ a_{\bs{w}}(0)=\pm \sqrt{a_{\bs{u}^{(-1)}}(0)\cdot a_{\bs{u}^{(1)}}(0)} $.
	
	Now, if $ \conv(C_1)\subseteq \Gamma_{\bs{w}} $, the fact that $ Ini(f) $ is singular quadratic polynomial (in one variable) means that it is a squared binomial. 
	The point condition that corresponds to $ \conv(C_1) $ fixes the coefficients of this binomial (up to proportionality) which fixes the lowest term of the coefficients $ \left(a_{\bs{v}}(0)\right)_{\bs{v}\in C_0\cup C_1} $ up to proportionality.
	The rest of the lower terms of coefficients of $ f $ determined by the point conditions corresponding to edges of $ \Gamma_{\bs{w}} $ with lattice length 1.
\end{proof}

\begin{lemma}\label{lemma:patch_not_part_of_path}
	If $ \conv(C_1) \nsubseteq \Gamma_{\bs{w}} $, then there exists $ m(C) $ algebraic surfaces passing through $ \overline{\bs{p}} $ that tropicalize to $ S_0 $ and having a singularity that tropicalize to $ \bs{0} $. 
\end{lemma}

\begin{proof}
	Pick $ \bs{v}^{(0)}\in \Delta\cap\mbb{Z}^n\setminus \{\bs{w}\} $.
	Since $ a_{\bs{v}^{(0)}}\ne 0 $, we can scale all the coefficients and assume, without loss of generality, that $ a_{\bs{v}^{(0)}}=1 $.
	For every $ 0< i \le n $ let $ s_i\in \mbb{R} $ be the maximal s.t. $ t^{-s_i}\der{z_i}f $ contains no coefficients with negative valuation and $ s_0\in\mbb{R} $ be the maximal s.t. $ t^{-s_0}f $ contains no coefficients with negative valuation.
	Consider the map 
	\begin{align*}
	\Psi : \mbb{C}\times\mbb{C}^N\times \mbb{C}^n & \to \mbb{C}^{N+n} \\
	\left( t, \left( a_{\bs{v}} \right)_{\bs{v}\in \Delta\cap\mbb{Z}^n\setminus \{\bs{v}^{(0)}\}}, \left( \bs{q}_i \right)_{1\le i\le n}   \right) & \mapsto \left( \left(t^{-s'_k} f(\bs{p}^{(k)}) \right)_{1\le k\le N-1}, t^{-s_0}f(\bs{q}), \left( t^{-s_i}\bs{q}_i\der{z_i}f|_{z=\bs{q}} \right)_{1\le i \le n}   \right) 
	\end{align*}
	where $ N:=|\Delta\cap\mbb{Z}^n|-1 $ and $ s'_k $ as in \eqref{eq:lemma_mult_point_cond}.
	We will show that the fiber $ \Psi^{-1}(\boldsymbol{0}) $ consisting of $ m(C) $ points with $ t=0 $ which are all regular points of $ \Psi $.
	Thus, by the implicit function theorem, there exists Puiseux series $ a_{\bs{v}},\bs{q}_i\in \mbb{K} $ s.t. $ \Psi(t, \left( a_{\bs{v}}(t) \right), \left( \bs{q}_i(t) \right)  )=\bs{0} $ for small enough (in absolute value of $ \mbb{C} $) values of $ t $.
	Since the conditions of a hypersurface to pass through $ (\bs{p}^{(k)})_{k\le N-1} $ encoded in the first coordinates of $ \Psi $, and the condition of $ \bs{q} $ to be a singular point of $ V(f) $, where $ f(z):=\sum_{\bs{v}\in\Delta\cap\mbb{Z}^n}a_{\bs{v}}t^{\nu(\bs{v})}z^{\bs{v}} $ is encoded in the last coordinates of $ \Psi $, we will get the desired result.
	
	By Lemma \ref{lemma:patchwork_first_term}, there exist $ m(C_1) $ possible solutions for $ \left( a_{\bs{v}} \right)_{\bs{v}\in \Delta\cap\mbb{Z}^{n+1}}  $.
	First, we will find one solution of $ \left( \bs{q}_i(0) \right)_{i\le r_1}  $ for each of them, where we denote by $ r_1:=\dim C_1 $, the dimension of the first level of $ C $.
	If $ C_1 $ is of type II, III or IV, then it is a regular simplex by Remark  \ref{remark:C_1_is_regular}, so we can change coordinates s.t. $ C_1=\{ 0,e_1,\dots,e_{r_1} \} $.
	Then 
\begin{equation}\label{eq:lemma_patch_not_in_path_derivatives}
	0 =\left (z_i \der{z_i} f\right )|_{z=\bs{q}} = a_{\bs{e}^{(i)}}\bs{q}_i+\bs{w}_i a_{\bs{w}} \bs{q}^{\bs{w}} + o(1)
\end{equation}
	for all $ i\le r_1 $, and so, 
\begin{equation}\label{eq:lemma_patch_not_in_path_func}
	0 = f(\bs{q}) = a_{\bs{0}} + \sum_{i=1}^{r_1} a_{\bs{e}^{(i)}}\bs{q}_i + a_{\bs{w}} \bs{q}^{\bs{w}} + o(1) = a_{\bs{0}} + \left (1-\sum_{i=1}^{r_1}  \bs{w}_i\right )a_{\bs{w}}\bs{q}^{\bs{w}} + o(1)
\end{equation}
	meaning that 
	\[ \bs{q}_i(0) = -\bs{w}_i \frac{a_{\bs{w}}(0)}{a_{\bs{e}^{(i)}}(0)}\bs{q}(0)^{\bs{w}} = \frac{a_{\bs{0}}(0)}{a_{\bs{e}^{(i)}}(0)}\left (1-\sum_{i=1}^{r_1} \bs{w}_i\right )^{-1}. \]
	Note that $ \sum_{i=1}^{r_1} \bs{w}_i\ne 1 $ since $ w $ does not lie on a hyperplane with $ \bs{e}^{(1)},\dots,\bs{e}^{(r_1)} $.
	
	If $ C_1 $ is of type I then it is equivalent to $ \{0, 2\bs{e}^{(1)} \} $ with $ C_0=\{ \bs{e}^{(1)} \} $.
	Since \linebreak $ {a_{\bs{e}^{(1)}}^2-4a_{\bs{0}}a_{2\bs{e}^{(1)}}=0} $ by \eqref{eq:lemma_mult_discriminant}, $ f(z)=a_{2\bs{e}^{(1)}}(z_1+\frac{a_{\bs{e}^{(1)}}}{2 a_{2\bs{e}^{(1)}}})^2+o(1) $ meaning that the unique solution for the first term is $ \bs{q}_1(0)=-\frac{a_{\bs{e}^{(1)}}}{2 a_{2\bs{e}^{(1)}}} $.
	
	We thus showed that for every solution of equations \eqref{eq:lemma_mult_point_cond},\eqref{eq:lemma_mult_discriminant} at $ t=0 $ there exists a unique solution of the equations $ t^{-s_0}f(\bs{q})=t^{-s_1}\bs{q}_1\der{z_1}f|_{z=\bs{q}}=\dots=t^{-s_{r_1}}\bs{q}_{r_1}\der{z_{r_1}}f|_{z=\bs{q}}=0 $ at $ t=0 $ (in the variables $ \bs{q}_1,\dots,\bs{q}_{r_1} $).
	
	By Lemma \ref{lemma:patch_high_levels} below, for every such a solution, there exist $ \prod_{h>1}m(C_h)=\frac{m(C)}{m(C_1)} $ solutions $ (\bs{q}_i(0))_{i>r_1} $ for $ (t^{-s_i}\der{z_i}f|_{t=0,z=\bs{q}}=0 )_{i>r_1} $, together with the $ m(C_1) $ solutions for $ a_{\bs{w}}(0) $ we get the desired $ m(C) $ points in the fiber $ \Psi^{-1}(\boldsymbol{0}) $, so we are left to show that for all of them $ \Psi $ linearizes to a non-degenerate system of linear equations.
	We arrange the variables and equations into blocks as follows: 
	\begin{itemize}
		\item The first block contains the variables $ (a_{\bs{v}})_{\bs{v}\in \Delta\cap\mbb{Z}^n\setminus \{ \bs{w}, \bs{v}^{(0)} \} } $ and the equations \linebreak $ {(t^{-s_j'}f(\bs{p}^{(j)})=0)_{j\le N-1}} $.
		
		\item The second block contains the variables $ a_{\bs{w}}, (\bs{q}_i)_{i\le r_1} $ and the equations 
		\[t^{-s_0}f(\bs{q})=t^{-s_1}\bs{q}_1\der{z_1}f=\dots=t^{-s_{r_1}}\bs{q}_{r_1}\der{z_{r_1}}f=0 .\]
		
		\item A block for every $ 1<h\le ht(C) $ consisting the variables $ \bs{q}_i $ and the equations $ t^{-s_i}\bs{q}_i\der{z_i}f=0 $ for all $ r_{h-1} < i \le r_h $.
	\end{itemize}
	The equations \eqref{eq:lemma_mult_point_cond} do not depend on $ \bs{q} $ and $ a_{\bs{w}} $ (in case $ \conv(C_1)\nsubseteq \Gamma $) meaning that only the first block in the first row is non-zero and it is (suitably arranged) triangular matrix with non-zero entries on the diagonal, thus invertible.
	By Lemma \ref{lemma:patch_high_levels} again, all other blocks above the diagonal are zero as well and the blocks on the diagonal that correspond to levels $ C_h $ with $ h>1 $ are invertible.
	We thus are left to show that the block on the diagonal corresponding to $ C_1 $ is invertible.
	
	If $ C_1 $ is of type I, there are 2 equations and 2 variables and the corresponding block of the Jacobian is
	\[ \begin{pmatrix}
	\bs{q}_1(0) & 2a_{\bs{u}^{(1)}}(0)\bs{q}_1(0)+a_{\bs{w}}(0) \\
	1 & 2a_{\bs{u}^{(1)}}
	\end{pmatrix} \]
	where $ \bs{u}^{(1)} $ is one of the point of $ C_1 $, see Example \ref{ex:type_I_circ}.
	So the determinant is $ -a_{\bs{w}}(0)\ne 0 $ and the block is invertible.
	If $ C_1 $ is of type II,III or IV, we set as before $ C_1=\{ 0, \bs{e}^{(1)},\dots,\bs{e}^{(r_1)} \} $ so the equations are \eqref{eq:lemma_patch_not_in_path_derivatives}, \eqref{eq:lemma_patch_not_in_path_func} and the corresponding block in the Jacobian is (here the first column corresponds to the variable $ a_{\bs{w}} $ and the first row to the equation $ t^{-s_0}f(\bs{q})=0 $) 
	\[  \begin{pmatrix}
	\bs{q}^{\bs{w}} & a_{\bs{e}^{(1)}}+\bs{w}_1a_{\bs{w}} \bs{q}^{\bs{w}-\bs{e}^{(1)}} &  & a_{\bs{e}^{(r_1)}}+\bs{w}_{r_1} a_{\bs{w}} \bs{q}^{\bs{w}-\bs{e}^{(r_1)}} \\
	\bs{w}_1 \bs{q}^{\bs{w}} & a_{\bs{e}^{(1)}}+\bs{w}_1^2a_{\bs{w}} \bs{q}^{\bs{w}-\bs{e}^{(1)}} & \cdots & \bs{w}_1\bs{w}_{r_1}a_{\bs{w}}\bs{q}^{\bs{w}-\bs{e}^{(r_1)}} \\
	\bs{w}_2\bs{q}^{\bs{w}} & \bs{w}_1\bs{w}_2a_{\bs{w}}\bs{q}^{\bs{w}-\bs{e}^{(1)}} &  & \bs{w}_2\bs{w}_{r_1}a_{\bs{w}}\bs{q}^{\bs{w}-\bs{e}^{(r_1)}} \\
	& \vdots & \ddots & \vdots \\
	\bs{w}_{r_1}\bs{q}^{\bs{w}} & \bs{w}_1\bs{w}_{r_1}a_{\bs{w}}\bs{q}^{\bs{w}-\bs{e}^{(1)}} & \cdots & a_{\bs{e}^{(r_1)}} + \bs{w}_{r_1}^2 a_{\bs{w}}\bs{q}^{\bs{w}-e^{(r_1)}}
	\end{pmatrix} \]
	Subtracting all the rows from the first one and scaling the first row and column accordingly, we get 
	\[ 
	\begin{pmatrix}
	1 & \bs{w}_1a_{\bs{w}}\bs{q}^{\bs{w}-\bs{e}^{(i)}} & & \bs{w}_{r_1}a_{\bs{w}}\bs{q}^{\bs{w}-\bs{e}^{(r_1)}} \\
	\bs{w}_1 & a_{\bs{e}^{(1)}}+\bs{w}_1^2a_{\bs{w}} \bs{q}^{\bs{w}-\bs{e}^{(1)}} & \cdots & \bs{w}_1\bs{w}_{r_1}a_{\bs{w}}\bs{q}^{\bs{w}-\bs{e}^{(r_1)}} \\
	\bs{w}_2 & \bs{w}_1\bs{w}_2a_{\bs{w}}\bs{q}^{\bs{w}-\bs{e}^{(1)}} &  & \bs{w}_2\bs{w}_{r_1}a_{\bs{w}}\bs{q}^{\bs{w}-\bs{e}^{(r_1)}} \\
	& \vdots & \ddots & \vdots \\
	\bs{w}_{r_1} & \bs{w}_1\bs{w}_{r_1}a_{\bs{w}}\bs{q}^{\bs{w}-\bs{e}^{(1)}} & \cdots & a_{\bs{e}^{(r_1)}} + \bs{w}_{r_1}^2 a_{\bs{w}}\bs{q}^{\bs{w}-\bs{e}^{(r_1)}}
	\end{pmatrix}
	\]
	Subtracting $ \bs{w}_i $ multiple of the first row from the $ (i+1) $'th (for $ i=1,\dots,r_1 $) we get an upper triangular matrix with non-zero entries on the diagonal, meaning that the Jacobian is invertible.
\end{proof}

For the treatment of higher levels of $ C $ we will introduce some notation.
For $ 1 \le h \le ht(C) $, denote by $ r_h:\dim \overline{C}_h $ and choose coordinates with $ \bs{w}=\bs{0} $ and $ {\Sp(\overline{C}_{r_{h}})=\Sp \{e^{(i)} \; | \; i\le r_h \}} $ for all $ 1\le h \le ht(C) $.
For $ 1\le i\le n $ let $ s_i\in \mbb{R} $ be the maximal s.t. $ t^{-s_i}\der{z_i}f $ contains no coefficients with negative valuation, and let $ s_0\in \mbb{R} $ be the maximal s.t. $ t^{-s_0}f $ contains no coefficients with negative valuation.
For $ 1 \le h\le ht(C) $, consider the map 
\begin{align*}
\Phi_h : \mbb{C}^{n+1} & \to \mbb{C}^{r_h+1} \\ 
\left (t, \bs{q}_1,\dots, \bs{q}_n \right ) & \mapsto (t^{-s_0}f(\bs{q}), t^{-s_1}\der{z_1}f|_{z=\bs{q}},\dots,t^{-s_{r_h}}\der{z_{r_h}}f|_{z=\bs{q}}).
\end{align*}

\begin{lemma}\label{lemma:patch_high_levels}
	Let $ (a_{\bs{v}}(0))_{\bs{v}\in \Delta\cap\mbb{Z}^n} $ be a solution to the point conditions and the discriminantal equation of $ C_1 $, as found in Lemma \ref{lemma:patchwork_first_term} and let $ 1<h\le ht(C) $.
	Then the value of $ \Phi_{h-1}(0, \bs{q}) $ does not depend on $ \bs{q}_{r_{h-1}+1},\dots, \bs{q}_n $, and for every solution $ \left( \bs{q}_i\in \mbb{C}^{*}\right)_{i\le r_{h-1}} $ of 
	\begin{equation}\label{eq:lemma_patch_high_levels_small_system} 
	\Phi_{h-1}(0,\bs{q})=0,
	\end{equation}
	there exist $ m(C') $ values for $ \left( \bs{q}_i \in \mbb{C}^{*} \right)_{r_{h-1}< i \le r_h}  $ s.t. 
	\begin{equation}\label{eq:lemma_patch_high_levels_big_system} 
	\Phi_{h}(0,\bs{q})=0.
	\end{equation}
	Moreover, the Jacobian of 
	\begin{equation*}
	\begin{split}
	\Phi_h' : (\bs{q}_i)_{r_{h-1}<i\le r_h} & \mapsto (t^{-s_i}z_i\der{z_i}f|_{z=\bs{q}})_{r_{h-1}<i\le r_h}
	\end{split}
	\end{equation*} 
	is invertible at each of those solutions and $ t=0 $.
	
\end{lemma}

\begin{proof}
	Let $ \bs{q} $ be a solution of $ \Phi_h(0, \bs{q})=0 $.
	For $ r_{h-1}< i\le r_h $ we have,
	\[ 0 = (z_i\der{z_i} f)|_{z=\bs{q}} = t^s\sum_{\bs{v}\in C_h} \bs{v}_i a_{\bs{v}} \bs{q}^{\bs{v}}+ o(t^s) \]
	where $ s:=\nu(C_h) $.
	Combining those equations together and comparing the coefficient of $ t^s $, we get that $ \sum_{\bs{v}\in C_h}(a_{\bs{v}}\bs{q}^{\bs{v}}+o(1))\bs{v}\in \Sp(\overline{C}_{h-1}) $.
	Since the linear combination with $ {\sum_{\bs{v}\in C_h}m_{\bs{v}} \bs{v}\in \Sp(\overline{C}_{h-1})} $ is unique up to re-scaling, $ a_{\bs{v}}\bs{q}^{\bs{v}}+o(1) = \lambda m_{\bs{v}} $ for all $ \bs{v}\in C_h $ (with the same $ \lambda $) where $ m_{\bs{v}} $ as in Corollary \ref{cor:graded_circuit_linear_combination_with_sum_1}.
	Then, for every $ \bs{u}\in C_h $,
	\[ \prod_{\bs{v}\in C_h}\left( \frac{a_{\bs{v}}\bs{q}^{\bs{v}}+o(1)}{m_{\bs{v}}} \right)^{m_{\bs{v}}}m_{\bs{u}} = \lambda^{\sum m_{\bs{v}}}m_{\bs{u}}=\lambda m_{\bs{u}}=a_{\bs{u}}\bs{q}^{\bs{u}}+o(1)  \]
	so
	\begin{equation}\label{eq:lemma_mult_high_level_no_double}
	\bs{q}^{\bs{u}} = \frac{m_{\bs{u}}}{a_{\bs{u}}}\prod\left( \frac{a_{\bs{v}}}{m_{\bs{v}}} \right)^{m_{\bs{v}}}\cdot \bs{q}^{\sum m_{\bs{v}}\bs{v}} + o(1)
	\end{equation}
	is a function of $ (a_{\bs{v}})_{\bs{v}\in\Delta\cap \mbb{Z}^n} $ and $ (\bs{q}_i(0))_{i\le r_{h-1}} $.
	
	If $ C_h $ is of types II,III or IV, it is a regular simplex, and we can extract $ \bs{q}_i $ from equations \eqref{eq:lemma_mult_high_level_no_double}, so they have a unique solution.
	Moreover, since the equations that we get from \eqref{eq:lemma_mult_high_level_no_double} after extracting $ \bs{q}_i $ are equivalent to the original equations, the Jacobian we are after is invertible if and only if the Jacobian with respect to $ (\bs{q}_i)_{r_{h-1}< i \le r_h} $ of those equations is invertible, but the later is the identity.
	
	In case of $ C_h $ of type I, we have 2 solutions for $ \bs{q}_{r_h}(0) $.
	Indeed, $ C_h= \{ \bs{u}^{(-1)}, \bs{u}^{(1)} \} $ with $ \bs{u}^{(-1)}+\bs{u}^{(1)}\in \Sp(\overline{C}_{h-1}) $ and we can pick coordinates s.t. $ \bs{u}^{1}=e_{r_h} $ and $ \bs{w}_{r_h}=0 $.
	Then 
	\[ 0 = t^{-s}\der{z_{r_h}}f(\bs{q}) = a_{\bs{u}^{(1)}} - a_{\bs{u}^{(-1)}}\prod_{i=1}^{r_h-1}\bs{q}_i^{\bs{u}^{(-1)}_i}\bs{q}_{r_h}^{-2}+o(1)  \]
	meaning that 
	\begin{equation}\label{eq:lemma_patch_high_levels_typeI_sol}
	\bs{q}_{r_h}(0) = \pm \sqrt{ \frac{a_{\bs{u}^{(-1)}}(0)}{a_{\bs{u}^{(1)}}(0)}\prod_{i=1}^{r_h-1}\bs{q}_i(0)^{\bs{u}^{(-1)}_i} }.
	\end{equation}
	The Jacobian in this case is a scalar (order 1 matrix) and it is non-zero since it is equal to the derivative of a separable quadratic at its root.
\end{proof}

We will now consider the case when $ \conv(C_1)\subseteq \Gamma_{\bs{w}} $.
When this happens, $ C_1 $ is of type I and index $ 0 $ (see Example \ref{ex:type_I_circ}), denote $ C_1=\{ \bs{u}^{(-1)}, \bs{u}^{(1)} \} $ with $ \bs{u}^{(-1)}\prec \bs{w}\prec \bs{u}^{(1)} $.
It would be convenient to choose coordinates with $ \bs{u}^{(-1)} $ at the origin and, as before, $ {\Sp(\overline{C}_{h})=\Sp\{\bs{e}^{(i)} \mid i\le r_h\}} $.
In the proof we would need the following property of Construction \ref{con:Mikhalkin_position}.
Recall that for $ j=1,\dots, N-1 $ and $ i=1,\dots,n $, we write $ \bs{p}^{(j)}_i = \xi_{j,i}t^{\mu_{j,i}} $ where $ \xi_{j,i}\in\mbb{K} $ with $ \val(\xi_{j,i})=0 $ and $ \mu_{j,i} = -M_j\cdot \bphi_i + \bs{\psi}_i $.

\begin{lemma}\label{lemma:patch_part_of_path_f_0_trop}
	Suppose that $ \conv(C_1)\subseteq \Gamma_{\bs{w}} $.
	Let $  $
	\[ f_0(z) = \sum_{\bs{v}\in \Delta\cap\mbb{Z}^{n+1} \setminus \Sp(C_1)} a_{\bs{v}}t^{\nu(\bs{v})}z^{\bs{v}}, \]
	$ \bs{q} $ be a singular point of 
	\[ f(z) = \sum_{\bs{v}\in \Delta\cap\mbb{Z}^{n+1}}a_{\bs{v}}t^{\nu(\bs{v})}z^{\bs{v}} \]
	and let $ k $ be the unique index s.t. $ \trop(\bs{p}^{(k)}) $ contained in the $ (n-1) $-cell that correspond to the edge $ Conv(C_1) $ in $ \Gamma_{\bs{w}} $.
	Then ${ \val(f_0(\bs{q})) > \val(f_0(\bs{p}^{(k)}))} $ and the only term of $ f_0(\bs{p}^{(k)}) $ with minimal valuation is $ a_{\bs{w}^{\pr}}t^{\nu(\bs{w}^{\pr})}\xi_k^{\bs{w}^{\pr}} $ where $ \bs{w}^{\pr} $ is the last (w.r.t. $ \prec $) point of $$ {\{ \bs{v}\in \Delta\cap\mbb{Z}^n \; | \; \bs{v}\prec \bs{w}, \bs{v}\notin \Sp(C_1) \} } $$.
\end{lemma}

\begin{proof}
	Denote by $ \nu_0:\Delta\to\mbb{R} $ the convex piece-wise linear function Legendre dual to the defining polynomial of the tropical hypersurface before the translation of $ \bs{\psi} $ to $ \boldsymbol{0} $, i.e. $ {\nu_0 = \nu + \bs{\psi}} $.
	
	We begin by computing the valuation $ \val(f_0(\bs{p}^{(k)})) $
	\[ f_0(\bs{p}^{(k)})=\sum_{\bs{v}\in \Delta\cap\mbb{Z}^n\setminus \Sp(C_1)}a_{\bs{v}}t^{\nu_0(\bs{v})- M_k\pair{\bphi}{\bs{v}}}\xi_{k}^{\bs{v}} \]
	where $ M_k, \bphi $ are as in Construction \ref{con:Mikhalkin_position}.
	The value of $ \nu_0(\bs{v}) $ depends with positive coefficient on all $ M_j $ that correspond to edges of $ \Gamma_{\bs{w}} $ preceding $ \bs{v} $, meaning that 
	\[\nu_0(\bs{v}) \gg M_k\pair{\bphi}{\bs{v}}+\nu_0(\bs{w}^{\pr})-M_k\pair{\bphi}{\bs{w}^{\pr}} \] for all $ \bs{v}\succ \bs{w} $.
	Conversely, for every $ \bs{v}\prec \bs{w} $, $ \nu_0(\bs{v})\ll M_k\pair{\bphi}{\bs{v}} $, meaning that the minimum of $ \nu_0-M_k\bphi $ attained whenever $ \bphi $ attains its maximum, which is only at $ \bs{w}^{\pr} $.
	This means that the only term of $ f_0(\bs{p}^{(k)}) $ with the minimal valuation is $ a_{\bs{w}^{\pr}}t^{\nu(\bs{w}^{\pr})}\xi_k^{\bs{w}^{\pr}} $ where $ \bs{w}^{\pr} $ is the last (w.r.t. $ \prec $) point with $ \bs{w}^{\pr}\prec \bs{w} $ and $ \bs{w}^{\pr}\notin \Sp(C_1) $.
	
	To compute $ \val(f_0(\bs{q})) $, recall that $ \bs{\psi} $ is defined up to an additive constant, so we can choose $ \pair{\bs{\psi}}{\bs{u}^{(-1)}}=0 $.
	Then, since $ \nu(\bs{u}^{(-1)})=\nu(\bs{u}^{(1)}) $, $$ \pair{\bs{\psi}}{\bs{u}^{(1)}} = \nu(\bs{u}^{(1)})-\nu_0(\bs{u}^{(1)}) =  \nu_0(\bs{u}^{(-1)})-\nu_0(\bs{u}^{(1)}). $$
	Pick coefficients $ m_{\bs{v}} $ for $ \bs{v}\in C_1\cup C_2 $ such that $ \sum_{\bs{v}\in C_1}m_{\bs{v}}\bs{v}=\sum_{\bs{v}\in C_2}m_{\bs{v}}\bs{v} $ and \linebreak $ {\sum_{\bs{v}\in C_1}m_{\bs{v}}=\sum_{\bs{v}\in C_2}m_{\bs{v}}=1} $.
	Then
	\begin{align*}
	\val(f_0(\bs{q})) =\nu(C_2) &=\sum_{\bs{v}\in C_2}m_{\bs{v}}\nu(\bs{v})= \\
	=\sum_{\bs{v}\in C_2}m_{\bs{v}}(\nu_0(\bs{v})-\pair{\bs{\psi}}{\bs{v}}) &= \sum_{\bs{v}\in C_2}m_{\bs{v}}\nu_0(\bs{v})-\sum_{\bs{v}\in C_2}m_{\bs{v}}\pair{\bs{\psi}}{\bs{v}} = \\
	= \sum_{\bs{v}\in C_2}m_{\bs{v}}\nu_0(\bs{v})-\sum_{\bs{v}\in C_1}m_{\bs{v}} \pair{\bs{\psi}}{\bs{v}} &= \sum_{\bs{v}\in C_2}m_{\bs{v}}\nu_0(\bs{v})+m_{\bs{u}^{(1)}}\left (\nu_0(\bs{u}^{(1)})-\nu_0(\bs{u}^{(-1)})\right ).
	\end{align*}
	We will show that the last (w.r.t. $ \prec $) point appearing in the right-hand side of the above equation succeeds $ \bs{w} $ and that its coefficient, (i.e. $ m_{\bs{v}} $) is positive, which will imply that $ \val(f_0(\bs{q}))>\nu_0(\bs{w}^{\pr})-M_k\pair{\bphi}{\bs{w}^{\pr}} =\val(f_0(\bs{p}^{(k)})) $.
	For this, note that $ C_2 $ is either of type I or IV (since the 0'th homogeneous coordinate of the projection of $ \bs{w} $ along $ C_1 $ is positive) and thus $ \max_\prec C_2 = [\bs{w}_0+\bs{w}_1-1:0:\bs{w}_2+1:\bs{w}_3:\dots:\bs{w}_n] $ by direct computation (which uses the fact that the min-section to the projection along $ C_1 $ is monotone in this case, see the proof of Lemma \ref{lemma:parrallel_projection}).
	Moreover, its coefficient is $ \frac{1}{2} $ in the case that $ C_2 $ is of type I and $ 1 $ in the case that $ C_2 $ is of type IV.
\end{proof}

\begin{lemma}\label{lemma:patch_part_of_path}
	If $ \conv(C_1) \subseteq \Gamma_{\bs{{w}}} $ (in particular, $ C_1 $ is of type I), then there exist $ m(C) $ algebraic surfaces that pass through $ \overline{\boldsymbol{p}} $, tropicalize to $ S_0 $ and having a singularity that tropicalize to $ \bs{0} $. 
\end{lemma}

\begin{proof}
	The proof will work along similar lines to the proof of Lemma \ref{lemma:patch_not_part_of_path} with the only differnce being that the Jacobian at $ t=0 $ is not regular so we will need to calculate higher terms of $ \bs{q}_1 $.
	We will pick coordinates as before (i.e. $ \Sp(\overline{C}_h)=\Sp\{\bs{w}, \bs{e}^{(i)}-\bs{e}^{(j)} \; | \; i,j\le r_h \} $) with the additional property that $ \bs{w}=\bs{e}^{(1)} $, so $ C_1=\{ \bs{u}^{(-1)}, \bs{u}^{(1)} \} $ with $ \bs{u}^{(-1)}=\bs{0} $ and $ \bs{u}^{(1)}=2\bs{e}^{(1)} $ (see Example \ref{ex:type_I_circ}).
	We will assume, without loss of generality, that $ a_{\bs{u}^{(1)}}=1 $.
	
	The defining equation \eqref{eq:patchworking_def_polynomial} of the hypersurface becomes
	\[ f(z) = z_1^2+a_{\bs{w}}z_1+a_{\bs{u}^{(-1)}} + \sum_{\bs{v}\notin C_1\cup C_0}t^{\nu(\bs{v})}a_{\bs{v}}z^{\bs{v}} \]
	and so, the equations $ 0=f(\bs{q})=\der{z_1} f(\bs{q})=t^{-s'_k}f(\bs{p}^{(k)}) $ ($ s'_k $ is as in \eqref{eq:lemma_mult_point_cond}), where $ k $ is s.t. $ \bchi{k} $ lies on the $ (n-1) $-face dual to $ \conv(C_1) $, takes the form
	\begin{alignat}{2}\label{eq:lemma_mult_circ_in_path_level_1}
	\xi_{k,1}^2+a_{\bs{w}}\xi_{k,1}+a_{\bs{u}^{(-1)}}+\sum_{\bs{v}\notin C_1\cup C_0}t^{\nu(\bs{v})}a_{\bs{v}}\left (\bs{p}^{(k)}\right )^{\bs{v}} = && t^{-s'_k}f(\bs{p}^{(k)}) & = 0  \nonumber \\
	2\bs{q}_1 + a_{\bs{w}} + \sum_{\bs{v}\notin C_1\cup C_0}\bs{v}_1\cdot t^{\nu(\bs{v})}a_{\bs{v}} \bs{q}^{\bs{v}-\bs{e}^{(1)}} = && \der{z_1} f(\bs{q}) & = 0  \\
	\bs{q}_1^2 + a_{\bs{w}} \bs{q}_1 + a_{\bs{u}^{(-1)}} + \sum_{\bs{v}\notin C_1\cup C_0}t^{\nu(\bs{v})}a_{\bs{v}}\bs{q}^{\bs{v}} = && f(\bs{q}) & = 0 \nonumber
	\end{alignat}
	We can extract $ a_{\bs{w}} $ and $ a_{\bs{u}^{(-1)}} $ from the first 2 equations
	\begin{equation}\label{eq:lemma_patch_circ_in_path_extracted}
	\begin{split}
	a_{\bs{w}} &= -2\bs{q}_1-\sum_{\bs{v}\notin C_1\cup C_0}t^{\nu(\bs{v})}\bs{v}_1a_{\bs{v}}\bs{q}^{\bs{v}-\bs{e}^{(1)}}, \\
	a_{\bs{u}^{(-1)}} & = -a_{\bs{w}}\xi_{k,1}-\xi_{k,1}^2-\sum_{\bs{v}\notin C_1\cup C_0}t^{\nu(\bs{v})}a_{\bs{v}} \left (\bs{p}^{(k)} \right )^{\bs{v}} = \\
	& = 2\bs{q}_1\xi_{k,1}-\xi_{k,1}^2+\sum_{\bs{v}\notin C_1\cup C_0}t^{\nu(\bs{v})}a_{\bs{v}} \left (\bs{v}_1\xi_{k,1}\bs{q}^{\bs{v}-\bs{e}^{(1)}}- \left (\bs{p}^{(k)}\right )^{\bs{v}} \right )
	\end{split}
	\end{equation}
	and substitute it to the last equation of \eqref{eq:lemma_mult_circ_in_path_level_1} to get 
	\begin{equation}\label{eq:lemma_mult_circ_in_path_level_1_simpl}
	(\bs{q}_1-\xi_{k,1})^2 = \sum_{\bs{v}\notin C_1\cup C_0}a_{\bs{v}}t^{\nu(\bs{v})}\left[ (1-\bs{v}_1)\bs{q}^{\bs{v}} + \bs{v}_1\xi_{k,1}\bs{q}^{\bs{v}-\bs{e}^{(1)}}-\left (\bs{p}^{(k)}\right )^{\bs{v}} \right] = o(1)
	\end{equation}
	so $ \bs{q}_1(0) = \xi_{k,1} $. 
	
	Before we proceed and find both options for the next term of $ \bs{q}_1 $, we will note that by Lemma \ref{lemma:patch_high_levels}, there exist $ \prod_{h>1}m(C_h)=\frac{m(C)}{2} $ solutions of $ (\bs{q}_i(0))_{i>1} $.
	
	To find the next term of $ \bs{q}_1 $, we separate the terms in the right-hand side of \eqref{eq:lemma_mult_circ_in_path_level_1_simpl} to those with $ \bs{v}\in \Sp(C_1)\setminus C_1\cup C_0 $ and those with $ \bs{v}\notin \Sp(C_1) $.
	
	If $ \bs{v}\in \Sp(C_1) $ then $ \bs{v}=\bs{v}_1\cdot \bs{e}^{(1)} $ meaning that the corresponding term in \eqref{eq:lemma_mult_circ_in_path_level_1_simpl} contains the univariate polynomial (w.r.t. $ \bs{q}_1 $) $ (1-\bs{v}_1)\bs{q}_1^{\bs{v}_1}+\bs{v}_1\xi_{k,1}\bs{q}_1^{\bs{v}_1-1}-\xi_{k,1}^{\bs{v}_1} $.
	Note that $ \bs{q}_1=\xi_{k,1} $ is a root of both this polynomial and its derivative $ (1-\bs{v}_1)\bs{v}_1\bs{q}_1^{\bs{v}_1-1}+(1-\bs{v}_1)\bs{v}_1\xi_{k,1}\bs{q}_1^{\bs{v}_1-2} $, so there exists a polynomial $ g_{\bs{v}}\in \mbb{Q}[\bs{q}_1, \xi_{k,1}] $ s.t. 
	\[ (1-\bs{v}_1)\bs{q}_1^{\bs{v}_1}+\bs{v}_1\xi_{k,1}\bs{q}_1^{\bs{v}_1-1}-\xi_{k,1}^{\bs{v}_1} = (\bs{q}_1-\xi_{k,1})^2 g_{\bs{v}}. \]
	So we can rewrite \eqref{eq:lemma_mult_circ_in_path_level_1_simpl} as
	\begin{equation}\label{eq:lemma_mult_circ_in_path_level_1_simpl2}
	\begin{split}
	\left[ 1+\sum_{\bs{v}\in \Sp(C_1)\setminus C_1\cup C_0} a_{\bs{v}} t^{\nu(\bs{v})} g_{\bs{v}}  \right] (\bs{q}_1-\xi_{k,1})^2 =\sum_{\bs{v}\notin \Sp(C_1)} a_{\bs{v}}t^{\nu(\bs{v})} \left[ (1-\bs{v}_1)\bs{q}^{\bs{v}} + \bs{v}_1\xi_{k,1}\bs{q}^{\bs{v}-\bs{e}^{(1)}}-\left (\bs{p}^{(k)} \right )^{\bs{v}} \right].
	\end{split}
	\end{equation}
	
	Now consider terms with $ \bs{v}\notin \Sp(C_1) $, i.e. terms in the right-hand side of \eqref{eq:lemma_mult_circ_in_path_level_1_simpl2}.
	By rearranging, we get
	\[ (1-\bs{v}_1)\bs{q}^{\bs{v}} + \bs{v}_1\xi_{k,1}\bs{q}^{\bs{v}-\bs{e}^{(1)}}-\left (\bs{p}^{(k)}\right )^{\bs{v}} = \bs{q}^{\bs{v}}-\left (\bs{p}^{(k)}\right )^{\bs{v}} + \bs{v}_1(\xi_{k,1}-\bs{q}_1)\bs{q}^{\bs{v}-\bs{e}^{(1)}}=\bs{q}^{\bs{v}}-\left (\bs{p}^{(k)}\right )^{\bs{v}}+o(1). \]
	We thus can rewrite \eqref{eq:lemma_mult_circ_in_path_level_1_simpl2}, using Lemma \ref{lemma:patch_part_of_path_f_0_trop}, as 
	\begin{equation*}
	\begin{split}
	\left( 1+o(1) \right) \left( \bs{q}_1-\xi_{k,1} \right)^2 &= \sum_{\bs{v}\notin \Sp(C_1)}a_{\bs{v}}t^{\nu(\bs{v})}\left( \bs{q}^{\bs{v}}-\left (\bs{p}^{(k)} \right )^{\bs{v}}+o(1) \right) = \\
	= f_0(\bs{q})-f_0(\bs{p}^{(k)})\cdot (1+o(1)) &= a_{\bs{w}^{\pr}}t^{\nu_0(\bs{w}^{\pr})-M_k\pair{\bphi}{\bs{w}^{\pr}}}\xi_k^{\bs{w}^{\pr}}(1+o(1)).
	\end{split}
	\end{equation*}
	where $ \bs{w}^{\pr} $ as in Lemma \ref{lemma:patch_part_of_path_f_0_trop} and $ M_k,\bphi $ as in Construction \ref{con:Mikhalkin_position}.
	We thus get that $$ \val(\bs{q}_1-\xi_{k,1})=s:=\frac{\nu_0(\bs{w}^{\pr})-M_k\pair{\bphi}{\bs{w}^{\pr}} }{2} $$ and 
	\begin{equation}\label{eq:lemma_patch_in_path_res_level1}
	\bs{q}_1 = \xi_{k,1}\pm t^s\sqrt{a_{\bs{w}^{\pr}}\xi_k^{\bs{w}^{\pr}}}+o(t^s). 
	\end{equation}
	
	To finish the proof, consider the map (note the order of the variables and equations)
	\begin{align*}
	\Psi : \mbb{C}\times \mbb{C}^{N-2} \times \mbb{C}\times \mbb{C} \times \mbb{C}^{n-1} \times \mbb{C} &\to \mbb{C}^{N-2}\times \mbb{C} \times \mbb{C} \times \mbb{C}^{n-1}\times \mbb{C} \\
	\begin{pmatrix}
	t \\ (a_{\bs{v}})_{\bs{v}\in \Delta\cap\mbb{Z}^n\setminus (C_0\cup C_1)}  \\ a_{\bs{u}^{(-1)}} \\ a_{\bs{w}} \\ (\bs{q}_i)_{i>1} \\ \delta \bs{q}_1
	\end{pmatrix}  &\mapsto  \begin{pmatrix}
	(t^{-s'_j}f(\bs{p}^{(j)}))_{j\ne k} \\ t^{-s_k}f(\bs{p}^{(k)}) \\ t^{-s_1}\der{z_1}f|_{z=\bs{q}} \\  (t^{-s_i}\bs{q}_i\der{z_i}f|_{z=\bs{q}})_{i>1} \\ t^{-2s}f(q)  
	\end{pmatrix}
	\end{align*}
	where $ N:=|\Delta\cap\mbb{Z}^n|-1 $, $ \bs{q}_1=\xi_{k,1}+t^{s}\delta \bs{q}_1 $, $ s_i $ as in Lemma \ref{lemma:patch_high_levels} and $ s'_j $ as in \eqref{eq:lemma_mult_point_cond}.
	By the above calculations, $ \Psi^{-1}(\boldsymbol{0}) $ contains $ m(C) $ points with $ t=0 $ so, by the implicit function theorem, we are left to show that the Jacobian of $ \Psi $ with respect to all variables but $ t $ at those points is invertible.
	The only equation that depends on $ \delta \bs{q}_1 $ at $ t=0 $ is the last one and its derivative w.r.t. $ \delta \bs{q}_1 $ is $ 2\delta \bs{q}_1 \ne 0 $.
	After removing the last column and row from the Jacobian we are left with triangular block matrix with the blocks on the diagonal that correspond to $ C_h $ with $ h>1 $ are invertible by Lemma \ref{lemma:patch_high_levels}.
	Finally, the block that correspond to the coefficients $ a_{\bs{v}} $, the point conditions and $ \der{z_1}f $ is (suitably arranged) lower triangular matrix with non-zero entries on the diagonal by $ \eqref{eq:lemma_patch_circ_in_path_extracted} $ and thus invertible.
\end{proof}

\section{Count of algebraic singular hypersurfaces}\label{sec:results}

In this section we use the results of preceding sections to count the number of singular hypersurfaces satisfying point conditions. 
More precisely, this number is a polynomial in $ d $ of degree $ n $, we calculate its leading coefficient.
In Subsection \ref{sec:discriminant_degree} we work over $ \mbb{C} $.
Since in this setting the number of singular hypersurfaces in a pencil is a known invariant, this result can be viewed as a verification that we constructed almost all graded circuits admitting Mikhalkin condition.

In Sections \ref{sec:real_count} and \ref{sec:multinodal}, we work over $ \mbb{R} $.
As mentioned in the introduction, the number of real singular hypersurfaces in a pencil is not invariant and one can only obtain an existence result.
In Section \ref{sec:real_count} we prove Theorem \ref{thm:real_count_formula}, namely present a choice of points that yield $ \alpha_nd^n+O(d^{n-1}) $ real singular hypersurfaces where $ \alpha_n $ satisfy the recurrence \eqref{eq:lemma_real_count_coef_rec}.
In particular, for $ n\ge 3 $ we attain more than $ 2d^n $ real singular hypersurfaces, which, to the best of our knowledge, is the best known result.
In Section \ref{sec:multinodal} we prove Theorem \ref{thm:multinode_count} and give a lower bound for the maximal number of $ \delta $-nodal real hypersurfaces satisfying point conditions.

\subsection{Degree of discriminant}\label{sec:discriminant_degree}

\begin{lemma}\label{lemma:type_IV_count}
	For $ n>1 $, the sum of multiplicities of Mikhalkin graded circuits spanning $ \Delta $ and centered at $ \bs{w}=[\bs{w}_0:0:0:\dots:\bs{w}_n] $ is $ \bs{w}_0^{n-1}+O(\bs{w}_0^{n-2}) $.
\end{lemma}

\begin{proof}
	The only type of Mikhalkin circuits that can be centered at $ \bs{w} $ is type IV. 
	For $ \ell>0 $ define 
	\begin{align*}
	\pr_{\ell} : \mbb{Z}^{n+1} &\to \mbb{Z}^{n-\ell+1} \\
	[\bs{x}_0:\bs{x}_1:\dots:\bs{x}_n ] &\mapsto [ \sum_{i=0}^{\ell} \bs{x}_i:\bs{x}_{\ell+1}: \dots:\bs{x}_n].
	\end{align*}
	For a vector $ \bs{u}'\in \Delta^{(n-\ell)}\cap\mbb{Z}^{n-\ell+1} $ with $ \bs{u}'_{n-\ell}=\bs{w}_n-1 $, denote by $ \mathcal{C}_{\ell}(\bs{u}') $ the set of Mikhalkin graded circuits $ C\subseteq \Delta^{(n)} $, centered at $ \bs{w} $, with $ \Sp(C) = \Sp\{ \bs{w},\bs{u}, \bs{e}^{(i)}-\bs{e}^{(j)} \; | \; i,j\le \ell \} $ for some (and thus for all) $ \bs{u}\in \pr_{\ell}^{-1}(\bs{u}') $.
	Denote by $ s_{\ell}(\bs{u}) := \sum_{C\in\mathcal{C}_{\ell}(\bs{u})}m(C) $ the sum of their multiplicities.
	We will prove, by induction on $ \ell $ and $ n $, that $ s_{\ell}(\bs{u}') = \left( \bs{u}'_0 \right)^{\ell} + O((\bs{u}'_0)^{\ell-1}) $ which will obviously prove the lemma since the desired number is equal to $ s_{n-1}(\bs{u}') $ for $ \bs{u}'=[\bs{w}_0:\bs{w}_n] $.
	
	If $ \ell=1 $ then $ \mathcal{C}_{\ell}(\bs{u}') $ is just the collection of type IV circuits with $ \bs{u}_0+\bs{u}_1=\bs{u}'_0 $. By Example \ref{ex:type_IV_circ}, they are in bijection with $ \pr_{\ell}^{-1}(\bs{u}') $ and they are all of  multiplicity 1, so $ s_1(\bs{u}') = \bs{u}'_0 + O(1) $.
	
	Suppose now that $ \ell>1 $.
	Let $ \bs{v}'\in \Delta^{(n-\ell+1)} $ with $ \bs{v}'_1>0 $ and $ C' $ the unique type I circuit centered at $ \bs{v}' $.
	For any $ \widetilde{C}\in\mathcal{C}_{\ell-1}(\bs{v}') $, the glue $ C:=\widetilde{C}\Vert C' $ is in $ \mathcal{C}_{\ell}(\pr_1(\bs{v}')) $ and its multiplicity is double of $ m(\widetilde{C}) $.
	We thus get, using the induction hypothesis for $ \ell -1 $, that the contribution of those graded circuits to $ s_{\ell}(\bs{u}') $ is
	\begin{equation}\label{eq:lemma_type4_count_type1}
	\sum_{\bs{v}'\in \pr_1^{-1}(\bs{u}')}\sum_{\widetilde{C}\in\mathcal{C}_{\ell-1}(\bs{v}')} 2\cdot m(\widetilde{C})=\sum_{\bs{v}'_0=1}^{\bs{u}'_0}2\cdot ((\bs{v}'_0)^{\ell-1} + O((\bs{v}'_0)^{\ell-2})=\frac{2}{\ell}{\bs{u}'_0}^{\ell} + O((\bs{u}'_0)^{\ell-1}).
	\end{equation}
	
	Now, take $ 1\le \ell_0<\ell-1 $ and $ \bs{v}'\in \Delta^{(n-\ell_0)} $ with $ \bs{v}'_1=\dots=\bs{v}'_{\ell-\ell_0-1}=0 $ and $ \bs{v}'_{\ell-\ell_0}>0 $.
	For $ \widetilde{C}\in\mathcal{C}_{\ell_0}(\bs{v}') $ and $ C'\subseteq \Delta_{\bs{u}'_0}^{(\ell-\ell_0)} $ centered at $ (\bs{v}'_{\ell-\ell_0},0,\dots,0,\bs{v}'_0) $, Lemma \ref{lemma:gluing_type_IV} shows that $ C:=\widetilde{C}\Vert C' $ admits Mikhalkin condition and by definition $ m(C)=m(\widetilde{C})\cdot m(C') $.
	By the induction hypothesis for $ n=\ell-\ell_0 $, the sum of multiplicities over all possible choices for $ C' $ is $ {\bs{v}'_0}^{\ell-\ell_0-1} + O({\bs{v}'_0}^{\ell-\ell_0-2}) $. 
	We get, by induction hypothesis for $ \ell_0 $ (as $ \ell $, the dimension $ n $ is the same), that the contribution of such graded circuits to the monomial $ {\bs{u}'_0}^{\ell} $ in $ s_{\ell}(\bs{u}') $ is  
	\begin{equation}\label{eq:lemma_type4_count_type4}
	\begin{split}
	\sum_{\bs{v}'_0=1}^{\bs{u}'_0}\sum_{\widetilde{C}\in\mathcal{C}_{\ell_0}(\bs{v}')}m(\widetilde{C})\cdot\left ( {\bs{v}'_0}^{\ell-\ell_0-1} + O({\bs{v}'_0}^{\ell-\ell_0-2})\right )&=\sum_{\bs{v}'_0=1}^{\bs{u}'_0}\left ({\bs{v}'_0}^{\ell-\ell_0-1}+O({\bs{v}'_0}^{\ell-\ell_0-2})\right ) s_{\ell_0}(\bs{v}')=\\
	= \sum_{\bs{v}'_0=1}^{\bs{u}'_0}\left ({\bs{v}'_0}^{\ell-1} + O({\bs{v}'_0}^{\ell-2})\right )&=\frac{1}{\ell}{\bs{u}'_0}^{\ell}+O({\bs{u}'_0}^{\ell-1}) .
	\end{split}
	\end{equation}
	Combining those contributions we get 
	\begin{equation}\label{eq:lemma_type4_count_combine}
	s_{\ell}(\bs{u}')=\frac{2}{\ell}{\bs{u}'_0}^{\ell} + \sum_{\ell_0=1}^{\ell-2}\frac{1}{\ell}{\bs{u}'_0}^{\ell} + O({\bs{u}'_0}^{\ell-1})= {\bs{u}'_0}^{\ell} + O({\bs{u}'_0}^{\ell-1}).
	\end{equation}
\end{proof}

Before we proceed to count all singular algebraic hypersurfaces that were constructed in Sections \ref{sec:construction_graded_circuits} and \ref{sec:patchworking}, we will need the following combinatorial lemma:

\begin{lemma}\label{lemma:count_type_3_volume}
	For $ \bs{w}=[0:0:\dots:\bs{w}_{n-1}:\bs{w}_n]\in \Delta $, there exist $ \bs{w}_{n-1}^{n-1}+O(\bs{w}_{n-1}^{n-2}) $ elementary Mikhalkin circuits of type III, lower index 0 and upper index $ n $ centered at $ \bs{w} $.
\end{lemma}

\begin{proof}
	The lattice volume of $ F=\{ \bs{x}\in\Delta \; | \; \bs{x}_n=\bs{w}_n+1 \} $ is $ (\bs{w}_{n-1}-1)^{n-1}=\bs{w}_{n-1}^{n-1}+O(\bs{w}_{n-1}^{n-2}) $.
	Let $ \sigma $ be the Mikhalkin triangulation of $ F $.
	Every $ (n-1) $-dimensional cell of $ \sigma $ is regular lattice simplex and thus have lattice volume equal to 1, meaning that $ \sigma $ consist of $ \bs{w}_{n-1}^{n-1}+O(\bs{w}_{n-1}^{n-2}) $ maximal dimensional cells.
	The number of cells that contain a vertex with $ \bs{x}_{n-1}\ge \bs{w}_{n-1}-2 $ is bounded by a constant independent of $ \bs{w}_{n-1} $ and thus do not affect the asymptotics. 
	
	We are left to show that the number of cells that contain a proper face which spans $ \bs{e}^{(n-1)}-\bs{e}^{(n-2)} $ is $ O(\bs{w}_{n-1}^{n-2}) $.
	Note that any such proper face contained in a facet of the $ (n-1) $-cell we are considering, which in turn is an $ (n-2) $-cell of $ \sigma $.
	Moreover, every $ (n-2) $-cell of $ \sigma $ is contained in at most two $ (n-1) $-cells of $ \sigma $, so it is sufficient to show there are $ O(\bs{w}_{n-1}^{n-2}) $ $ (n-2) $-cells of $ \sigma $ that span $ \bs{e}^{(n-1)}-\bs{e}^{(n-2)} $.
	
	Let $ \rho $ be an $ (n-2) $-cell of $ \sigma $ that spans linearly $ \bs{e}^{(n-1)}-\bs{e}^{(n-2)} $.
	Considering the construction of $ \sigma $ in lemma \ref{lemma:mikhalkin_triangulation}, we see that $ \rho=\conv(\rho_1\cup \rho_2) $ where $ \rho_1 $ and $ \rho_2 $ are cells of $ \sigma $ that contained in
	\[ \rho_1\subseteq \left\{ \bs{x}\in\Delta \; | \; \bs{x}_0=\bs{x}_1=\dots=\bs{x}_{n-r-3}=0, \bs{x}_{n-1}=a, \bs{x}_n=\bs{w}_n+1 \right\} =: F_1 \]
	and 
	\[ \rho_2\subseteq \left\{ \bs{x}\in\Delta \; | \; \bs{x}_{n-r-1}=\bs{x}_{n-r}=\dots=\bs{x}_{n-2}=0, \bs{x}_{n-1}=a+1, \bs{x}_n=w_n+1 \right\} =: F_2  \]
	for some $ 0\le r \le n-3 $ and $ a\in\mbb{Z}_{\ge 0} $.
	Since $ \dim \rho_1 + \dim \rho_2 = n-3 $, either $ \dim \rho_1=r $ and $ \dim \rho_2=n-3-r $ or $ \dim \rho_1 = r-1 $ and $ \dim \rho_2 = n-2-r $.
	In the former case it follows that the last vertex of $ F_2 $ contained in the span of $ \rho_2 $, while in the latter case the first vertex of $ F_1 $ is contained in the span of $ \rho_1 $.
	It is thus sufficient to prove the following lemma.
\end{proof}

\begin{lemma}
	Let $ \rho $ be an $ (n-1) $-cell of the Mikhalkin triangulation of $ n $ dimensional simplex $ \Delta^{(n)} $, and assume that either the first or the last vertex of $ \Delta^{(n)} $ is contained in the span of $ \rho $.
	Then $ \rho\subseteq \partial\Delta $.
\end{lemma}

\begin{proof}
	The proof is by induction on $ n $, with the case of $ n=1 $ being trivial.
	Let $ n>1 $.
	If the $ n $'th coordinate is constant on $ \rho $ then it has to be equal to 0, meaning that $ \rho\subseteq \partial \Delta $.
	Otherwise, write $ \rho=\conv(\rho_1\cup \rho_2) $ with 
	\[ \rho_1\subseteq \left\{ \bs{x}\in\Delta \; | \; \bs{x}_0=\bs{x}_1=\dots=\bs{x}_{n-r-2}=0, \bs{x}_{n}=a \right\} =: F_1 \]
	and 
	\[ \rho_2\subseteq \left\{ \bs{x}\in\Delta \; | \; \bs{x}_{n-r}=\dots=\bs{x}_{n-1}=0, \bs{x}_{n}=a+1 \right\} =: F_2  \]
	for some $ 0\le r\le n-2 $ and $ a\in\mbb{Z}_{\ge 0} $.
	As in the proof of Lemma \ref{lemma:count_type_3_volume}, we will deal separately with the case $ \dim \rho_1=r $, $ \dim \rho_2=n-2-r $ and with the case $ \dim \rho_1=r-1 $, $ \dim \rho_2 = n-1-r $.
	
	As before, if the first condition holds, then the last vertex of $ F_2 $ is in the span of $ \rho_2 $.
	Since the restriction of the Mikhalkin triangulation of $ \Delta^{(n)} $ to $ F_2 $ is the same as the Mikhalkin triangulation of $ F_2 $ when viewed as $ \Delta^{(n-1-r)}_{d-a-1} $, by the induction hypothesis, $ \rho_2\subseteq \partial F_2 $, which implies that $ \rho\subseteq \partial \Delta $.
	
	Similarly, if $ \dim \rho_1 = r-1 $, then the first vertex of $ F_1 $ is in the span of $ \rho_1 $, which means that $ \rho_1\subseteq \partial F_1 $ and thus $ \rho\subseteq \partial \Delta $.
\end{proof}

\begin{prop}\label{prop:discriminant_degree}
	Let $ \overline{\bs{p}}\subseteq \mbb{P}^n $ with $ \trop(\overline{\bs{p}}) = \xpoints $. 
	Then almost all singular hypersurfaces in $ \sing(\Delta, \overline{\bs{p}}) $ have tropicalization with Mikhalkin graded circuit. That is, the ratio of those singular hypersurfaces passing through the points $ \overline{\bs{p}} $ and not satisfying this condition tends to 0 as $ d\to \infty $.
\end{prop}

\begin{proof}
	We will show, by induction on $ n $, that the sum of multiplicities of all the graded circuits is $ (n+1)d^n + O(d^{n-1}) $ which is asymptotically equal to the classically known degree of the discriminant $ (n+1)(d-1)^n $ (see, for example, \cite[Chapter 9, Example 2.10(a)]{Gelfand_kapranov}). 
	
	Let $ h < ht(C) $ be the maximal level s.t. $ \partial \overline{C}_h = \Sp\{ \bs{e}^{(i)}-\bs{e}^{(j)} \; | \; \ell_1\le i,j \le \ell_2 \} $ for some $ 0 \le \ell_1 < \ell_2 \le n $. 
	Denote $ \widetilde{C}:=\overline{C}_h $, $ \pr $ the projection along $ \widetilde{C} $, $ \bs{w}':=pr(\bs{w}) $ and $ C':=\pr(C) $.
	From the maximality of $ h $, $ C' $ is elementary Mikhalkin graded circuit (see Definition \ref{def:elementary_Mikhalkin_graded_circuit}), and we are in the situation of Lemma \ref{lemma:parrallel_projection}. We will separate into cases according to the type of $ C'_1 $:
	\begin{description}
		\item[Type I] In this case $ \dim C'=1 $ and $ \widetilde{C}\subseteq \Delta_{\bs{w}'_0}^{(n-1)} $. By the induction hypothesis, we have 
		\[ \sum_{C'_1\text{ of type I}}m(C)=\sum_{\bs{w}'_0=1}^{d}\sum_{\widetilde{C}\subseteq \Delta_{\bs{w}'_0}^{(n-1)}}2\cdot m(\widetilde{C})=2\sum_{\bs{w}_0'=1}^d\left [ n{\bs{w}_0'}^{n-1}+O({\bs{w}_0'}^{n-2})\right ]=2d^n+O(d^{n-1}). \]
		
		\item[Type II] We have $ \dim C'=2 $ and $ \widetilde{C}\subseteq \Delta_{\bs{w}'_1}^{(n-2)} $. By the induction hypothesis:
		\[ \sum_{C'_1\text{ of type II}}m(C)=\sum_{\bs{w}'_1=1}^{d}\sum_{\widetilde{C}\subseteq \Delta_{\bs{w}'_1}^{(n-2)}} m(\widetilde{C})=\sum_{\bs{w}_1'=1}^d \left[ (n-1){\bs{w}_1'}^{n-2}+O({\bs{w}_1'}^{n-3}) \right] =O(d^{n-1}). \]
		meaning that circuits of type II do not contribute to the asymptotic count. 
		
		\item[type III] Denote by $ r $ the dimension of $ C' $, so $ \widetilde{C}\subseteq \Delta_{\bs{w}_{r-1}'}^{(n-r)} $. The number of type III circuits of dimension $ r $ centered at $ \bs{w}':=[0:0:\dots:0:\bs{w}'_{r-1}:\bs{w}'_r] $ is $ \bs{w}_{r-1}'^{r-1}+O(\bs{w}_{r-1}'^{r-2}) $ by Lemma \ref{lemma:count_type_3_volume}.
		We get that
		\begin{align*} 
		\sum_{C'_{1}\text{ of type III}}m(C)=\sum_{r=2}^{n}\sum_{\bs{w}'_{r-1}=1}^{d}\sum_{\widetilde{C}\subseteq \Delta_{\bs{w}'_{r-1}}^{(n-r)}} m(\widetilde{C})=\\
		= \sum_{r=2}^n \sum_{\bs{w}_{r-1}'=1}^d \left [{\bs{w}_{r-1}'}^{r-1}\cdot (n-r+1){\bs{w}'_{r-1}}^{n-r}+O({\bs{w}'_{r-1}}^{n-2}) \right ] =\\
		= \sum_{r=2}^n (n-r+1)\sum_{\bs{w}'_{r-1}=1}^d \left [{\bs{w}'_{r-1}}^{n-1}+O({\bs{w}'_{r-1}}^{n-2}) \right ] =\\
		\frac{n(n-1)}{2}\cdot \frac{1}{n}d^n + O(d^{n-1})=\frac{n-1}{2}d^n+O(d^{n-1})
		\end{align*}
		
		\item[Type IV] Denote by $ r $ the dimension of $ C' $, so $ \widetilde{C}\subseteq \Delta_{\bs{w}'_0}^{(n-r)} $. By Lemma \ref{lemma:type_IV_count}, the sum of multiplicities of dimension $ r $ graded circuits centered at $ \bs{w}' $ with $ C'_1 $ of type IV that span $ \Delta^{(r)} $ is $ {\bs{w}'_0}^{r-1} + O{\bs{w}'_0}^{r-2}) $. Thus, by induction hypothesis,
		\begin{align*}
		\sum_{C'_1\text{ of type IV }}m(C) = \sum_{r=2}^{n}\sum_{\bs{w}'_0=1}^{d}\left( {\bs{w}'_0}^{r-1} \cdot (n-r+1){\bs{w}'_0}^{n-r} + O({\bs{w}'_0}^{n-1})\right ) = \frac{n-1}{2}d^n+O(d^{n-1})
		\end{align*}
	\end{description}
	
	Summing the contribution of different types gives us 
	\[ \sum_{\Sp(C)=\Delta^{(n)}}m(C) = 2d^n+\frac{n-1}{2}d^n+\frac{n-1}{2}d^n + O(d^{n-1})=(n+1)d^n + O(d^{n-1}). \]
\end{proof}

\subsection{Enumeration of real uninodal hypersurfaces}\label{sec:real_count}

\newcommand{\ind}{\text{ind}}

In this subsection we will compute the number of real singular hypersurfaces (i.e. defined over $ \mbb{R} $ and having a real singularity) for a particular arrangement of points.
Different from the complex case, the real count is dependent on the specific choice of the point conditions.
Specifically, it depends on the signs of $ \xi $, see the notation of Section \ref{sec:patchworking}.

We will present a choice of signs that ensure all the solutions for the equations appearing in the first level of the graded circuits are real.
We then show that with this choice of signs, one can expect one real solution on average for each successive levels, and thus get $ m(C_1) $ real singular hypersurfaces on average. 
We performed only initial investigation of the effect of the signs of $ \xi_{i,j} $ on the number of real singular hypersurfaces we obtain.
We believe a better result can be obtained using the technique described in this work.

\begin{Def}\label{def:real_multiplicity_with_xi}
	For a choice $ \xi_{j,i}\in\{\pm 1 \} \; (j=1,\dots,N-1; i=1,\dots,n) $ and a Mikhalkin graded circuit $ C $, denote by $ mt(C, (\xi_{j,i})) $ the number of real singular hypersurfaces passing through the points $ \xi_j t^{-M_j\bphi}\in (\mbb{K}_\mbb{R}^\times)^n $ for every $ j=1,\dots,N-1 $ and having a singularity tropicalizing to the tropical singular point corresponding to $ C $.
\end{Def}

It would be convenient to index the points $ \xi $ by their occurrence in the lattice path, which is the purpose of the following notation.

\begin{Def}
	For $ \bs{v}\in\Delta\cap\mbb{Z}^{n+1} $, denote by $ \ind(\bs{v}) $ its index in the order $ \prec $ (starting to count from $ 1 $), i.e. 
	\[ \ind(\bs{v}) = |\{\bs{u}\in\Delta\cap\mbb{Z}^{n+1} \; | \; \bs{u}\preceq \bs{v} \}|. \]
\end{Def}

\begin{Def}
	Denote by $ \Omega(\Delta) $ the set of choices for the signs of $ \xi $ s.t. for all $ C\in\mgc(\Delta) $ with $ C_1 $ of type I, the equations corresponding to the first level (either \eqref{eq:lemma_mult_circ_in_path_level_1_simpl} if $ \conv(C_1)\subseteq \Gamma $ or the discriminantal equation $ a_{\bs{w}}^2-4a_{\bs{u}^{(-1)}}a_{\bs{u}^{(1)}} $ if $ \conv(C_1) \nsubseteq \Gamma $) have two real solutions.
\end{Def}

\begin{lemma}\label{lemma:real_count_average}
	The set $ \Omega(\Delta) $ is not empty and	 the average of $ mt(C, (\xi_{j,i})) $ over all $ (\xi_{j,i})\in\Omega(\Delta) $ is $ m(C_1) $ for all $ C\in \mgc(\Delta) $.
\end{lemma}

\begin{proof}
	Since we will need the explicit conditions of Lemmas \ref{lemma:parrallel_projection} and \ref{lemma:gluing_type_IV}, it will be more convenient to work in homogeneous coordinates similar to Section \ref{sec:construction_graded_circuits} and not the affine coordinates of Section \ref{sec:patchworking}.
	To this end, set $ \xi_{j,0}=1 $ for all $ j=1,\dots,N-1 $.
	
	We begin by exploring the conditions to be an element of $ \Omega(\Delta) $.
	By changing to homogeneous coordinates for the coefficients of the polynomial (i.e. removing the assumption that $ {a_{\bs{u}^{(1)}}=1} $), we get that \eqref{eq:lemma_patch_in_path_res_level1} have 2 real solutions if
	\[ a_{\bs{u}^{(1)}}a_{\bs{w}^{\pr}}\xi_k^{\bs{w}^{\pr}-\bs{u}^{(-1)}}>0 \]
	where $ \bs{w}^{\pr} $ is, as in Lemma \ref{lemma:patch_part_of_path}, the last point preceding $ \bs{w} $ and not contained in $ \Sp(C_1) $.
	Calculating the dependence of the coefficients $ (a_{\bs{v}})_{\bs{v}\in \Delta\cap\mbb{Z}^n} $ in $ \xi $ by \eqref{eq:lemma_mult_point_cond}, we get that in order for \eqref{eq:lemma_mult_circ_in_path_level_1_simpl} to have two real solutions for all its possible occurrences, we need
	\begin{equation}\label{eq:lemma_real_count_level1_in_path_cond}
	(-1)^{\ind(\bs{u})-\ind(\bs{w}^{\pr})}\xi_{\ind(\bs{w}^{\pr})}^{\bs{v^{(p+1)}}-\bs{w}^{\pr}}\cdot \xi_{\ind(\bs{u})}^{\bs{w}^{\pr}-\bs{u}} \cdot \prod_{\ind(\bs{w}^{\pr})<j<\ind(\bs{u})} \xi_{j,1}>0
	\end{equation}
	for all points $ \bs{u}\in \Delta\cap\mbb{Z}^n $ with $ \bs{u}_0>1 $ ($ \bs{u} $ here acts as $ \bs{u}^{(-1)} $), where $ \bs{w}^{\pr} $ is the last point preceding $ \Sp(\bs{u}, \bs{e}^{(1)}-\bs{e}^{(0)}) $ and $ \bs{v^{(p+1)}} $ is its immediate successor, i.e. the first point of $ \Sp(\bs{u},\bs{e}^{(1)}-\bs{e}^{(0)})\cap \Delta\cap\mbb{Z}^{n+1} $.
	Note that none of the conditions \eqref{eq:lemma_real_count_level1_in_path_cond} depend on $ \xi_{\ind(\bs{v})} $ for $ \bs{v}\in\Delta\cap\mbb{Z}^{n+1} $ with $ \bs{v}_0=1 $.
	
	Next, let $ \bs{w}\in\Delta\cap\mbb{Z}^{n+1} $ with $ \bs{w}_i=0 $ for all $ i<\ell $ and $ w_{\ell}>0 $ where $ \ell>0 $ and let $ {C=\{ \{\bs{w}\}, \{\bs{u}^{(-1)}, \bs{u}^{(1)}\}\}} $ be a type I Mikhalkin elementary circuit of index $ \ell $ centered at $ \bs{w} $.
	The discriminantal equation has 2 real solutions if $ a_{\bs{u}^{(-1)}}a_{\bs{u}^{(1)}}>0 $ where \linebreak $ {\bs{u}^{(-1)}=\bs{w}-\bs{e}^{(\ell+1)}+\bs{e}^{(\ell)}} $ and $ \bs{u}_1=\bs{w}+\bs{e}^{(\ell+1)}-\bs{e}^{(\ell)} $. 
	Writing the dependence in $ \xi $ we get that once we set the first $ \ell-1 $ coordinates of $ \xi_j $ for all $ \ind(\bs{u}^{(-1)}) \le j < \ind(\bs{u}^{(1)})-1  $ and set $ \xi_{\ind(\bs{u}^{(-1)}), \ell} $, we have a unique sign for $ \xi_{\ind(\bs{w})-1, \ell} $ that ensures $ a_{\bs{u}^{(-1)}}a_{\bs{u}^{(1)}} $ would be positive.
	
	We will show those conditions can be satisfied simultaneously and deduce that $ \Omega(\Delta) $ is not empty.
	Indeed, we can first set $ \xi_{\ind(\bs{u})} $ for all $ \bs{u}\in\Delta\cap\mbb{Z}^n $ with $ \bs{u}_0=0 $ to be arbitrary.
	Next, set $ \xi_{\ind(\bs{u})} $ for $ \bs{u}\in\Delta\cap\mbb{Z}^n $ with $ \bs{u}_0>1 $ s.t. \eqref{eq:lemma_real_count_level1_in_path_cond} will hold.
	Finally, set $ \xi_{\ind(\bs{u})-1, \ell} $ for $ \bs{u} $ with $ \bs{u}_0=\dots=\bs{u}_{\ell-1}=0 $ and $ \bs{u}_{\ell}>0 $, inductively on $ \ell $, s.t. the discriminantal equation will have two real solutions.
	
	Now, fix $ C\in \mgc(\Delta) $.
	We will prove, by induction on $ h\ge 0 $ that the average number of solutions $ (a_{\bs{v}})_{\bs{v}\in \Delta},(\theta_i)_{i\le r_h} $ to the system \eqref{eq:lemma_patch_high_levels_big_system}, when averaging over $ \Omega(\Delta) $ is $ m(C_1) $.
	The basis case $ h=1 $ follows either from the definition of $ \Omega(\Delta) $ (if $ C_1 $ is of type I) or from the fact that whenever we have unique solution to real system of equations it has to be real (if $ m(C_1)=1 $).
	
	For $ h>1 $, if $ m(C_h)=1 $, then every real solution of \eqref{eq:lemma_patch_high_levels_small_system} corresponds to exactly one solution of \eqref{eq:lemma_patch_high_levels_big_system} so the result follows immediately from the induction hypothesis.
	To deal with the case when $ C_h $ is of type I, note that it can happen either in the conditions of Lemma \ref{lemma:parrallel_projection} or in the conditions of Lemma \ref{lemma:gluing_type_IV}.
	If the former holds, $ {\Sp(\overline{C}_{h-1})=\Sp(\bs{w},\bs{e}^{(i)}-\bs{e}^{(j)}\; | \;\ell\le i,j\le r)} $ for some $ 0<\ell<r $ and $ \bs{w}\in \Delta\cap\mbb{Z}^{n+1} $ and $ C_h=\{\bs{u}^{(-1)}, \bs{u}^{(1)}\} $ with $ \bs{u}^{(-1)} $ preceding $ \Sp(\overline{C}_{h-1}) $ and $ \bs{u}_1 $ succeeding it.
	Denote by $ \bs{u}_{\max} $ the last point of $ \Sp(\overline{C}_{h-1})\cap\Delta\cap\mbb{Z}^{n+1} $ and by $ \bs{u}' $ its immediate successor.
	Then, changing the sign of $ \xi_{\ind(\bs{u}_{\max})-1}^{\bs{u}'-\bs{u}_{\max}} $ would change the sign of $ a_{\bs{u}^{(-1)}}a_{\bs{u}^{(1)}} $ without affecting $ a_v $ for all $ \bs{v}\in \Sp(\overline{C}_{h-1}) $ meaning that it would change the sign under the radical in \eqref{eq:lemma_patch_high_levels_typeI_sol}.
	Unfortunately, only changing the sign of $ \xi_{\ind(\bs{u}_{\max})}^{\bs{u}'-\bs{u}_{\max}} $ would take us outside of $ \Omega(\Delta) $.
	To fix it, first note that we can change the sign of $ \xi_{\ind(\bs{u}_{\max})-1}^{\bs{u}'-\bs{u}_{\max}} $ without changing the sign of  $ \xi_{\ind(\bs{u}_{\max})-1}^{\bs{u}'-(\bs{u}_{\max}-\bs{e}^{(\ell)}+\bs{e}^{(\ell-1)})} $ and thus preserve the sign of the discriminantal equation corresponding to $ \bs{w}=\bs{u}_{\max} $.
	We may violated some of the conditions that correspond to type I circuits with index higher that $ \ell $, but they can be fixed by changing only coordinates higher that $ \ell $ of $ \xi_{\ind(\bs{v})} $ where $ \bs{v}_0=1 $, which fixes $ a_{\bs{v}} $ for all $ \bs{v}\in \Sp(\overline{C}_{h-1}) $ and does not affect \eqref{eq:lemma_real_count_level1_in_path_cond}. 
	We thus showed that there exists an involution of $ \Omega(\Delta) $ without fixed points that change the sign under the radical of \eqref{eq:lemma_patch_high_levels_typeI_sol}.
	This involution give us a partition of $ \Omega(\Delta) $ into pairs s.t. for every real solution of the system \eqref{eq:lemma_patch_high_levels_small_system} exactly one element in every pair has two real solutions of the system \eqref{eq:lemma_patch_high_levels_big_system}, meaning that, using the induction hypothesis, the average number of solutions to \eqref{eq:lemma_patch_high_levels_big_system} is $ m(C_1) $.
	
	Suppose now that we are in the conditions of Lemma \ref{lemma:gluing_type_IV}.
	Then \linebreak $ {\Sp(\overline{C}_{h-1}) = \Sp\{ \bs{w}, \bs{u}, \bs{e}^{(i)}-\bs{e}^{(j)} \; | \; \ell\le  i,j\le r \}} $ for some $ 0<\ell<r $, $ \bs{u}_r>0 $ and $ \bs{w}_i=0 $ for all $ i<\ell $ and for all $ \ell< i < r' $ where $ r'>r $.
	The proof of Lemma \ref{lemma:gluing_type_IV} shows that if we write $ C_h=\{ \bs{u}^{(-1)}, \bs{u}^{(1)} \}  $, then $ \bs{u}^{(1)} $ is the first point succeeding $ \Sp\{\bs{u},\bs{e}^{(i)}-\bs{e}^{(j)}\; | \; \ell \le i,j\le r\} $ with $ \bs{u}^{(1)}_i=0 $ for all $ i<\ell $.
	Denote the last point of $ \Sp\{ \bs{u}, \bs{e}^{(i)}-\bs{e}^{(j)} \; | \; \ell \le i,j\le r \} $ by $ \bs{u}_{\max} $ and its immediate successor by $ \bs{u}' $.
	We work in a similar way to the previous case and change the sign of $ \xi_{\ind(\bs{u}_{\max})}^{\bs{u}'-\bs{u}_{\max}} $ which would change the sign of $ a_{\bs{u}^{(1)}} $.
	To get an element of $ \Omega(\Delta) $, we first change the signs of $ \xi_{j,r+1} $ for all $ \ind(\bs{u}_{\max}) \le j \le \ind(\bs{u}_{\max})+(\bs{u}_{\max})_{r}-3 $, which ensures \eqref{eq:lemma_real_count_level1_in_path_cond} holds.
	Note that this would not change $ a_{\bs{u}^{(1)}} $ since its only dependent on the first coordinate of the corresponding $ \xi $.
	Next we might need to change $ \xi_{\ind(\bs{v}), r+1} $ for some $ \bs{v}\in\Delta\cap\mbb{Z}^{n+1} $ with $ \bs{v}_0=1 $ and $ \bs{v}_r=w_r-1 $.
	This can result in change of sign for $ a_{\bs{w}+\bs{e}^{(\ell+1)}-\bs{e}^{(\ell)}} $ that we can fix (in order to not change the signs of the $ \bs{q}_i $'s ) by changing the sign of some coordinate of $ \xi_{\ind(\bs{w})-2} $ without changing the sign of the power of $ \xi_{\ind(\bs{w})-2} $ that appear in any discriminantal equation.
	Finally, we change signs of high enough coordinates of $ \xi_{\ind(\bs{v})} $ for $ \bs{v}\in\Delta\cap\mbb{Z}^{n+1} $ with $ \bs{v}_0=1 $ to fix discriminantal sign conditions that correspond to level 1 circuits of type I and index higher than $ r' $.
	
\end{proof}

We can now prove Theorem \ref{thm:real_count_formula} which gives a recursive formula for a lower bound on the number of real singular hypersurfaces of given degree, satisfying point conditions. 

\begin{proof}[Proof of Theorem \ref{thm:real_count_formula}]
	We prove the assertion similarly to Lemma \ref{lemma:type_IV_count} and Proposition \ref{prop:discriminant_degree}.
	In both cases, the computation is similar to the one over the complex numbers, but with different coefficients. 
	For $ r>0 $ and $ \bs{u}'\in\Delta^{(n-r)}\cap\mbb{Z}^{n-r+1} $ denote, using the notation of Lemma \ref{lemma:type_IV_count}, \linebreak $ {s_r^{\mbb{R}}(\bs{u}'):=\sum_{C\in \mathcal{C}_r(\bs{u}')}m(C_1)} $.
	Since all $C\in \mathcal{C}_r(\bs{u}') $ have $C_1$ of type IV, $m(C_1)=1$ so \linebreak ${s_r^{\mbb{R}}(\bs{u}') =|\mathcal{C}_r(\bs{u}') |}$. 
	The recurrence relation for $ s_r^{\mbb{R}}(\bs{u}') $, analogously to \eqref{eq:lemma_type4_count_type1} and \eqref{eq:lemma_type4_count_type4}, is
	\begin{equation}\label{eq:lemma_real_count_type4_rec_s}
	s_r^{\mbb{R}}(\bs{u}') = \sum_{\bs{v}'\in pr_1^{-1}(\bs{u}')} s_{r-1}^{\mbb{R}}(\bs{v}') + \sum_{r_0=1}^{r-2} \sum_{\bs{v}'\in pr_{r_0}^{-1}(\bs{u}')} s_{r_0}^{\mbb{R}}(\bs{v}')\cdot  s_{r-r_0-1}^{\mbb{R}}([\bs{v}'_0:\bs{v}'_{r-r_0}])
	\end{equation}
	with $ s_1^{\mbb{R}}(\bs{u}')=s_1(\bs{u}')=\bs{u}'_0+O(1) $.
	It is easy to see that $ s_r^{\mbb{R}}(\bs{u}') $ is a polynomial in $ \bs{u}'_0 $ of degree $ r $, denote by $ \beta_{r} $ its leading coefficient.
	We thus get from \eqref{eq:lemma_real_count_type4_rec_s} a recurrence relation for $ \beta_{r} $,
	\begin{equation}\label{eq:lemma_real_count_type4_rec_coef}
	\beta_r = \frac{1}{r}\left( \beta_{r-1} + \sum_{r_0=1}^{r-2}\beta_{r_0}\cdot \beta_{r-r_0-1} \right) .
	\end{equation}
	The initial condition is given by 
	\[ s_1^{\mbb{R}}(\bs{u}')=s_1(\bs{u}')=\bs{u}'_0+O(1) \]
	meaning that $ \beta_1=1 $.
	
	We now turn to prove the analog of Proposition \ref{prop:discriminant_degree}.
	We will show, by induction on $ n $, that the sum of $ m(C_1) $ over all Mikhalkin graded circuits $ C $ is $ \alpha_nd^n +O(d^{n-1}) $ and that $ \alpha_n $ satisfies the recurrence relations \eqref{eq:lemma_real_count_coef_rec}.
	By Lemma \ref{lemma:real_count_average}, this would imply that the average of $ \sum_{C\in\mgc(\Delta)}mt(C, (\xi_{j,i})) $ over $ (\xi_{j,i})\in \Omega(\Delta) $ is $ \alpha_nd^n+O(d^{n-1}) $ meaning that there exists $ (\xi_{j,i})\in\Omega(\Delta) $ with $ \sum_{C\in\mgc(\Delta)}mt(C, (\xi_j,i))\ge \alpha_nd^n+O(d^{n-1}) $.
	
	Similarly to the proof of Proposition \ref{prop:discriminant_degree}, we pick $ h<ht(C) $ to be maximal s.t. \linebreak $ {\partial \overline{C}_h=\Sp\{ \bs{w}, \bs{e}^{(i)}-\bs{e}^{(j)} \; | \; \ell_1\le i,j\le \ell_2 \}} $ for some $ 0\le \ell_1<\ell_2\le n $.
	Denote $ \widetilde{C}:=\overline{C}_h $ and $ C':=\pr(C) $ where $ \pr $ is a projection along $ \widetilde{C} $.
	Next, split into cases according to the type of $ C'_1 $.
	For $ C'_1 $ of type I, we get, using the induction hypothesis,
	
	\begin{align*}		
	\sum_{C'_1\text{ is of type I}}m(C_1) & = \sum_{\bs{w}_0'=1}^d \sum_{\widetilde{C}\in\mgc(\Delta_{\bs{w}'_0}^{(n-1)})}m(\widetilde{C}_1) = \\
	= \sum_{\bs{w}_0'=1}^d \left[\alpha_{n-1}{\bs{w}_0'}^{n-1}+O(d^{n-2})\right]  & = \frac{\alpha_{n-1}}{n}d^n+O(d^{n-1}).
	\end{align*}
	
	Since the number of real solutions is at most the number of complex solutions, type II circuits do not contribute to the asymptotic.
	For $ C'_1 $ of type III, using Lemma \ref{lemma:count_type_3_volume}, 
	\begin{align*} 
	\sum_{C'_{1}\text{ of type III}}m(C_1)&=\sum_{r=2}^{n}\sum_{\bs{w}'_{r-1}=1}^{d}\sum_{\widetilde{C}\subseteq \Delta_{\bs{w}'_{r-1}}^{(n-r)}} m(\widetilde{C}_1)= \\
	&= \sum_{r=2}^n \sum_{\bs{w}_{r-1}'=1}^d \left [{\bs{w}_{r-1}'}^{r-1}\cdot \alpha_{n-r}{\bs{w}'_{r-1}}^{n-r}+O({\bs{w}'_{r-1}}^{n-2})\right ] = \\
	&= \sum_{r=2}^n \alpha_{n-r}\sum_{\bs{w}'_{r-1}=1}^d \left [{\bs{w}'_{r-1}}^{n-1}+O({\bs{w}'_{r-1}}^{n-2}) \right ]= \frac{1}{n}\left ( \sum_{r=0}^{n-2}\alpha_{r} \right) d^n+O(d^{n-1}).
	\end{align*}
	
	For $ C'_1 $ of type IV, we get
	\begin{align*}
	\sum_{C'_1\text{ of type IV }}m(C_1) =\sum_{r=2}^{n}\sum_{\bs{w}'_0=1}^{d} s_{r-1}^{\mbb{R}}(\bs{w}')\cdot\sum_{\widetilde{C}\in\mgc(\Delta_{\bs{w}_0'}^{(n-r)})}m(\widetilde{C}_1) = \\
	\sum_{r=2}^{n}\sum_{\bs{w}'_0=1}^{d}\left [\beta_{r-1}{\bs{w}'_0}^{r-1}\cdot \alpha_{n-r}{\bs{w}'_0}^{n-r} + O(d^{n-2})\right ] =\frac{1}{n} \left( \sum_{r=0}^{n-2}\beta_{n-1-r}\alpha_r \right)d^n+O(d^{n-1}). 
	\end{align*}
	
	Combining the contributions from the different types we get exactly \eqref{eq:lemma_real_count_coef_rec}.
	We thus got the desired result for $ \mbb{K}_\mbb{R} $. 
	By Tarski's transfer principle \cite[Theorem 1.16]{Jensen_Lenzing_model_theory}, the same result is applicable for $ \mbb{R} $.
\end{proof}

We present the first values of $ \alpha_n $ and $ \beta_n $ in Table \ref{table:first_real_coount_coefs}.
As for the asymptotic behavior, computation of the first 10000 values of $ \alpha $ suggests that $ \alpha_n $ tends to a constant (approximately 4.228) as $ n\to\infty $.
Since the choice of the signs of $ \xi_{j,i} $ is probably not optimal, the exact value of the limit is of little interest.
Instead, we prove a lower bound in Corollary \ref{cor:at_least_4dn_sols} and an upper bound in Remark \ref{remark:alpha_upper_bound}.

\begin{cor}\label{cor:at_least_4dn_sols}
	For $ n\ge 14 $ and $ d\gg 0 $ there exists a generic real pencil of hypersurfaces of degree $ d $ in $ \mbb{P}^n $ that contains at least $ 4d^n $ singular hypersurfaces.
\end{cor}

\begin{proof}
	Computing the first 100 values of $ \alpha_n $ we verify that $ \alpha_n>4 $ for all $ 14\le n\le 100 $.
	Denote by $ \widetilde{\alpha}_n:=\frac{1}{n}\left( \sum_{k=0}^{n-1}\alpha_k \right )  $ the sequence of averages.
	A direct computation shows that $ \widetilde{\alpha}_{100}>4 $.
	We can now show, by induction on $ n $, that both $ \alpha_n $ and $ \widetilde{\alpha}_n $ are greater than $ 4 $ for all $ n\ge 100 $.
	Indeed, for $ n\ge 100 $,
	\begin{equation*}
	\alpha_n = \frac{1}{n}\left( \alpha_{n-1}+\sum_{r=0}^{n-2}(1+\beta_{n-1-r})\alpha_r \right)  \ge \frac{1}{n}\sum_{r=0}^{n-1}\alpha_r=\widetilde{\alpha}_n>4.
	\end{equation*}
	and so,
	\begin{equation*}
	\widetilde{\alpha}_{n+1} = \frac{1}{n+1}\left[ \sum_{k=0}^n \alpha_k  \right]  > \frac{1}{n+1}\left[ 4n+\alpha_n \right] > 4.  
	\end{equation*}
\end{proof}

\begin{table}[h]
	\begin{subtable}[h]{\textwidth}
		\begin{tabular}{| c | c | c | c | c | c | c | c | c | c | c | c | c | c |}
			\hline
			$ n $ & $ 0 $ & $ 1 $ & $ 2 $ & $ 3 $ & $ 4 $ & $ 5 $ & $ 6 $ & $ 7 $ & $ 8 $ & $ 9 $ & $ 10 $ & $ 11 $ & $ 12 $ \\ \hline \hline
			$ \alpha_n $ & $ 1.000 $ & $ 2.000 $ & $ 2.000 $ & $ 2.500 $ & $ 2.750 $ & $ 3.025 $ & $ 3.225 $ & $ 3.402 $ & $ 3.543 $ & $ 3.662 $ & $ 3.760 $ & $ 3.841 $ & $ 3.908 $ \\ \hline 
			$ \beta_n $ & $ 1.000 $ & $ 1.000 $ & $ 0.500 $ & $ 0.500 $ & $ 0.375 $ & $ 0.325 $ & $ 0.262 $ & $ 0.220 $ & $ 0.181 $ & $ 0.150 $ & $ 0.124 $ & $ 0.102 $ & $ 0.085 $ \\ \hline 
			$ \gamma_n $ & $ 1.000 $ & $ 1.000 $ & $ 1.500 $ & $ 1.667 $ & $ 1.917 $ & $ 2.075 $ & $ 2.226 $ & $ 2.343 $ & $ 2.443 $ & $ 2.524 $ & $ 2.592 $ & $ 2.648 $ & $ 2.694 $ \\ \hline
		\end{tabular}
	\end{subtable}
	\hfill\\
	\hfill\\
	\hfill\\
	\begin{subtable}[h]{\textwidth}
		\begin{tabular}{| c | c | c | c | c | c | c | c | c | c | c | c | c | c |}
			\hline
			$ n $ & $ 13 $ & $ 14 $ & $ 15 $ & $ 16 $ & $ 17 $ & $ 18 $ & $ 19 $ & $ 20 $ & $ 21 $ & $ 22 $ & $ 23 $ & $ 24 $ & $ 25 $ \\ \hline \hline
			$ \alpha_n $ & $ 3.963 $ & $ 4.009 $ & $ 4.047 $ & $ 4.078 $ & $ 4.104 $ & $ 4.125 $ & $ 4.143 $ & $ 4.158 $ & $ 4.170 $ & $ 4.180 $ & $ 4.188 $ & $ 4.195 $ & $ 4.201 $ \\ \hline 
			$ \beta_n $ & $ 0.070 $ & $ 0.058 $ & $ 0.048 $ & $ 0.040 $ & $ 0.033 $ & $ 0.027 $ & $ 0.022 $ & $ 0.019 $ & $ 0.015 $ & $ 0.013 $ & $ 0.010 $ & $ 0.009 $ & $ 0.007 $ \\ \hline 
			$ \gamma_n $ & $ 2.732 $ & $ 2.764 $ & $ 2.790 $ & $ 2.812 $ & $ 2.830 $ & $ 2.845 $ & $ 2.857 $ & $ 2.867 $ & $ 2.875 $ & $ 2.882 $ & $ 2.888 $ & $ 2.893 $ & $ 2.897 $ \\ \hline
		\end{tabular}
	\end{subtable}
	\hfill\\
	\hfill
	\caption{The first values of $ \alpha_n $, $ \beta_n $ and $ \gamma_n $ as defined in Theorem \ref{thm:real_count_formula} and Theorem \ref{thm:multinode_count}. All values specified up to 3 significant digits.}\label{table:first_real_coount_coefs}
\end{table}

\begin{remark}\label{remark:alpha_upper_bound}
	We can use similar technique to show that $ \alpha_n<5 $ for all $n\ge 0$.
	To that end, we begin by bounding $ \beta_r $ by an exponentially decreasing sequence. 
	The exact parameters of the bound are of little importance\footnote{We can actually show that the generating function of $ \beta $ is $ \sum_{r=1}^{\infty}\beta_rz^r = \frac{\sqrt{3}}{2}\tan(\frac{\sqrt{3}}{2}z+\frac{\pi}{6}) $ which implies that $ \limsup_{r\to \infty}\sqrt[r]{\beta_r}=\frac{3\sqrt{3}}{2\pi}\approx 0.827 $ so this is the "correct" ration of this exponential sequence.}, let us show, by induction on $ r $, that $ \beta_r\le (0.9)^{r+1} $ for all $ r\ge 2 $.
	This inequality holds for $ \beta_2 = \frac{1}{2} \le (0.9)^3 $ and $ \beta_3=\frac{1}{2}\le (0.9)^4 $.
	For $ r\ge 4 $, 
	\begin{equation*}
	\begin{split}
	\beta_r &= \frac{1}{r}\left[ \beta_{r-1} + \sum_{r_0=1}^{r-2} \beta_{r_0} \beta_{r-1-r_0} \right] \le \frac{1}{r}\left[ (0.9)^{r} + 2\beta_1\cdot (0.9)^{r-1} + \sum_{r_0=2}^{r-3} (0.9)^{r+1} \right]  = \\
	&= (0.9)^{r-1}\frac{2.9+(0.9)^2\cdot (r-4)}{r} = (0.9)^{r-1}\left( (0.9)^2 - \frac{0.34}{r} \right)  \le (0.9)^{r+1}.
	\end{split}
	\end{equation*}
	
	Now, we can use this to prove the bound $ \alpha_n< 5 $. 
	Denote by $ \delta\alpha := \sum_{n=0}^{99}(5-\alpha_n) $ the sum of the discrepancies between the bound and the actual value of $ \alpha $ for the first 100 places.
	Computing explicitly the first values of $\alpha$ shows that $ \delta\alpha \approx 95.1 > 60 $ and that $ \alpha_n<5 $ for all $ n<100 $.
	For $ n\ge 100 $,
	\begin{equation*}
	\begin{split}
	\alpha_n &=  \frac{1}{n} \left[ \sum_{r=0}^{n-1}\alpha_r + \sum_{r=1}^{n-1}\beta_r\alpha_{n-1-r}\right] < \frac{1}{n}\left[ 5n-\delta\alpha + 5 + \sum_{r=2}^{n-1}5\cdot (0.9)^{r+1} \right] = \\
	&= 5 + \frac{1}{n}\left[5 +  5\cdot \frac{(0.9)^3-(0.9)^n}{1-0.9} - \delta\alpha \right] \le 5 + \frac{5 + 50\cdot (0.9)^3 - \delta\alpha}{n} \le 5.
	\end{split}
	\end{equation*}
\end{remark}

\subsection{Enumeration of real multi-nodal hypersurfaces}\label{sec:multinodal}

In this subsection we use the results of the current work to count the number of real multi-nodal hypersurfaces satisfying certain point conditions.
It was observed in \cite{Markwig2019_floor_plan} that the tropical approach is particularly well suited for the multi-nodal problem.

We first extend Definition \ref{def:real_multiplicity_with_xi} to the multi-nodal case.
\begin{Def}
	Let $ \delta\in\mbb{N} $. 
	For a choice $ \xi_{j,i}\in\{\pm 1 \} \; (j=1,\dots,N-\delta; i=1,\dots,n) $ and Mikhalkin graded circuits $ C^{(1)},\dots,C^{(\delta)} $, denote by $ mt(C^{(1)},\dots,C^{(\delta)}; (\xi_{j,i})) $ the number of real singular hypersurfaces passing through the points $ \xi_j t^{-M_j\bphi}\in (\mbb{K}_\mbb{R}^\times)^n $ and having a singularities tropicalizing to the tropical singular points corresponding to $ C^{(1)},\dots,C^{(\delta)} $.
\end{Def}

It is a non trivial task to extend the choice of signs made in the previous subsection to the multi-nodal case, namely picking signs that ensure the equations corresponding to the first  levels of all the graded circuits have two real solutions.
Instead, we focus in the current work on the calculation of the average number of $ \delta $-nodal hypersurfaces satisfying point conditions.

\begin{lemma}
	Let $ 0\le k_1<\dots<k_\delta\le d  $, and let $ C^{(1)},\dots,C^{(\delta)}\in \mgc(\Delta) $ with $ C^{(i)} $ centered at a point $ w\in\Delta\cap\mbb{Z}^{n+1} $ with $ w_n=k_i $.
	Then the average of $ mt(C^{(1)},\dots, C^{(\delta)}; (\xi_{i,j})) $ over all possible choices of $ (\xi_{i,j})\in \{\pm 1\}^{(N-\delta)\cdot n} $ is 1.
\end{lemma}

\begin{proof}
	The condition that the $ n $'th coordinate of the center of each graded circuit differs, ensures that for each $ j=1,\dots,N-\delta $, $ \xi_{j,i} $ affects the equations of at most one graded circuit.
	We then can mimic the proof of Lemma \ref{lemma:real_count_average}. 
	In particular, we can show by induction on $ \delta_0 $ and $ h $ that the average number of real solutions to the equations corresponding to $ C^{(1)},\dots,C^{(\delta_0-1)} $ and to the first $ h $ levels of $ C^{(\delta_0)} $ is 1.
	If $ m\left( C^{(\delta_0)}_h \right) =1  $, the induction step is obvious.
	
	Otherwise, if $ 1\le j\le N-\delta $ is the last s.t. $ \xi_j $ affects the equations corresponding to $ C^{(\delta_0)}_h $, half of the choices of $ \xi_j $ give two real solutions to those equations, and half of the choices give no solutions.
	By the above considerations, $ \xi_j $ does not affect the number of real solutions corresponding to $ C^{(1)},\dots,C^{(\delta_0-1)} $. 
	If we are in the situation of Lemma \ref{lemma:parrallel_projection}, $ \xi_j $ can not appear in the equations corresponding to lower levels of $ C^{(\delta_0)} $ and does not affect the number of real solutions of those equations.
	Suppose now that we are in the situation of Lemma \ref{lemma:gluing_type_IV} and $ C^{(\delta_0)}_h $ is of type I.
	This means that $ C^{(\delta_0)}_{h_0}=:\{\bs{u}^{(-2)},\bs{u}^{(-1)},\bs{u}^{(1)}\} $ is of type IV and index $ r $, and all the points of $ C^{(\delta_0)}_{h_0+1},\dots,C^{(\delta_0)}_{h} $ have $ (r+1) $'th coordinate equal to $ \bs{w}_{r+1}-1 $.
	Then changing $ \xi_j $ change only the coefficients $ a_{\bs{u}^{(1)}},a_{\bs{w}} $ and does not affect the coefficients $ a_{\bs{v}} $ for $ \bs{v}\in \bigcup_{h_0<i<h}C^{(\delta_0)}_i $.
	Moreover, if $ \bs{q} $ is the corresponding singular point of the resulting algebraic hypersurface, changing $ \xi_j $ only changes $ \bs{q}_{r+1} $ which does not affect the number of real solutions to equations corresponding $ C^{(\delta_0)}_{h_0+1},\dots ,C^{(\delta_0)}_{h-1} $.
\end{proof}

Similarly to the uni-nodal case, we can use this average get a count of real singular hypersurfaces with $ \delta $ nodes satisfying point conditions, i.e. to prove Theorem \ref{thm:multinode_count}.

\begin{proof}[Proof of Theorem \ref{thm:multinode_count}]
	By the proof of Theorem \ref{thm:real_count_formula}, there exist $ \gamma_n d^n+O(d^{n-1}) $ graded circuits in $ \Delta_d^{(n)} $ admitting Mikhalkin condition.
	For $ 0 \le a \le d $, denote by $ K^{(n)}_{d,a} $ the number of graded circuits in $ \Delta^{(n)}_d $ admitting Mikhalkin condition and centered at a point $ \bs{w} $ with $ \bs{w}_n=d-a $.
	Let $ C\subseteq \Delta_d^{(n)} $ centered at $ \bs{w} $ with $ \bs{w}_n=d-a $ and suppose that $ C=\widetilde{C}\Vert C' $ with $ C' $ elementary Mikhalkin graded circuit, see Definition \ref{def:elementary_Mikhalkin_graded_circuit}.
	As usual, separate into cases according to the type of $ C'_1 $.
	
	If $ C'_1 $ is of type I, then $ \widetilde{C}\in\mgc(\Delta^{(n-1)}_{a}) $ and for every such graded circuit, $ \widetilde{C}\Vert C' $ contribute to the count of $ K^{(n)}_{d,a} $.
	As in Proposition \ref{prop:discriminant_degree}, type II circuits do not contribute to the asymptotic of the count.
	The contribution to $ K^{(n)}_{d,a} $ of graded circuits with $ C'_1 $ of type III and lower index $ r $, is $ \gamma_{n-1-r} a^{n-1} + O(a^{n-2}) $.
	For $ C'_1 $ of type IV and $ \dim C'=r $, we get $ \gamma_{n-r}\beta_{r-1} a^{n-1}+O(a^{n-2}) $ such circuits.
	Combining those together we get that
	\begin{equation*}
	K^{(n)}_{d,a} = \left( \gamma_{n-1} + \sum_{r=1}^{n-1}\gamma_{n-1-r} + \sum_{r=0}^{n-2}\gamma_{n-r}\beta_{r-1} \right)a^{n-1} + O(a^{n-2}) = n\gamma_n a^{n-1}+O(a^{n-2})
	\end{equation*}
	
	We conclude that, for fixed $ 0<  k_1 < k_2 <\dots < k_\delta < d  $, the number of $ \delta $-tuples of graded circuits $ C^{(1)},\dots,C^{(\delta)} $ s.t. $ C_i $ is centered at a point with $ n $'th coordinate equal to $ k_i $ is 
	\begin{equation*}
	\prod_{i=1}^\delta \left( n\gamma_n k_i^{n-1} + O(k_i^{n-2}) \right) = n^\delta\gamma_n^\delta \left (\prod_{i=1}^\delta k_i\right )^{n-1} + O\left( \left (\prod_{i=1}^\delta k_i\right )^{n-2} \right).
	\end{equation*}
	Summing over all choices of coordinates $ 0<k_1<\dots<k_\delta<d $ and using Faulhaber's formula $ \sum_{a=1}^d a^p = \frac{1}{p+1}d^{p+1}+O(d^p) $, $\delta$ times we get
	\begin{align*}
	\sum_{0<k_1<\dots<k_\delta<d}\prod_{i=1}^\delta K^{(n)}_{d,k_i} = \sum_{k_\delta=1}^{d}\sum_{k_{\delta-1}=1}^{k_\delta}\dots \sum_{k_2=1}^{k_3} \sum_{k_1=1}^{k_2} \left [ n^\delta\gamma_n^\delta \left (\prod_{i=1}^\delta k_i\right )^{n-1} + O\left( \left (\prod_{i=1}^\delta k_i\right )^{n-2} \right) \right ] &= \\
	= n^\delta \gamma_n^\delta\cdot (1+o(1)) \sum_{k_\delta=1}^{d} k_\delta^{n-1}\sum_{k_{\delta-1}=1}^{k_\delta} k_{\delta-1}^{n-1}\dots \sum_{k_2=1}^{k_3}k_2^{n-1}\sum_{k_1=1}^{k_2}k_1^{n-1} &=  \\
	= n^\delta \gamma_n^\delta\cdot (1+o(1)) \sum_{k_\delta=1}^{d} k_\delta^{n-1}\sum_{k_{\delta-1}=1}^{k_\delta} k_{\delta-1}^{n-1}\dots \sum_{k_2=1}^{k_3}\frac{1}{n}k_2^{2n-1} &= \\
	= n^\delta \gamma_n^\delta\cdot (1+o(1)) \sum_{k_\delta=1}^{d} k_\delta^{n-1}\sum_{k_{\delta-1}=1}^{k_\delta} k_{\delta-1}^{n-1}\dots \sum_{k_{i+1}=1}^{k_{i+2}}\frac{1}{i!\cdot n^i}k_2^{i\cdot n-1} 
	= \frac{1}{\delta!}\gamma_n^\delta d^{\delta\cdot n}+O(d^{\delta\cdot n-1}) &
	\end{align*}
	
	Which proves the desired count for the field $ \mbb{K_R} $ of real Puiseux series.
	The result over $ \mbb{R} $ follows from Tarski's transfer principle.
\end{proof}

\begin{remark}
	We neglected in Theorem \ref{thm:multinode_count} all $ \delta $-tuples of graded circuits with repeating $ n $'th coordinate of $ \bs{w} $ and graded circuits centered at $ bs{w}
	 $ with $ \bs{w}_n=0 $.
	As can be seen from Proposition \ref{prop:discriminant_degree} (or by a direct computation), they do not contribute to the asymptotic of the count.
\end{remark}

The first values of $ \gamma_n $ are shown in Table \ref{table:first_real_coount_coefs}. 
Computation of the first 10000 values suggests that $ \gamma_n $ tends to a constant (approximately 2.915) as $ n\to \infty $.
Similarly to our analysis of $ \alpha_n $ we will prove lower and upper bounds of $ \gamma_n $.

\begin{cor}\label{cor:gamma_lower_bound}
	Let $ \delta\in\mbb{N} $. 
	For all $ n\ge 9 $ and $ d\gg 0 $ there exists configuration $ \overline{\boldsymbol{w}} $ of $ N-\delta = \binom{n+d}{n}-\delta $ real points in $ \mbb{P}^n $ s.t. there are at least $ \frac{1}{\delta!}\left( \frac{5}{2} d^n \right)^\delta + O(d^{n\delta-1}) $ real hypersurfaces passing through $ \overline{\boldsymbol{w}} $ and having $ \delta $ real nodes.
\end{cor}

\begin{proof}
	The proof is identical to the proof of Corollary \ref{cor:at_least_4dn_sols}.
	Denote by $ \widetilde{\gamma}_n:=\frac{1}{n}\sum_{k=0}^{n-1}\gamma_k $, the sequence of averages. 
	A direct computation shows that $ \gamma_n>\frac{5}{2} $ for all $ 9\le n \le 50 $ and that $ \widetilde{\gamma}_{50} > \frac{5}{2} $.
	Then, the induction proof of Corollary \ref{cor:at_least_4dn_sols} is applicable here, since the recurrence relation of $ \gamma $ is the same as of $ \alpha $ (and only the initial conditions differ).
\end{proof}

\begin{remark}
	Similarly to Remark \ref{remark:alpha_upper_bound}, we can show that $ \gamma_n< 3 $ for all $ n $.
	Indeed, computing explicitly the first values of $ \gamma $ we see that the bounds holds for $ n<300 $ and the discrepancy up to this place is  $ \delta\gamma:=\sum_{n=0}^{299}(3-\gamma_n)\approx37.76 > 33 $.
	Thus, for $ n\ge 300 $,	
	\begin{equation*}
	\begin{split}
	\gamma_n &= \frac{1}{n}\left[ \sum_{r=0}^{n-1}\gamma_r + \sum_{r=1}^{n-1}\beta_r \gamma_{n-1-r} \right] < \frac{1}{n}\left[ 3n-\delta\gamma + 3+\sum_{r=2}^{n-1}3\cdot (0.9)^{r+1} \right] = \\
	&= 3+\frac{1}{n}\left[ 3 + 3\cdot \frac{(0.9)^3-(0.9)^n}{1-0.9}-\delta\gamma \right]  \le 3+\frac{3+30\cdot (0.9)^3-\delta\gamma}{n} < 3.
	\end{split}
	\end{equation*}
\end{remark}

\insertbibliography

\todos


\providecommand{\bysame}{\leavevmode\hbox to3em{\hrulefill}\thinspace}
\providecommand{\MR}{\relax\ifhmode\unskip\space\fi MR }
\providecommand{\MRhref}[2]{%
  \href{http://www.ams.org/mathscinet-getitem?mr=#1}{#2}
}
\providecommand{\href}[2]{#2}
\begin{thebibliography}{10}

\bibitem{Dickenstein2012}
A.~Dickenstein and L.~F. Tabera, \emph{{Singular Tropical Hypersurfaces}},
  Discrete and Computational Geometry \textbf{47} (2012), no.~2, 430--453.

\bibitem{Gelfand_kapranov}
I.~M. Gelfand, M.~M. Kapranov, and A.~V. Zelevinsky, \emph{{Discriminants,
  Resultants and Multidimensional Determinants}}, Modern Birkhauser Classics,
  1994.

\bibitem{Grunbaum2003}
B.~Grunbaum, \emph{{Convex Polytopes}}, Graduate Texts in Mathematics,
  Springer, 2003.

\bibitem{Santos_Valino_class_dim4_simpl}
{\'{O}}.~Iglesias~Vali{\~{n}}o and F.~Santos, \emph{{The complete
  classification of empty lattice 4-simplices}}, Electronic Notes in Discrete
  Mathematics \textbf{68} (2018), 155--160.

\bibitem{Mikhalkin_Shustin_Ittenberg}
I.~Itenberg, G.~Mikhalkin, and E.~Shustin, \emph{Tropical algebraic geometry},
  Oberwolfach Seminars, Birkh{\"a}user Basel, 2009.

\bibitem{Jensen_Lenzing_model_theory}
C.~U. Jensen and H.~Lenzing, \emph{Model-theoretic algebra with particular
  emphasis on fields, rings, modules}, Algebra, Logic and Applications, vol.~2,
  Gordon and Breach Science Publishers, New York, Switzerland, 1989 (English).

\bibitem{Kerner_complex_multinodal}
D.~Kerner, personal communication, 2017.

\bibitem{Maclagan2014}
D.~Maclagan and B.~Sturmfels, \emph{{Introduction to Tropical Geometry}},
  Graduate Studies in Mathematics, American Mathematical Society, 2015.

\bibitem{Markwig2019_floor_plan}
H.~Markwig, T.~Markwig, K.~Shaw, and E.~Shustin, \emph{{Tropical floor plans
  and enumeration of complex and real multi-nodal surfaces}}, arXiv preprint
  arXiv:1910.08585 (2019).

\bibitem{Markwig2012_trop_curve_sing}
H.~Markwig, T.~Markwig, and E.~Shustin, \emph{{Tropical curves with a
  singularity in a fixed point}}, Manuscripta Mathematica \textbf{137} (2012),
  no.~3, 383--418.

\bibitem{Markwig2012_trop_surf_sings}
\bysame, \emph{{Tropical Surface Singularities}}, Discrete and Computational
  Geometry \textbf{48} (2012), no.~4, 879--914.

\bibitem{Shustin2017_enum}
\bysame, \emph{{Enumeration of complex and real surfaces via tropical
  geometry}}, Advances in Geometry \textbf{18} (2018), no.~1, 69--100.

\bibitem{Mikhalkin2003_lattice_path}
G.~Mikhalkin, \emph{{Counting curves via lattice paths in polygons}}, Comptes
  Rendus Mathematique \textbf{336} (2003), no.~8, 629--634.

\bibitem{Mikhalkin2003}
\bysame, \emph{{Enumerative tropical algebraic geometry in $ \mbb{R}^2 $}},
  Journal of the American Mathematical Society \textbf{18} (2005), no.~2,
  313--377.

\bibitem{Shustin2005}
E.~Shustin, \emph{{A Tropical Approach to Enumerative Geometry}}, Algebra i
  Analiz \textbf{17} (2005), 170--214.

\bibitem{Viro1876}
O.~Viro, \emph{{Patchworking Real Algebraic Varieties}}, {Appendix to G.-M.
  Greuel et al. Singular Algebraic Curves, Springer, Switzerland, 2018}.

\bibitem{Ziegler2007}
G.~M. Ziegler, \emph{{Lectures On Polytopes}}, Graduate Texts in Mathematics,
  Springer, 2007.

\end{thebibliography}
\end{document}